\documentclass[12pt]{article}
\usepackage[english]{babel}
\usepackage{sfmath}

\usepackage{amssymb}
\usepackage{natbib}
\usepackage[leqno]{mathtools}

\usepackage{amsthm}

\usepackage{xr-hyper}
\usepackage[bookmarks,pdfnewwindow,plainpages=false]{hyperref}

\hypersetup{
colorlinks=true,
linkcolor=black,
citecolor=black,
pdfauthor={Victor Ivrii},
pdftitle={Sharp Spectral Asymptotics for four-dimensional Schr\"odinger operator with a strong but degenerating magnetic field. II},
pdfsubject={Sharp Spectral Asymptotics},
pdfkeywords={Microlocal Analysis, Magnetic Schr\"odinger Operator},
}

\usepackage{graphicx,color}

\definecolor{lightgray}{gray}{.90}
\definecolor{verylightgray}{gray}{.93}

\usepackage{subfigure}

\usepackage{enumitem}

\newcommand\W{{\rm W}}
\newcommand\MW{{\rm {MW}}}
\newcommand\Def{{\overset {\rm {def}}{\ =\ }}}
\newcommand\redef{{\overset {\rm {redef}}{\ =\ }}}

\newcommand\out{{\rm {out}}}
\newcommand\inn{{\rm {inn}}}
\newcommand\per{{\rm {per}}}
\newcommand\const{{\rm {const}}}
\newcommand\dist{{\rm {dist}}}

\newcommand\corr{{\rm {corr}}}

\newcommand\w{{\rm w}}

\newcommand\bx{{\bar {\mathstrut x}}}

\newcommand\bz{{\bar {\mathstrut z}}}

\newcommand\bxi{\bar {\mathstrut \xi}}

\newcommand\Span{\operatorname{Span}}
\newcommand\codim{\operatorname{codim}}

\newcommand\Hess{\operatorname{Hess}}

\renewcommand\Im{\operatorname{Im}}

\newcommand\Ker{\operatorname{Ker}}

\newcommand\rank{\operatorname{rank}}
\renewcommand\Re{\operatorname{Re}}

\newcommand\sing{{\operatorname{sing}}}

\newcommand\bR{{\mathbb R}}

\newcommand\bK{{\mathbb K}}

\newcommand\bZ{{\mathbb Z}}
\newcommand\cA{{\mathcal A}}
\newcommand\cB{{\mathcal B}}

\newcommand\cE{{\mathcal E}}
\newcommand\cF{{\mathcal F}}

\newcommand\cL{{\mathcal L}}

\newcommand\cO{{\mathcal O}}

\newcommand\cQ{{\mathcal Q}}

\newcommand\cY{{\mathcal Y}}
\newcommand\cZ{{\mathcal Z}}

\newtheorem{theorem}{Theorem}[section]

\newtheorem{corollary}[theorem]{Corollary}
\newtheorem{proposition}[theorem]{Proposition}
\newtheorem{lemma}[theorem]{Lemma}

\theoremstyle{definition}

\theoremstyle{remark}
\newtheorem{remark}[theorem]{Remark}
\newtheorem{example}[theorem]{Example}

\numberwithin{equation}{section}

\newenvironment{claim}[1][{\rm(\theequation)}]{\refstepcounter{equation}\vglue10pt
\begin{trivlist}
\item[{\hskip\labelsep#1}]}{\vglue10pt\end{trivlist}}

\newcounter{note}

\newenvironment{claim*}[1]{\smallskip
\begin{trivlist}
\item[{\hskip\labelsep#1}]}{\smallskip\end{trivlist}}

\newenvironment{phantomequation}[1][]{\refstepcounter{equation}}{}

\renewcommand\subsubsection{\paragraph{\thesubsubsection}\refstepcounter{subsubsection}}

\begin{document}



\title{%
Sharp spectral asymptotics for four-dimensional Schr\"odinger operator with a strong degenerating magnetic field.}

\author{%
Victor Ivrii
\footnote{Work was partially supported by NSERC grant OGP0138277.}
}
\maketitle


{\abstract%
I continue analysis of the Schr\"odinger operator with the strong degenerating magnetic field, started in \cite{IRO6}. Now I consider 4-dimensional case, assuming that magnetic field is generic degenerated and under certain conditions
I derive spectral asymptotics with the principal part $\asymp h^{-4}$ and the remainder estimate $O(\mu^{-1/2}h^{-3})$ where $\mu\gg 1$ is the intensity of the field and $h\ll 1$ is the Plank constant; $\mu h\le 1$.

These asymptotics can contain correction terms of magnitude
$\mu ^{5/4}h^{-3/2}$ corresponding to the short periodic trajectories.
\endabstract}

\section{Introduction}
\setcounter{subsection}{-1}
\label{sect-0}

\subsection{Preface}
\label{sect-0-0}

I continue analysis of he Schr\"odinger operator with the strong degenerating magnetic field, started in \cite{IRO6}
\begin{equation}
A= {\frac 1 2}\Bigl(\sum_{j,k}P_jg^{jk}(x)P_k -V\Bigr),\qquad P_j=D_j-\mu V_j
\label{0-1}
\end{equation}
where $g^{jk}$, $V_j$, $V$ are smooth real-valued functions of $x\in \bR^2$ and
$(g^{jk})$ is positive-definite matrix, $0<h\ll 1$ is a Planck parameter and
$\mu \gg1$ is a coupling parameter. I assume that $A$ is a self-adjoint operator.

Now I consider 4-dimensional case, assuming that magnetic field is generic degenerated, i.e. of Martinet-Roussarie type \cite{Ma,Rou} which will be described in details in subsection \ref{sect-1-1}. Here I just mention that in this case magnetic field $F$ (which currently is considered as a closed 2-form) degenerates on the manifold $\Sigma$ of dimension 3 and $\dim\Ker F =2$ at $\Sigma$. Further, $\dim (\Ker F \cap T\Sigma)=1$ at $\Sigma\setminus \Lambda$ and $\Ker F \subset T\Sigma$ at $\Lambda$ where $\Lambda\subset \Sigma$ is 1-dimensional manifold; furthermore, an angle between $\Ker F$ and $T\Sigma$ is exactly of magnitude $\dist (x,\Lambda)$. Finally, in an appropriate coordinates $\Sigma=\{x_1=0\}$, $\Lambda=\{x_1=x_3=x_4=0\}$ and \emph{magnetic lines} ${\frac {dx} {dt}}\in \Ker F\cap T\Sigma$ are circles $\{x_1=0, x_2=\const, x_3^2+x_4^2=\const\}$.

My goal is to find asymptotics of $\int e(x,x,0)\psi (x)\, dx$ with respect to $h,\mu$, where $e(x,y,\tau)$ is the Schwartz kernel of the spectral projector of $A$ and $\psi$ is a smooth function supported in the vicinity of $\Sigma$. I assume that $\mu h\le \const$ (otherwise $e(x,y,0)= O(\mu^{-\infty})$).

In the next paper I am gong to analyze the case of the generic non-degenerated field.

\subsection{Assumptions, notations and results}
\label{sect-0-1}
So, magnetic field is characterized by $F_{jk}=\partial_jV_k-\partial_kV_j$.

However, from the point of view of the classical and quantum dynamics and spectral asymptotics properties of $(F^j_k)=(g^{jl})(F_{lk})$ are more important than those of $(F_{jk})$. Let
$\pm if_1,\pm if_2$ be eigenvalues of $(F^j_k)$, $f_j\ge 0$. Then with the correct choice of notations $f_1 \asymp \dist (x,\Sigma)$ and $f_2\asymp 1$.

My first statement holds almost without any further assumptions:

\begin{theorem}\label{thm-0-1} Let $F$ be of Martinet-Roussarie type and
\begin{equation}
V\ge \epsilon_0>0.
\label{0-2}
\end{equation}
Let $\psi$ be supported in the small vicinity of $\Sigma$. Then
\begin{equation}
|\int \bigl( e (x,x,0)-\cE^\MW (x,0)\bigr)\psi (x)\,dx| \le C\mu^{-1/2}h^{-3}+ C\mu^2 h^{-2}
\label{0-3}
\end{equation}
where
\begin{multline}
\cE^\MW (x,\tau)= \\(2\pi)^{-2}\mu^2h^{-2}\sum _{(m,n)\in \bZ^{+\,2} } \theta \bigl(2\tau+V- (2m+1)\mu h f_1 - (2n+1)\mu hf_2\bigr) f_1 f_2 \sqrt g
\label{0-4}
\end{multline}
is Magnetic Weyl Expression (see \cite{IRO3} for the general case),
$g=\det (g^{jk})^{-1}$.
\end{theorem}

\begin{remark}\label{rem-0-2} Without any additional conditions remainder estimate (\ref{0-3}) is sharp. An extra term $O(\mu^2h^{-2})$ in its right-hand expression is due to the possibility that for some pair $(m,n)\in \bZ^{+\,2}$
function $V-(2m+1)\mu h f_1 - (2n+1)\mu hf_2$ identically vanishes in the vicinity of some point of $\Sigma$. However it is not as bad as it could be in the case of constant comeasurable $f_1,f_2$ in which case
$V-(2m+1)\mu h f_1 - (2n+1)\mu hf_2$ could vanish up to $\asymp (\mu h)^{-1}$ pairs $(m,n)\in \bZ^{+\,2}$
and the remainder estimate would be the worst possible $O(\mu h^{-3})$.
\end{remark}

In the next few statements I am improving remainder estimate (\ref{0-3}).

\begin{theorem}\label{thm-0-3} Let $F$ be of Martinet-Roussarie type, conditions $(\ref{0-2})$ and
\begin{equation}
|\bigl({\frac V {f_2}}\bigr)- (2n+1)\mu h| + |\nabla_\Sigma \bigl({\frac V {f_2}}\bigr)|\le \epsilon_0 \implies |\det \Hess_\Sigma \bigl({\frac V {f_2}}\bigr) |\ge \epsilon_0\qquad\forall n\in \bZ^+
\label{0-5}
\end{equation}
be fulfilled. Let $\psi$ be supported in the small vicinity of $\Sigma$. Then
\begin{equation}
|\int \bigl( e (x,x,0)-\cE^\MW (x,0)\bigr)\psi (x)\,dx -\int_\Sigma \cE^\MW_\corr (x',\tau)\psi (x')\,d x'| \le C\mu^{-1/2}h^{-3}
\label{0-6}
\end{equation}
where $dx'$ is the standard density on $\Sigma$ and $\cE^\MW_\corr (x',\tau)$ is a correction term (see $(\ref{4-51})$ which is $O(\mu^{5/4}h^{-3/2})$.
\end{theorem}

\begin{remark}\label{rem-0-4} (i) Correction terms $\cE^\MW_\corr$ and $\cE^\MW_\corr$ are of completely different nature;

\smallskip\noindent
(ii) The first term in the left-hand expression of (\ref{0-5}) is meaningful only in the case of the \emph{very strong magnetic field} $\mu h\ge \epsilon$; so as $\mu h\le \epsilon$ condition (\ref{0-5}) is equivalent to
\begin{equation}
|\nabla_\Sigma \bigl({\frac V {f_2}}\bigr)|\le \epsilon_0 \implies |\det \Hess_\Sigma \bigl({\frac V {f_2}}\bigr) |\ge \epsilon_0;
\label{0-7}
\end{equation}

\smallskip\noindent
(iii) One can weaken condition (\ref{0-5}) to
\begin{phantomequation}\label{0-8}\end{phantomequation}
\begin{multline}
\exists n\in \bZ^+:\ |{\frac V {f_2}} -(2n+1)f_2\mu h|+ |\nabla_\Sigma \bigl({\frac V {f_2}} \bigr)|\le \epsilon_0 \implies\\ \Hess_\Sigma \bigl({\frac V {f_2}} \bigr)\, \text{has at least $q$ eigenvalues with absolute values greater than $\epsilon$}_0
\tag*{$(0.8)_q$}\label{0-8-q}
\end{multline}
with $q=1,2,3$ depending on the magnitude of $\mu$; as $\mu h\le \epsilon$ it is equivalent to
\begin{phantomequation}\label{0-9}\end{phantomequation}
\begin{multline}
|\nabla_\Sigma \bigl({\frac V {f_2}} \bigr)|\le \epsilon_0 \implies\\ \Hess_\Sigma \bigl({\frac V {f_2}} \bigr)\, \text{has at least $q$ eigenvalues with absolute values greater than $\epsilon$}_0
\tag*{$(0.9)_q$}\label{0-9-q};
\end{multline}
obviously (\ref{0-5}), (\ref{0-7}) are equivalent to $(\ref{0-8})_3$, $(\ref{0-9})_3$ respectively;

\smallskip\noindent
(iv) Under stronger condition
\begin{equation}
|\bigl({\frac V {f_2}}\bigr)- (2n+1)\mu h| + |\nabla_\Sigma \bigl({\frac V {f_2}} \bigr)|\ge \epsilon_0\qquad\forall n\in \bZ^+
\label{0-10}
\end{equation}
the proof of theorem \ref{thm-0-3} is much easier; as $\mu h\ll 1$ this condition is equivalent to
\begin{equation}
|\nabla_\Sigma \bigl({\frac V {f_2}} \bigr)|\ge \epsilon_0.
\label{0-11}
\end{equation}

\smallskip\noindent
(v) In the very special cases (see Appendix \ref{sect-A}) one can derive remainder estimate (\ref{0-6}) without condition (\ref{0-5}) or even $(\ref{0-8})_1$. Then the correction term could be larger, up to
$O(\mu^{1/4}h^{-5/2})$.
\end{remark}

I will also prove

\begin{theorem}\label{thm-0-5} Let $F$ is of Martinet-Roussarie type and condition $(\ref{0-10})$ be fulfilled. Then as $\psi$ is supported in the small vicinity of $\Sigma\cap \{V=0\}$ asymptotics $(\ref{0-6})$ holds with
$ \cE^\MW_\corr (x',\tau)=0$.
\end{theorem}

\subsection{Plan of the paper}
\label{sect-0-2}

First of all, in section \ref{sect-1} I study the weak magnetic case $\mu \le h^{-\delta}$ with sufficiently small exponent $\delta>0$ and derive remainder estimate $O(\mu^{-1/2}h^{-3})$ as the main part is given by the standard Weyl formula. I also analyze there the geometry of the degenerate magnetic field and the corresponding classical dynamics.

Then, assuming that magnetic field is not weak, $h^{-\delta}\le \mu \le ch^{-1}$ in section \ref{sect-2} I consider the different canonical forms of the operator in question; this canonical forms contain powers of $\mu^{-1}$. There are universal canonical forms and also specific canonical forms as $\dist (x,\Sigma)\gg \mu^{-1/2}$ or as $\dist (x,\Lambda)\gg \mu^{-1/2}$.

Then in section \ref{sect-3} I derive remainder estimates $O(\mu^{-1/2}h^{-3}+\mu^2h^{-2})$ in the general case and
$O(\mu^{-1/2}h^{-3})$ as $\mu \ge h^{-2/5}$ and some non-degeneracy condition (depending on the magnitude of $\mu$) is fulfilled. The main part of this estimate is given by the \emph{standard intermediate formula}
\begin{equation}
h^{-1}\int_{-\infty}^0 \Bigl( F_{t\to h^{-1}\tau} {\bar\chi}_T(t) \Gamma uQ_y^t\Bigr)\,d\tau
\label{0-12}
\end{equation}
where $u$ is the Schwartz kernel of $e^{ih^{-1}tA}$, $Q=I$ and ${\bar\chi}_T(t)$ equals 1 as $|t|\le {\frac 1 2}T$ and vanishes as $|t|\ge T$, $T$ is rather arbitrary from interval $[T_0,T_1]$ and the main fight is to make $T_0$ as small and $T_1$ as large as possible. More precisely, the main part of asymptotics is given by the sum of expressions (\ref{0-12}) with $Q=Q_{(\iota)}$ making partition of unity and $T=T_\iota \in [T_{(\iota) 0},T_{(\iota) 1}]$.

After this in section \ref{sect-4} I calculate (\ref{0-12}) in more explicit way.

Finally, Appendices \ref{sect-A} are devoted to the analysis of some special cases.

\section{Weak magnetic field}
\label{sect-1}
\subsection{Geometry of degenerating magnetic field}
\label{sect-1-1}
So let $F_{jk}(x)=\partial_j V_k-\partial_k V_j$ be components of the matrix intensity of magnetic field; then
$\omega = \sum_{j,k}F_{jk}dx_j\wedge dx_k=d (\sum_k V_kdx_k)$ is the corresponding magnetic 2-form. According to \cite{Ma}, the generic form in dimension 4 never degenerates completely and has rank 2 on the submanifold $\Sigma$ of codimension 1; more precisely, let
\begin{equation}
\Sigma \Def \bigl\{x: \rank \omega (x)=2\bigr\};
\label{1-1}
\end{equation}
then

\begin{claim}\label{1-2}
For a generic 2-form $\omega$ \ $\Sigma$ is the smooth manifold of codimension 1 and if $\pm if_1$, $\pm if_2$ are eigenvalues of the corresponding matrix $(F_{jk})$, $f_2\ge f_1\ge 0$\,\footnote{\label{foot-1} This is a temporary notation, corresponding to Euclidean metrics $g^{jk}$.} then $f_2(x)\ge \epsilon$ and $f_1(x)\asymp \dist (x,\Sigma)$.
\end{claim}

Let us consider $\Ker F(x) \cap T_x\Sigma$ in points of $\Sigma$. At each point it can be of dimension 1 or 2; let
\begin{equation}
\Lambda \Def \bigl\{x\in \Sigma: \dim \bigl(\Ker F(x) \cap T_x\Sigma\bigr)=2 \};
\label{1-3}
\end{equation}
according to \cite{Rou}

\begin{claim}\label{1-4}
For a generic 2-form $\omega$ \ $\Lambda$ is the smooth manifold of dimension 1; moreover in the appropriate coordinates
\begin{equation}
\omega = dx_1\wedge dx_2 +x_3dx_2\wedge dx_3 +
d\bigl(x_1x_3+x_2x_4-{\frac 1 3}x_3^3\bigr)\wedge dx_4.
\label{1-5}
\end{equation}
\end{claim}

One can rewrite (\ref{1-5}) as
\begin{equation*}
d\Bigl(\bigl(x_1-{\frac 1 2}(x_3^2+x_4^2)\bigr)d(x_2+{\frac 1 2}x_3x_4) +
(x_3^2+x_4^2)x_1-{\frac 1 4}(x_3^2+x_4^2)^2\Bigr);
\end{equation*}
after transformation $x_2\mapsto x_2- {\frac 1 2}x_3x_4$ I get instead
\begin{multline}
\omega = dx_1\wedge dx_2 - x_4dx_1\wedge dx_3 + x_3dx_1\wedge dx_4 +x_3dx_2\wedge dx_3 +x_4dx_2\wedge dx_4 +\\
2\bigl(x_1 -{\frac 1 2}(x_3^2+x_4^2)\bigr)dx_3\wedge dx_4=\\
d\bigl(x_1-{\frac 1 2}r^2\bigr)\wedge dx_2 +r^2 dx_1\wedge d\theta +
2\bigl(x_1-{\frac 1 2}r^2\bigr)rdr\wedge d\theta =
d \Bigl( (x_1-{\frac 1 2}r^2)dx_2 + (x_1-{\frac 1 4}r^2)r^2 d\theta\Bigr)
\label{1-6}
\end{multline}
where $x_3=r\cos \theta$, $x_4=r\sin \theta$. In contrast to the original Roussarie 2-form (\ref{1-5}) this latter form is obviously $(x_3,x_4)$-circular symmetric.

Then \begin{equation}
\Sigma=\{x_1=0\},\quad \Lambda=\{x_1=x_3=x_4=0\}
\label{1-7}
\end{equation}
and \emph{magnetic lines} defined by
\begin{equation}
{\frac {dx}{dt}}\in \Ker F(x)\cap T_x\Sigma, \quad x\in \Sigma
\label{1-8}
\end{equation}
are helices $\{x_1=0, x_3=r\cos \theta, x_4=r\sin \theta, x_2= \const - r^2\theta/2 \}$ (with $r=\const$) winging around $\Lambda$.

Further, away from $\Lambda$ one can rewrite
\begin{equation}
\omega = 2x_1dx_1\wedge d\theta + dy\wedge dz, \qquad y=x_1-{\frac 1 2}r^2,\quad
z=x_2-2y\theta
\label{1-9}
\end{equation}
and therefore there I have just a direct sum of two 2-dimensional magnetic fields; the first one is a generic degenerating field and second one is nondegenerate.

\subsection{Classical dynamics}
\label{sect-1-2}

Classical dynamics described by the Hamiltonian
\begin{equation}
a(x,\xi)= {\frac 1 2}\Bigl(\sum_{j,k}g^{jk}(x)p_jp_k -V\Bigr),\qquad p_j=\xi_j-\mu V_j
\label{1-10}
\end{equation}
corresponding to (\ref{0-1}), is much more complicated than the geometry because it also depends on the metrics $(g^{jk})$ and the scalar potential $V$; however it appears that it depends mainly on $f_1,f_2$ which are eigenvalues of $-i(F^j_k)$\,\footnote{\label{foot-2} In contrast to footnote \ref{foot-1} this is the final definition.} (where $(F^j_k)=(g^{jl})(F_{lk})$) and $V$ and even on their ratios rather than these functions themselves. It also depends on the partition of energy (see below).

\subsubsection{}
\label{sect-1-2-1} Let us start from considering few examples.

\begin{example}\label{ex-1-1} As $F_{jk}$, $g^{jk}$ and $V$ are constant (which is not the case in our article) $\bR^4 =\bK_1\oplus \bK_2$ where $\bK_j$ here and below are eigenspaces of $-(F^j_k)^2$ corresponding to eigenvalues $f_j^2$,
$\dim \bK_j=2$. Then the kinetic part of the Hamiltonian splits into
$a^0= a^0_1 +a^0_2$ where $a^0_j$ are Hamiltonians on $\bK_j$ and $a^0_j$ are movement integrals and the classical dynamics splits into two \emph{cyclotron} movements along elliptical orbits in $\bK_j$ with the angular velocities $\mu f_j$; the sizes of these ellipses are $\asymp (\mu f_j)^{-1}(a_j^0)^{1/2}$.
\end{example}

\begin{example}\label{ex-1-2} (i) Let in the frames of the previous example $V$ be variable. Then $a_j$ are no more movement integrals but they evolve with the average speeds\footnote{\label{foot-3} While the speed of spatial movement is $\asymp 1$ the drift is much slower and one can calculate all functions of $x$ in the instant cyclotron centers.} $O(\mu^{-1})$ and in addition to the fast cyclotron movements there appear also slow \emph{drift} movements along $\bK_j$ with the speeds of $\asymp \mu^{-1}a_j^0 |\nabla_{\bK_j} (V/f_j)|$; these drifts are orthogonal to $\nabla_{\bK_j} (V/f_j)$. So these drifts depend on the instant energy partition.

\smallskip
\noindent
(ii) As $(g^{jk})$ and $(F_{jk})$ and thus $f_j$ become variable, the spaces $\bK_j$ rotate with the average speeds $O(\mu^{-1})$ but the previous statements remain true.
\end{example}

\begin{example}\label{ex-1-3} At points disjoint from $\Lambda$ one can assume that magnetic field corresponds to
$\omega = x_1dx_1\wedge dx_2 + dx_3\wedge dx_4$; I used 2-form (\ref{1-9}) and redefined coordinates $x'=(x_2,x_3,x_4)$ in the obvious way. Assume that
the metrics in these coordinates is Euclidean $g^{jk}=\delta_{jk}$ and that $V=V(x_1,x_2)$. Then $a_1=a_1^0-{\frac 1 2}V$ and $a_2=a^0_2$ are movement integrals as in example \ref{ex-1-1} and the movement splits into two movements again: one of them is the cyclotron movement in $(x_3,x_4)$ and another one is described by Hamiltonian $a_1$ in coordinates $(x_1,x_2)$. It was studied in details in \cite{IRO6}.

Then there are \emph{outer} and \emph{inner} zones $\cZ_\out$ and $\cZ_\inn$ respectively; in the outer zone $x_1$ should be much larger than the cyclotron radius (associated with $a_1$) which is $O(\mu^{-1} |x_1|^{-1})$ and this is the case as $\cZ_\out = \{ |x_1|\ge C \mu^{-1/2}\}$. However now $a^0_1$ is not necessarily disjoint from 0 (which was an assumption in \cite{IRO6}) and I will need to take it into account and redefine zones.

Still I can conclude from \cite{IRO6} that if a trajectory $(x(t),\xi(t))$ on an energy level $\le c$ starts from the point $(x(0),\xi(0))$ with
$x(0)\in B(0,{\frac 1 2})\cap \{x:|x_1|=\gamma \ge C^2 \mu^{-1/2}\}$, then
$x(t)\in B(0,{\frac 1 2})\cap \{x:C^{-1}\gamma \le |x_1|\le C\gamma \}$ as
$|t|\le \epsilon \mu \gamma^2 $. Further the drift speed is
$O(\mu^{-1}\gamma^{-2})$; this estimate is sharp as long as $a^0_1$ is disjoint from 0. So again I have the cyclotron movement with angular velocities $\asymp \mu |x_1|$ and $\mu$ and the drift movement mainly along $x_2$.

On the other hand, if a trajectory $(x(t),\xi(t))$ on an energy level $\le c$ starts from the point $(x(0),\xi(0))$ with
$x(0)\in B(0,{\frac 1 2})\cap \{x:|x_1|\le C \mu^{-1/2}\}$ then
$x(t)\in B(0,{\frac 1 2})\cap \{x: |x_1|\le C^2\mu^{-1/2}\}$ as
$|t|\le \epsilon $ and the averaged propagation speed is $O(1)$; this estimate is sharp for the ``typical'' (in the heuristic sense) trajectory. In this case there are no separate cyclotron and drift movements along $\bK_1$ but the cyclotron movement along $\bK_2$ remains.
\end{example}

\begin{example}\label{ex-1-4} Assume that in some coordinates $(x_1,\dots,x_4)$ 2-form $\omega$ is defined by (\ref{1-6}) and also $g^{jk}=\delta_{jk}$. Let $V=\const$. I am going to show that the conclusion of the previous example remains true (with the drift directed mainly along magnetic lines). Now however the angle between $\bK_1(x)$ and $T_x\sigma$ is of magnitude $r$ and therefore to prevent cyclotron movements from hitting $\Sigma$ it is sufficient to assume that $C\mu^{-1}\gamma^{-1}r \le \gamma$ and to keep $r$ of the same magnitude along trajectory it is sufficient to assume that $r\ge C\mu^{-1}\gamma^{-1}$. Alternatively I can assume that $r\le C\mu^{-1}\gamma^{-1}$ and then it will remain this way and the first inequality should be fulfilled with $r=C\mu^{-1}\gamma^{-1}$. So I get zone
\begin{multline}
\cZ_\out=\bigl\{C\mu^{-1}\gamma^{-1}\le r\le \epsilon \mu \gamma^2\bigr\}\cup
\bigl\{r\le C\mu^{-1}\gamma^{-1}, C\mu^{-1}\gamma^{-1}\le \epsilon \mu \gamma^2\bigr\}=\\
\bigl\{\gamma \ge C\max (\mu^{-2/3},\mu^{-1/2}r^{1/2})\bigr\}.
\label{1-11}
\end{multline}

\subsubsection{}
\label{sect-1-2-2} Let instead of rather heuristic arguments above apply precise calculations. With $\omega$ defined by (\ref{1-6}) and $g^{jk}=\delta_{jk}$
\begin{multline}
a= {\frac 1 2}\Bigl(\xi_1^2 +\xi_r^2 +W(x_1,r)\Bigr),\\
W(x_1,r)=\Bigl(\xi_2 -\mu (x_1-{\frac 1 2}r^2)\Bigr)^2 + r^{-2}\Bigl(\xi_\theta-\mu (x_1-{\frac 1 2}r^2)r^2\Bigr)^2= \\
(r^2+1)\Bigl(\mu \bigl(x_1-{\frac 1 2}r^2\bigr) -
{\frac {\xi_\theta -\xi_2r^2}{r^2+1}}\Bigr)^2
+{\frac 1 {r^2(r^2+1)} }\bigl(\xi_\theta -\xi_2 r^2\bigr)^2.
\label{1-12}
\end{multline}
Then $\xi_\theta$ and $\xi_2$ are movement integrals and the dynamics is restricted to the zone
$\cY_{\xi_2,\xi_\theta}=\bigl\{(x_1,r): W(x_1,r)\le c\bigr\}$. While $W$ is nonnegative and convex in the vicinity of $(0,0)$ it is a degenerate as $\xi_2=\xi_\theta=0$: $W= \mu^2(1+r^2) \bigl(x_1-{\frac 1 2}r^2\bigr)^2$ and so particle seems to move far away if time is unbounded. So I will need to bound time (but this bound will be large enough to sustain a nice remainder estimate).
\end{example}

Now I can consider rigorously the general case.

Let $\Bbbk_\alpha= (\Bbbk_\alpha ^m )$ be eigenvectors of $(F^j_{,k})=(g^{jl})(F_{lk})$ corresponding to eigenvalues $i f_\alpha$ ($\alpha =1,2$); then $\Bbbk_\alpha^\dag$ are eigenvectors corresponding to eigenvalues $-i f_\alpha$:
\begin{equation}
\sum_k F^j_{,k}\Bbbk^k_\alpha = if_\alpha \Bbbk^j_\alpha, \qquad
\sum_k F^j_{,k}\Bbbk^{\dag k}_\alpha = -if_\alpha \Bbbk^{\dag j}_\alpha.
\label{1-13}
\end{equation}
I normalize them so that
\begin{equation}
\sum_{j,k} g_{jk}\Bbbk^j_\alpha \Bbbk^k_\beta =0, \qquad \sum_{j,k} g_{jk}\Bbbk^{\dag j}_\alpha \Bbbk^{\dag k}_\beta =0,\qquad
\sum_{j,k} g_{jk}\Bbbk^j_\alpha \Bbbk^{\dag k}_\beta =2\delta_{\alpha\beta}.
\label{1-14}
\end{equation}
where the second equations in (\ref{1-13}),(\ref{1-14}) follow from the first ones.

\begin{claim}\label{1-15} One can select $\Bbbk_j\in C^\infty$ satisfying
(\ref{1-13}), (\ref{1-14}) in domain $\Omega= B(0,1)\cap \{|x_1\le \epsilon\}$\,\footnote{\label{foot-4} Provided I redefine temporarily $f_1$ so that $f_1/x_1 >0$.}.\end{claim}

Really, this statement is true for $\Bbbk_2$; therefore $\bK_1\in C^\infty$ and $(F^j_{,k})$ transforms $\bK_1\in C^\infty$ into itself and is skew-symmetric with respect to $(g^{jk})$. Then in orthonormal real base on $\bK_1$ $F=\begin{pmatrix} 0 &f_1\\ -f_1 &0\end{pmatrix}$ with $f_1\in C^\infty$ and one can select $\Bbbk_1\in C^\infty$.

Let
\begin{multline}
Z_\alpha(x,\xi) = \sum_m \Bbbk_\alpha ^m(x) p_m(x,\xi), \qquad
Z_\alpha ^\dag (x,\xi) = \sum_m \Bbbk_\alpha ^{\dag m}(x)p_m (x,\xi),\\
p_m=\xi_m-\mu V_m(x),\qquad \{p_j,p_k\}=-\mu F_{jk}
\label{1-16}
\end{multline}
where the second equation follows from the first one. Then since

\begin{claim}\label{1-17}
$|p_j(x,\xi)|\le c_1$ on the energy levels $\{(x,\xi):a (x,\xi)\le c_0\}$
\end{claim}
\noindent
I arrive to
\begin{equation}
\{ Z_\alpha, Z_\beta\} \equiv \{ Z_\alpha^\dag, Z^\dag _\beta\} \equiv 0, \quad
\{ Z_\alpha, Z_\beta^\dag\} \equiv 2 i\mu f_\alpha\delta_{\alpha\beta}\qquad \mod O(1).
\label{1-18}
\end{equation}

\subsubsection{}
\label{sect-1-2-3} Further, note that corrected $x_j$:
\begin{equation}
x'_j=x_j-\mu^{-1} \sum_k {\check F}^{jk}p_k, \qquad
\sum_k{\check F}^{jk}F_{kl}=\delta_{jl}
\label{1-19}
\end{equation}
satisfy
\begin{align}
&\{x'_j, p_k\}= O(\mu^{-1}\gamma^{-2}),\label{1-20}\\
&\{x'_j,x'_k\}=\mu^{-1}{\check F}^{jk}+ O(\mu^{-2}\gamma^{-3}).\label{1-21}
\end{align}

Also note that
\begin{equation}
a^0 (x,\xi)\Def {\frac 1 2}\sum_{k,l} g^{kl}p_kp_l =
{\frac 1 2}\Bigl( |Z_1|^2 + |Z_2|^2\Bigr)
\label{1-22}
\end{equation}
and therefore
\begin{equation}
\{ a^0, Z_\alpha\} \equiv \{ a, Z_\alpha\} \equiv - i\mu f_\alpha Z_\alpha,
\qquad
\{ a^0, Z_\alpha^\dag\} \equiv \{ a, Z_\alpha^\dag\} \equiv i\mu f_\alpha Z_\alpha^\dag \mod O(1)
\label{1-23}
\end{equation}
and
\begin{equation}
\{ a^0, b_\alpha \} \equiv \{ a, b_\alpha \} \equiv 0 \mod O(1), \qquad b_\alpha\Def |Z_\alpha|^2.
\label{1-24}
\end{equation}
More precisely, (\ref{1-18}) are fulfilled modulo linear forms\footnote{\label{foot-5} With respect to $(Z_1,Z_1^\dag,Z_2,Z_2^\dag)$ with smooth complex coefficients depending on $x$.} and therefore $\{a^0, b_\alpha\}$ is a cubic form$^{\ref{foot-5}}$ and
$\{a, b_\alpha\}$ is a cubic form$^{\ref{foot-5}}$ plus a linear form$^{\ref{foot-5}}$. Moreover, these cubic terms should contain at least one factor $Z_1$ or $Z_1^\dag$ and also at least one factor $Z_2$ or $Z_2^\dag$.

\begin{claim}\label{1-25}
Then correcting $b_\alpha$ ($\alpha=1,2$) by cubic terms$^{\ref{foot-5}}$ plus linear terms$^{\ref{foot-5}}$, both multiplied by $\mu^{-1}$, one can eliminate from
$\{a, b_\alpha\}$ all cubic terms but those with
$ Z_1Z_2 Z_2^\dag$ and $ Z_1^\dag Z_2 Z_2^\dag$ and all linear terms but those with $Z_1$, $Z_1^\dag$. Then these remaining terms produce $O(\mu^{-1})$ error.
\end{claim}

Further, correcting in the zone $\cZ_\out=\{C\mu^{-1/2}\le |x_1|\le \epsilon\}$ by cubic and linear terms of the above type but with coefficients of the type
``$x_1^{-1}\times \underline{\textsf{smooth function}}$'' one can eliminate the rest of the terms but now instead of $O(\mu^{-1})$ error I get $O(\mu^{-1}x_1^{-2})$ error:

\begin{claim}\label{1-26}
In the classical dynamics one can correct $b_\alpha=|Z_\alpha|^2$ by $O(\mu^{-1}x_1^{-1})$ so that their propagation speeds would be $O(\mu^{-1}x_1^{-2})$. Therefore as long as $x(t)\in \Omega$ and
$|x_1(t)|\asymp \gamma$, $|t|\le T$ the following is fulfilled:
\begin{equation*}
|Z_\alpha (t)|^2 = |Z_\alpha (0)|^2 + O\bigl(\mu^{-1}\gamma^{-1}+\mu^{-1}\gamma^{-2}T\bigr).
\end{equation*}
\end{claim}

Similarly, consider $\{a, x_j\}$ which is a linear form$^{\ref{foot-5}}$. One can correct $x_j$ in the same way arriving to

\begin{claim}\label{1-27}
In the classical dynamics one can correct $x_j$ by $O(\mu^{-1}x_1^{-1})$ so that their propagation speeds would be $O(\mu^{-1}x_1^{-2})$.
Therefore as long as $x(t)\in \Omega$ and $|x_1(t)|\asymp \gamma$, $|t|\le T$
\begin{equation*}
x_j(t)= x_j(0) + O\bigl(\mu^{-1}\gamma^{-1}+\mu^{-1}\gamma^{-2}T\bigr).
\end{equation*}
\end{claim}

\vskip-20pt
Moreover, (\ref{1-24}) and (\ref{1-25}) with $b_\alpha$ replaced by $x_j$ corrected yields that
\begin{equation}
{\frac {dx}{dt}} \in \bK_1(x)+O(\mu^{-1}x_1^{-1})
\qquad \text{for } x, \text{ corrected modulo } O(\mu^{-1}x_1^{-1}).
\label{1-28}
\end{equation}
Still it is not yet what I want which is to calculate $\{a, x_j\}$ modulo $O(\mu^{-1}x_1^{-1})$. However, let us consider first
\begin{equation}
b'_3\Def-\Re \Bigl(i\mu^{-1}\{Z_1^\dag,x_1\}\phi(x)^{-1} Z_1\Bigr) + x_1^2,\qquad \phi = f_1/x_1
\label{1-29}
\end{equation}
and note that $\{Z_1Z_1^\dag , b'_3\}=O(\mu^{-1})$. Also note that $\{a,b'_3\}$ is a linear form of $Z_2,Z_2^\dag$ modulo $O(\mu^{-1})$ and therefore it could be corrected by a linear form of $Z_2,Z_2^\dag$ multiplied by $\mu^{-1}$ so that after correction $\{a,b'_3\}=O(1)$. Thus

\begin{claim}\label{1-30}
In the classical dynamics one can correct $b_3=x_1^2$ by $O(\mu^{-1})$ so that its propagation speed would be $O(\mu^{-1})$. Therefore
$x_1^2(t)= x_1^2(0) + O\bigl(\mu^{-1}+\mu^{-1}T\bigr)$ and statements $(\ref{1-26})$, $(\ref{1-27})$ hold without precondition ``as long as $x(t)\in \Omega$ and $|x_1(t)|\asymp \gamma$'' which is fulfilled automatically for
$x(0)\in B(0,{\frac 1 2})\cap \{{\frac 1 2}\gamma\le |x_1(0)|\le 2\gamma\}$,
$T=\epsilon \mu \gamma^2$, $\gamma\ge {\bar\gamma}_1=C\mu^{-1/2}$.
\end{claim}

\vskip-10pt
\begin{claim}\label{1-31}
Also, trajectory originated in $B(0,{\frac 1 2})\cap \{|x_1|\le \gamma)$ remains in $B(0,1)\cap \{|x_1|\le C\gamma\}$ as $\gamma={\bar\gamma}_1$, $T=\epsilon$.
\end{claim}

Further, as long as $|x_1|\le c\gamma$ and $|\{Z_1,x_1\}|\le cr$, correction is $O\bigl(\mu^{-1}(r+\gamma)\bigr)$ and therefore
$x_1^2(t)=x_1^2(0) + O\bigl(\mu^{-1}(r+\gamma)+\mu^{-1}T\bigr)$. Then I arrive to

\begin{claim}\label{1-32}
As long as $|\{Z_1,x_1\}|\le cr$ along trajectory, statements (\ref{1-30}), (\ref{1-31}) hold with ${\bar\gamma}_1 = C\mu^{-1}+ C\mu^{-1/2}r^{1/2}$.
Moreover $|\{Z_1,x_1\}(t)|\le cr$ provided it was fulfilled (with constant $c/2$) at $t=0$ and $T=\epsilon \mu \gamma^2 r$, $r\ge \mu^{-1/2}+\gamma$.
\end{claim}

\subsubsection{}
\label{sect-1-2-4} Now let us return to calculation of $\{a, x_j\}$ modulo $O(\mu^{-1})$ where I will correct $x_j$. Any error, larger than this, comes from eliminating linear terms containing $Z_1$ or $Z_1^\dag$, namely
$\Re \bigl(\{Z_1^\dag ,x_j\}Z_1\bigr)$. This expression was eliminated by adding
$ \mu ^{-1}f_1^{-1}\Re (i \{Z_1 ,x_j\}Z_1^\dag)$ to $x_j$ which in turn generates an error
\begin{equation}
\mu^{-1}\Re i\Bigl( \bigl\{a, f_1^{-1}\{Z_1 ,x_j\}\bigr\}Z_1 ^\dag +
f_1^{-1}\{Z_1^\dag ,x_j\}\bigl(i\{a,Z_1^\dag\}-\mu f_1Z_1\bigr) \Bigr).
\label{1-33}
\end{equation}
In the expression (\ref{1-33}) the part containing factors $Z_2$ and $Z_2^\dag$ in the combinations $Z_2$, $Z_2^2$, $Z_2^\dag$, $Z_2^{\dag\,2}$ can be eliminated
by the same way as above; the corresponding correction has an extra factor $O(\mu^{-1})$ and the corresponding error just acquires factor $O(\mu^{-1}x_1^{-1})$; therefore resulting error will be $O\bigl(\mu^{-2}x_1^{-3}|Z_1| + \mu^{-2}x_1^{-2}\bigr)= O\bigl(\mu^{-1}+\mu^{-1}\gamma^{-1}\rho\bigr)$ provided
\begin{equation}
|Z_1|\le c\rho, \qquad \gamma \ge c\mu^{-1/2}.
\label{1-34}
\end{equation}
This leaves us with expression (which I do not call an error anymore)
\begin{equation}
{\frac 1 2}\Bigl(-\mu ^{-1}x_1^{-2}\phi^{-1}\{Z_1 Z_1^\dag,x_1\}\Re (i\{Z_1,x_j\}Z_1^\dag)+
\mu^{-1}x_1^{-1}\Re \bigl(\alpha \{Z_1,x_j\}\bigr) |Z_2|^2\Bigr)
\label{1-35}
\end{equation}
with the complex-valued coefficient $\alpha=\alpha (x')$ provided
\begin{equation}
f_1 = \phi (x')x_1 +O(x_1^2),\qquad x'=(x_2,x_3,x_4).
\label{1-36}
\end{equation}
Note that the first term in (\ref{1-35}) is equal to
\begin{multline}
-\mu ^{-1} x_1^{-2}\phi ^{-1}\cdot
\Re( \{Z_1 ,x_1\}Z_1^\dag) \cdot \Re (i \{Z_1 ,x_j\}Z_1^\dag)= \\
{\frac 1 4}\mu ^{-1} x_1^{-2}\phi ^{-1}\cdot
\Re i\biggl(2\{Z_1^\dag, x_1\}\{Z_1^\dag, x_j\}Z_1^2 +
\bigl( -\{Z_1^\dag, x_1\}\{Z_1,x_j\}+
\{Z_1^\dag, x_j\}\{Z_1,x_1\}\bigr)|Z_1|^2\biggr)
\label{1-37}
\end{multline}
and one needs to calculate coefficients at $Z_1^2$, $|Z_1|^2$ as $x_1=0$ only.
Thus I need to calculate $\Bbbk_1$ as $x_1=0$. Let us decompose $F_{jk}$, $g_{jk}$ and $k_1$ into powers of $x_1$; then $F_{jk}\Bbbk_1= if_1\Bbbk_1$ implies that
\begin{equation}
F^0_{jk}\Bbbk^0_1=0,\qquad
(F^1_{jk}-i\phi g^0_{jk}) \Bbbk^0_1+F^0_{jk}\Bbbk^1_1=0,
\label{1-38}
\end{equation}
where $F_{jk}$ corresponds to symplectic form (\ref{1-6}) and therefore
\begin{equation}
F=\begin{pmatrix}
0 &1 &-x_4 &x_3\\
\\
-1 & 0 &x_3 & x_4\\
\\
x_4 & -x_3 & 0 &2\bigl(x_1 -{\frac 1 2}(x_3^2+x_4^2)\bigr)\\
\\
-x_3 &-x_4 &-2\bigl(x_1 -{\frac 1 2}(x_3^2+x_4^2)\bigr) &0
\end{pmatrix}.
\label{1-39}
\end{equation}
Then (\ref{1-38}), (\ref{1-39}) imply for the ``differentiation'' part of $Z_1$
\begin{multline}
Z_{1,\ \text{diff}}|_{x_1=0}=\alpha\Bigl(\bigl(r^{-1}\partial_\theta -r\partial_2\bigr) +
\beta \bigl(r\partial_1 +\partial_r\bigr)\Bigr)\\
\text{with}\quad c^{-1}\le |\alpha|\le c,\;
c^{-1}\le \Re i\beta \le |\beta |\le c.
\label{1-40}
\end{multline}
Then
\begin{equation*}
\Re i\bigl(\{Z_1^\dag, x_1\}Z_1 \bigr)\bigr|_{x_1=0}= 2|\alpha|^2\Re (i\beta)\bigl(\partial_\theta -r^2\partial_2\bigr)
\end{equation*}
is the differentiation along magnetic lines. Further, without any loss of the generality one can assume that $\alpha=1$.

\subsubsection{}
\label{sect-1-2-5} Consider first $\psi=x_1-r^2/2$ which is a regular function of $x$; note that $Z_1\psi=O(x_1)$; therefore $\psi$ can be corrected by $O(\mu^{-1})$ so that ${\frac {d\ }{dt}}\psi = O(\mu^{-1})$ and therefore original $\psi$ is preserved modulo $O\bigl(\mu^{-1}(T+1)\bigr)$. Since I already know that $x_1$ is preserved modulo $O\bigl(\mu^{-1}\gamma^{-1}r (T+1) \bigr)$ I conclude that

\begin{proposition}\label{prop-1-5} Consider classical dynamics in
$\Omega = B(0,1)\cap \{|x_1|\le C\epsilon\}$, originated in
$\Omega' = B(0,{\frac 1 2})\cap \{|x_1|\le \epsilon\}$.

\smallskip\noindent
(i) If in the original point $|x_1|=\gamma$ with
$\gamma \ge C\max(\mu^{-1}r^{-1}, r^2)$ then both $|x_1(t)|$ and $r(t)$ remain of the same magnitudes as $T=\epsilon \mu \gamma r$.

\smallskip\noindent
(ii) If in the original point
$r^2\ge \gamma\ge C\mu^{-1/2}r^{1/2}$ then both $|x_1(t)|$ and $r(t)$ remain of the same magnitudes as $T=\epsilon \mu \gamma^2 r^{-1}$.

\smallskip\noindent
(iii) If in the original point $\gamma\ge C\mu^{-2/3}$ and $r\le c\mu^{-1}\gamma^{-1}$ then $|x_1(t)|$ remains of the same magnitude while $r(t)\le 2c\mu^{-1}\gamma^{-1}$ as $T=\epsilon$.

\smallskip\noindent
(iv) If in the original point $|x_1|\le \gamma = c\mu^{-1/2}r^{1/2}$ and
$r\ge C\mu^{-1/3}$ then $r(t)$ remains of the same magnitude while $|x_1(t)|\le C\gamma$ as $T=\epsilon $.

\smallskip\noindent
(v) If in the original point $|x_1|\le \gamma = c\mu^{-2/3}$ and
$r\le c\mu^{-1/3}$ then $r(t)\le C\mu^{-1/3}$ and $|x_1(t)|\le C\gamma$ as $T=\epsilon $.
\end{proposition}

On the other hand, at the energy level $\tau$ one can rewrite the last term in (\ref{1-35}) as \newline$\mu^{-1}x_1^{-1}\Re (\sigma\{Z_1,x_j\})(\tau+V)$ modulo $O(\mu^{-1}+\mu^{-1}\gamma^{-1}\rho^2)$. Consider (\ref{1-38}). Using the same technique as before one can eliminate
\begin{equation*}
\mu ^{-1} x_1^{-2}\phi ^{-1}\cdot
\Re i\Bigl(2\{Z_1^\dag, x_1\}\{Z_1^\dag, x_j\}Z_1^2\Bigr)
\end{equation*}
by $\mu ^{-2} x_1^{-3}\phi ^{-2}\cdot
\Re \Bigl(\{Z_1^\dag, x_1\}\{Z_1^\dag, x_j\}Z_1^2\Bigr)$ correction
(which is $O(\mu^{-2}\gamma^{-3}r\rho^2)=O(\mu^{-1}\gamma^{-1}\rho^2)$) as $\gamma \ge \mu^{-1/2}r^{1/2}$ adding $O(\mu^{-2}\gamma^{-4}r^2\rho^2+\mu^{-2}\gamma^{-3}\rho^2)$ error; this error is less than $\epsilon \mu^{-1}\gamma^{-2}r\rho^2$ as
$\gamma \ge C\mu^{-1/2}r^{1/2}$.

Also, any $Z_1$ or $Z_1^\dag$ unbalanced in the product could be treated in the same way as before.

Then I arrive to

\begin{proposition}\label{prop-1-6}
As $\gamma \ge C\mu^{-{1/2}}r^{1/2}$ the propagation speed with respect to $\theta$ is
\begin{equation}
v(\rho ;x,\mu^{-1})=\mu^{-1}\gamma^{-2}\kappa (x,\rho^2, \mu) \Bigl(\rho^2 +
r^{-1} \kappa_2 (x,\rho^2,\mu)\gamma \Bigr) +
O\Bigl( \mu^{-2}\gamma^{-3} \rho r^{-1}\Bigr)
\label{1-41}
\end{equation}
with $\kappa,\kappa_2$ bounded and $\kappa$ disjoint from $0$.
\end{proposition}

I also need to consider dynamics of $b_2=|Z_2|^2$ more precisely. Note that
\begin{equation*}
\bigl \{|Z_2|^2,a\}={\frac 1 2}\Re \biggl( \Bigl( \{Z_2, Z_1\}Z_1 ^\dag + \{Z_2, Z_1^\dag\}Z_1 -\{Z_2,V\}\Bigr)Z_2^\dag\biggr)
\end{equation*}
is a combination of cubic and linear terms (see footnote $^{\ref{foot-5}}$) and the only terms where factor $Z_1$ is not compensated by $Z_1^\dag$ and v.v. come from
\begin{align}
\{Z_2,Z_1\}= &\alpha_{11} Z_2 + \alpha_{12} Z_2^\dag + \beta _{11}Z_1 + \beta_{12}Z_1^\dag,\label{1-42}\\
\{Z_2,Z_1^\dag\}= &\alpha_{21} Z_2 + \alpha_{22} Z_2^\dag + \beta _{21}Z_1 + \beta_{22}Z_1^\dag;\label{1-43}
\end{align}
these terms are equal to $\Re (\sigma Z_1)|Z_2|^2$ with $\sigma=\alpha_{11}^\dag +\alpha_{21}$.
I need to calculate $\beta$ as $x_1=0$. All other terms could be corrected by $O(\mu^{-1})$ leading to $O(\mu^{-1})$ error while this one is corrected by
$\mu^{-1}x_1^{-1}\phi \Re (i\beta Z_1)|Z_2|^2$ leading to $O(\mu^{-1}x_1^{-2})$ error; I will calculate this latter error more precisely below. Since
\begin{align*}
&\alpha_{11}=(2i\mu f_2)^{-1}\bigl\{\{Z_2,Z_1\},Z_2^\dag\bigr\},\\
&\alpha_{21}=(2i\mu f_2)^{-1}\bigl\{\{Z_2,Z_1^\dag\},Z_2^\dag\bigr\}\\
\end{align*}
calculated as $Z_1=Z_2=0$, I conclude that
\begin{equation*}
\sigma =(2i\mu f_2)^{-1} \bigl( -\bigl\{\{Z_2,Z_1\},Z_2^\dag\bigr\} +
\bigl\{\{Z^\dag_2,Z_1\},Z_2\bigr\}=
-(i\mu f_2)^{-1}\bigl\{\{Z_2,Z_2^\dag\},Z_1\bigr\}.
\end{equation*}

Therefore correction is $O(\mu^{-1}\gamma^{-1}\rho)$ and then the part of an error, with factors $Z_2$ and $Z_2^\dag$ balanced is $O(\mu^{-1}\gamma^{-1}+\mu^{-1}\gamma^{-2}\rho r)$.

Moreover, the part of $\{|Z_2|^2,a\}$ where $Z_2$ and $Z_2^\dag$ are balanced is $\kappa |Z_2|^2$ with
\begin{equation*}
\kappa=(2\mu f_2)^{-1}\Re \Bigl(\bigl\{i \{Z_2,Z_2^\dag\}, Z_1\bigr\}\bigr|_{Z_1=Z_2=0 }Z_1^\dag\Bigr)
\end{equation*}
since $i \{Z_2,Z_2^\dag\}$ is real; further, this expression is equal to
\begin{multline*}
(2f_2)^{-1}\{f_2,|Z_1|^2\} + O(\rho\gamma)+O(\mu^{-1}\gamma^{-1})=\\
(2f_2)^{-1}\{f_2,a\} - (2f_2)^{-1}\{f_2,|Z_2|^2\}+O(\rho\gamma)+O(\mu^{-1}\gamma^{-1}).
\end{multline*}
Note that $f_2^{-1}\{f_2,|Z_2|^2\}$ has unbalanced $Z_2$ or $Z_2^\dag$ and thus is eliminated by $O(\mu^{-1})$ correction. So I arrive to an error
$-(2f_2)^{-1}\{f_2,a\}|Z_2|^2$.

Thus

\begin{claim}\label{1-44} (i)
One can correct $|Z_2|^2$ by $O(\mu^{-1}+\mu^{-1}\gamma^{-1}\rho)$ term so that the propagation speed after correction would be $O(\mu^{-1}\gamma^{-1}+\mu^{-1}\gamma^{-2}\rho r)$.

\smallskip\noindent
(ii) Moreover, one can correct $f_2^{-1}|Z_2|^2$ by $(\mu^{-1}+\mu^{-1}\gamma^{-1}\rho^2)$ term so that the propagation speed after correction would be $O(\mu^{-1}+\mu^{-1}\gamma^{-2}\rho^2)$.
\end{claim}

I will also need the shorter term classical dynamics result:

\begin{proposition}\label{prop-1-7} Let us consider dynamics described in proposition \ref{prop-1-5}(i)-(v). Then \newline $\dist (x(t),(x(0))\ge \epsilon \rho t$ as $|t|\le T_1 =\epsilon \min(\mu ^{-1}\gamma^{-1},\rho)$.
\end{proposition}

\begin{remark}\label{rem-1-8}
It effectively makes $T_1=\epsilon \mu ^{-1}\gamma^{-1}$ because contribution to the remainder estimate of the subzone $\{ |x_1|\asymp \gamma, |Z_1|\le \rho, \rho \le \mu^{-1}\gamma^{-1}\}$ is $O(\mu h^{-3} \times \mu^{-2}\gamma^{-2}\times \gamma)=O(\mu^{-1}h^{-3}\gamma^{-1})$ which after summation over zone $\{\gamma\ge \mu^{-1/2}\}$ results in $O(\mu^{-1/2}h^{-3})$.
\end{remark}

\subsection{Quantum dynamics. I. Outer zone}
\label{sect-1-3}

In this subsection begin to deal with the quantum dynamics and the new parameter $h$ appears; so one can now compare $\mu$ and $h$. Here I assume that magnetic field is rather weak; more precise requirements should vary from statement to statement but actually I do not care very much now and will just assume as needed that
\begin{equation}
C\le \mu \le h^{-\delta}
\label{1-45}
\end{equation}
with small enough $\delta>0$; I do not want to push it up so far.

\subsubsection{}
\label{sect-1-3-1} As I have mentioned there are few time scales and the shortest one is
$\asymp \mu^{-1}$. Let us deal with it first. Scaling $x\mapsto \mu x$, $t\mapsto \mu t$, $h \mapsto \mu h$, $\mu \mapsto 1$ I find ourselves with the standard Schr\"odinger operator with the propagator which is a standard Fourier integral operator. Returning to the original scale I get

\begin{proposition}\label{prop-1-9} (i) As $|t|\le T_0=C\mu^{-1}$ propagator $e^{ih^{-1}A t}$ is $h$-FIO corresponding to the classical dynamics $\Phi_t=e^{t H_a}$ with Hamiltonian $a(x,\xi)$.
\smallskip\noindent
(ii) In particular, if $\psi \in C_0^\infty (B({\bar x}, \epsilon\gamma))$
\begin{equation}
e^{2\pi i\mu^{-1}h^{-1}{\bar f}_2^{-1}A}\psi
\equiv e^{2\pi i\mu^{-1}h^{-1}{\bar f}_2^{-1} B}\psi
\label{1-46}
\end{equation}
with ${\bar f}_2=f_2({\bar x})$ and the standard $h$-pseudodifferential operator $\gamma ^{-1}B$.
\end{proposition}

\begin{proof} Let us consider classical dynamics $\Phi_t$ first; in $(x,p)$ variables it is described by
\begin{align}
&{\frac {dx_j}{dt}}= \mu^{-1}\sum_k K_{jk}(x)p_k,\label{1-47}\\
&{\frac {dp_j}{dt}}=\sum_k L_{jk}(x) p_k +\mu^{-1} \sum_{kl} L_{jkl}(x) p_kp_l +
\mu^{-1} L_j(x)
\label{1-48}
\end{align}
(after rescaling $t\mapsto \mu^{-1}t$) with uniformly smooth $K_*,L_*$. Therefore in these coordinates $\Phi_t: (x,Z)\mapsto \bigl(x+ \mu^{-1} X(x,p,t), Y(x,p,t))$. Thus $e^{ih^{-1}A t}$ is a product of FIO corresponding to the symplectomorphism $\Phi_t$ (a standard quantization, phase function is defined in the standard way and the symbol is just 1) and some PDO.

In particular one can see easily that $\Phi_t -I =O(\gamma)$
as $t=2\pi {\bar f}_2^{-1}$; therefore I conclude that $\Phi_t = e^{\gamma H_b}$ with some symbol $b=b(\mu^{-1}x,Z)$. Quantizing it I get that (\ref{1-46}) holds but with an extra PDO factor $\cQ$ between FIO and $\psi$. However then one can perturb $B$ by operator $\gamma^{-1}h B_1$ so that
\begin{equation*}
e^{2\pi i\mu^{-1}h^{-1}{\bar f}_2^{-1} B}\cQ \equiv e^{2\pi i\mu^{-1}h^{-1}{\bar f}_2^{-1} (B+\gamma^{-1}hB)};
\end{equation*}
I leave this easy exercise to the reader.
\end{proof}

\begin{proposition}\label{prop-1-10} Let condition
\begin{equation}
V\ge \epsilon_0
\label{1-49}
\end{equation}
be fulfilled. Then
\begin{equation}
|F_{t\to h^{-1}\tau}\chi_T(t)\Gamma \bigl(u \psi Q\bigr) |\le C'h^s
\qquad \forall \tau :|\tau|\le\epsilon
\label{1-50}
\end{equation}
as $\mu \ge K$, ${\bar T}\Def Ch|\log h|\le T\le 2\pi (1-\epsilon_0)\mu^{-1}{\bar f}_2^{-1}$, $\gamma \le \epsilon$ where $\psi $ is supported in $B({\bar x},\gamma)$, $\gamma $ is calculated in ${\bar x}$, $Q$ is arbitrary $h$-pdo \footnote{\label{foot-6} In our usual manner here and until the end of the paper
$\chi\in C_0^\infty( [-1,-{\frac 1 2}]\cup[{\frac 1 2},1])$, ${\bar\chi}\in C_0^\infty ([-1,1])$ are functions of \cite{IRO1} type, $\chi_T(t)=\chi (t/T)$ etc. and $\Gamma_x (u)=u(x,x)$, $\Gamma u = \int u(x,x)\,dx$.}.
\end{proposition}

Now let us consider in the outer zone the intermediate scale dynamics:

\begin{proposition}\label{prop-1-11} Consider
${\bar x}\in \cZ_\out\Def \{\gamma \ge C\mu^{-1/2}\}$.
Then as $|t|\le C\mu^{-1}\gamma^{-1}$ \ $e^{ih^{-1}A t}$ is $h$-FIO corresponding to the classical dynamics $\Phi_t=e^{t H_a}$ with Hamiltonian $a$
and
\begin{equation}
e^{ih^{-1}tA }\psi\equiv e^{ih^{-1}t''A }e^{ih^{-1} t'B }\psi
\label{1-51}
\end{equation}
where $t'=2\pi \mu^{-1} {\bar f}_2^{-1} n$,
$n=\bigl\lfloor (2\pi)^{-1}{\bar f} \mu t \bigr\rfloor$,
$t''=t-t'$.
\end{proposition}

\begin{proof} Formula (\ref{1-51}) follows from (\ref{1-46}) and the fact that $|x_1|$ retains its magnitude in the classical dynamics in the zone $\cZ_\out$.
\end{proof}

\begin{proposition}\label{prop-1-12} (i) Let
$\mu \le \epsilon h^{-1}|\log h|^{-1}$,
${\bar x}\in \cZ_\out$ and $Q$ be $h$-PDO with the symbol supported in $\{|Z_1| \ge \rho \}$ with
\begin{equation}
\rho \ge C\mu^{-1}\gamma^{-1} + C(\mu h \gamma |\log h|)^{1/2}.
\label{1-52}
\end{equation}
Then $(\ref{1-50})$ holds as
$C\rho^{-2}h|\log h| \le T\le T_1\Def \epsilon_1\mu^{-1}\gamma^{-1}$.

\smallskip\noindent
(ii) In particular, let
\begin{equation}
\rho \ge C\mu^{-1}\gamma^{-1} + C(\mu h |\log h|)^{1/2}.
\label{1-53}
\end{equation}
Then estimate $(\ref{1-50})$ holds as
${\bar T}=Ch|\log h| \le T\le \epsilon_1\mu^{-1}\gamma^{-1}$.
\end{proposition}

\begin{proof} Let us temporarily direct $\bK_1$ as a coordinate plane\footnote{\label{foot-7} For $T\le c\mu^{-1}\gamma^{-1}$ and $\gamma\ge C\mu^{-1/2}$ with $C=C(\epsilon,c)$ our propagation is confined to $B({\bar x},\epsilon\gamma)$.} $\{x_1,x_2\}$ and consider $\rho$-admissible partition in $(\xi_1,\xi_2)$. Then one can see easily that the propagation along $\bK_1$ has the speed $\asymp \rho$. Further, condition (\ref{1-52}) ensures that $\rho$ retains its magnitude since the propagation speed of $|Z_1|^2$ does not exceed $C_0\mu^{-1}\rho ^2 \gamma^{-2}+C_0\mu^{-1}$ and $|Z_1|^2$ is corrected by $O(\rho \mu^{-1})$. Furthermore, the last term in (\ref{1-52}) ensures the logarithmic uncertainty principle.

Then for time $T$ the shift along $\bK_1$ is exactly of magnitude $\rho T$ and the logarithmic uncertainty principle
$\rho T \times \rho \ge Ch|\log h|$ is fulfilled as long as $T\ge T'=C\rho^{-2}h|\log h|$. This implies (i).

Now it follows from (i) and proposition \ref{prop-1-10} that estimate (\ref{1-50}) holds as ${\bar T}=Ch|\log h|\le T\le T_1=\epsilon \mu^{-1}\gamma^{-1}$ provided
$T'\le \epsilon \mu^{-1}$ which means upgrade of condition (\ref{1-52}) to (\ref{1-53}).
\end{proof}

\subsubsection{}
\label{sect-1-3-2} On the other hand one can see easily that the contribution to the remainder estimate of the rather thin subzone where $|x_1|\asymp \gamma$ and condition (\ref{1-53}) is violated does not exceed
$C\mu h^{-3} \bigl(\mu^{-2}\gamma^{-2}+\mu h|\log h|\bigr)\gamma$;
further, summation of this expression over $\cZ_\out$ results in
$C\mu^{-1/2}h^{-3}+C\mu^2h^{-2}|\log h|$. Therefore I have proven

\begin{proposition}\label{prop-1-13} Contribution to the remainder estimate of the subzone of $\cZ_\out$ where condition $(\ref{1-52})$ is violated is
$O\Bigl(\mu^{-1/2}h^{-3}+\mu^2h^{-2}|\log h|\Bigr)$. In particular this contribution is $O(\mu^{-1/2}h^{-3})$ as
$\mu \le C(h|\log h|)^{-2/5}$.
\end{proposition}

Actually, one can get rid off logarithmic factors here but under condition (\ref{1-45}) it is not needed; this estimate (without logarithmic factor) is the
best possible if no non-degeneracy condition is assumed.

From now on I will consider only the \emph{main subzone} of $\cZ_\out$ where condition (\ref{1-53}) is fulfilled and thus estimate (\ref{1-50}) holds with
${\bar T}\le T\le T_1$. I need to estimate its contribution to the remainder estimate.

First of all one needs the quantum version of the results of the previous subsection.

\begin{proposition}\label{prop-1-14} Let condition $(\ref{1-45})$ be fulfilled. Let function $\psi $ be properly supported in $B({\bar x},\gamma')$, and $Q$ be properly supported in $\{|Z_1|\le \rho\}$ with $\rho \ge C_0 \mu^{-1/2}$ and
\begin{equation}
\gamma ' \ge C\rho ^{-1}h|\log h|.
\label{1-54}
\end{equation}
Further, let $\gamma \ge C_0\mu^{-1/2}$ and $T\le T_2\Def \epsilon \mu \gamma^2$ and let $x''=(x_3,x_4)$. Then

\smallskip\noindent
(i) If
\begin{equation}
r =|{\bar x}''|\ge r'\Def C_0\mu^{-1}\gamma^{-1}+C_0\gamma'
\label{1-55}
\end{equation}
then estimate
\begin{equation}
|F_{t\to h^{-1}\tau } {\bar \chi}_T (t)(1-\psi _1) (u \psi_yQ_y^t) |\le Ch^s
\label{1-56}
\end{equation}
holds provided $\psi_1=1$ in the domain
\begin{multline}
\Bigl\{x: C_0^{-1}\gamma \le |x_1| \le C_0\gamma,\
C_0^{-1}r\le |x''| \le C_0r,\\
\dist (x,{\bar x})\ge C_0\gamma' + C_0\mu^{-1}\gamma^{-1}+ C_0\mu^{-1}(r +\gamma)\gamma^{-2}T\Bigr\} .
\label{1-57}
\end{multline}

\smallskip\noindent
(ii) On the other hand, if condition $(\ref{1-55})$ is violated then estimate $(\ref{1-56})$ holds provided $\psi_1=1$ in the domain
\begin{multline}
\Bigl\{x: C_0^{-1}\gamma \le |x_1| \le C_0\gamma,\ |x''|\le C_0r',\\
\dist (x,{\bar x})\ge C_0\gamma' + C_0\mu^{-1}\gamma^{-1}+ C_0\mu^{-1}(r'+\gamma)\gamma^{-2}T\Bigr\} .
\label{1-58}
\end{multline}
\end{proposition}

\begin{proof} Proof is the standard one (see e.g. \cite{IRO1}) based on the analysis of the symbol
\begin{equation}
\varpi \biggl( {\frac 1 T}\Bigl( t \pm \epsilon v^{-1}\bigl(X (x,\xi) - X(y,\xi)\bigr)\Bigr)\biggr)
\label{1-59}
\end{equation}
where $\varpi$ is the same function as in \cite{IRO1}, $X(x,\xi)$ is one of the corrected symbols $x_1^2$,
$x_1-{\frac 1 2}|x''|^2$ and $x_j$ with $j=1,\dots,4$ and $v$ is the corresponding speed, namely $\mu^{-1}$, $\mu^{-1}$ and $\mu^{-1}x_1^{-2}r$ respectively.

Obviously, it is sufficient to consider $T$ such that the last term in definitions of (\ref{1-57}),(\ref{1-58}) is dominant: $\mu^{-1}(r+\gamma)\gamma^{-2}T \ge \gamma'+ \mu^{-1}\gamma^{-1}$. One can see easily that under this condition and (\ref{1-45}) symbol (\ref{1-59}) is quantizable; here also one can be more specific about exponent $\delta>0$ in condition (\ref{1-45}).
\end{proof}

Now I want to be more precise for smaller $r$, $\rho$:

\begin{proposition}\label{prop-1-15} In frames of proposition \ref{prop-1-14} as $T\le \epsilon \mu \gamma^2 \rho$ estimate
\begin{equation}
|F_{t\to h^{-1}\tau } {\bar \chi}_T (t)(1-Q_1) \psi_1(u \psi_yQ_y^t) |\le Ch^s
\label{1-60}
\end{equation}
holds provided $Q_1=I$ in the domain $\{c^{-1}\rho \le |Z_1|\le c\rho\}$ and estimate $(\ref{1-56})$ holds provided $\psi_1=1$ in the domain
\begin{multline}
\Bigl\{x: C_0^{-1}\gamma \le |x_1| \le C_0\gamma,\
C_0^{-1}r\le |x''| \le C_0r,\\
\dist (x,{\bar x})\ge C_0\gamma' + C_0\mu^{-1}\gamma^{-1}+
C_0\mu^{-1}\gamma^{-2} (r\rho^2+\gamma) T\Bigr\}
\label{1-61}
\end{multline}
as $r\ge r'$, and in the domain
\begin{multline}
\Bigl\{x: C_0^{-1}\gamma \le |x_1| \le C_0\gamma,\ |x''|\le C_0r',\\
\dist (x,{\bar x})\ge C_0\gamma' + C_0\mu^{-1}\gamma^{-1}+ C_0\mu^{-1}\gamma^{-2}(r'\rho^2+\gamma)T\Bigr\}.
\label{1-62}
\end{multline}
as $r\le r'$.
\end{proposition}

\begin{proof} Proof is standard, based on the symbol (\ref{1-59}) which is quantizable due to assumption (\ref{1-45}); here again $X$ is one of the corrected symbols $|Z_1|^2$, $x'$ and $v$ is the corresponding upper bound for a speed, namely $\mu^{-1}\gamma^{-2}\rho$, $\mu^{-1}\gamma^{-2}(\rho^2 + \gamma)$ respectively.
\end{proof}

Note that under condition (\ref{1-45}) $r'\asymp \mu^{-1}\gamma^{-1}$ and I want to eliminate the corresponding subzone from the future analysis. Really, let us note that the contribution of the subzone
$\bigl\{\gamma \le |x_1|\le 2\gamma, |x''|\le r, |Z_1|\le \rho \bigr\}$ to the remainder estimate is
$O(T_1^{-1}h^{-3}r^2 \rho^2 \gamma)= O(\mu h^{-3}r^2\rho^2 \gamma^2)$ (since I already eliminated subzone where condition (\ref{1-53}) is violated) and plugging $r\rho =\mu^{-1}\gamma^{-1}$ results in $O(\mu^{-1}h^{-3})$. Then summation with respect to $\gamma$ results in $O(\mu^{-1}h^{-3}|\log \mu|)$. So, I conclude that

\begin{claim}\label{1-63} Contribution to the remainder estimate of $\cZ_\out\cap \{r \rho \le C\mu^{-1}\gamma^{-1}\}$ is $o(\mu^{-1/2}h^{-3})$.
\end{claim}

\subsubsection{}
\label{sect-1-3-3} From now on I restrict myself to the analysis of subzone
$\cZ_\out \cap \{r \rho \ge C\mu^{-1}\gamma^{-1} \}$.

Now one can see that in $\cZ_\out$ function $v(\rho;x,\mu^{-1})$ defined by (\ref{1-41}) has no more than one root $\rho \ge C\mu^{-1}r^{-1} \gamma^{-1}$; I denote it by $w=w(x,\mu^{-1})$. Otherwise (if such root does not exist) I set $w=0$. One can see easily that
\begin{equation}
v(\rho;x,\mu^{-1}) \asymp \mu^{-1}\gamma^{-2}(\rho^2 -w^2).
\label{1-64}
\end{equation}

If $\rho -w \asymp \Delta$ then the propagation speed for corrected symbol $\theta$ is $v(\rho,.)$; then the linear shift for time $\epsilon \mu^{-1}\gamma^{-1}$ is
$\asymp r \mu^{-2}\gamma^{-3}\rho\Delta $ (where $r$ appears because $\theta$ is an angle). This shift is observable if the logarithmic uncertainty principle
$r \mu^{-2}\gamma^{-3}\rho\Delta \times \Delta \ge Ch|\log h|$ holds\footnote{\label{foot-8} Now $\Delta$ is the scale in $\rho$.} i.e.
\begin{equation}
\Delta \ge C\mu \bigl( h|\log h| \bigr)^{1/2} \rho^{-1/2}\gamma^{3/2}r^{-1}.
\label{1-65}
\end{equation}
Since the propagation speed of symbol $|Z_1|^2$ corrected is $O(\mu^{-1}\gamma^{-2}\rho^2)$, the magnitude of $\Delta$ is preserved during the time interval $T''=\epsilon \mu \gamma^2 \rho^{-1}\Delta$ which is larger than $C\mu^{-1}\gamma^{-1}$ as
\begin{equation}
\Delta \ge C\mu^{-2}\rho^{-1}\gamma^{-3}.
\label{1-66}
\end{equation}
This condition is stronger than (\ref{1-65}) as $r\ge \mu^{-1}$ and (\ref{1-45}) is fulfilled with small enough $\delta>0$.

To justify this assertion properly one needs to operate with the quantizable symbols and one can see easily that this is the case under assumptions (\ref{1-45}) with small enough $\delta>0$ and (\ref{1-65}).
Then the contribution to the remainder estimate of the corresponding subzone (where magnitudes of $\gamma$, $\rho$, $r$ and $\Delta$ are fixed) does not exceed $C r^2h^{-3}\rho\gamma\Delta /T''=C\mu^{-1}h^{-3}r^2\rho^2 \gamma^{-1}$.

Let us note that as magnitudes of $\gamma,r,\rho$ are fixed and $\Delta$ ranges from $C\mu^{-2}\rho^{-1}\gamma^{-3}$ dictated by (\ref{1-66}) to $\epsilon\rho$, the number of such elements does not exceed $C\log \mu$. So the summation with respect to this zone results in $C\mu^{-1}h^{-3}r^2\rho^2 \gamma^{-1}\log \mu$. Further,
as $\Delta \ge \epsilon\rho$, there are no more than $C$ such elements and therefore the contribution to the remainder estimate of the part of the zone
\begin{equation}
\cY_{\gamma,\rho, r}\Def
\bigl\{{\frac 1 2}\gamma \le |x_1|\le 2\gamma,\ {\frac 1 2}\rho\le |Z_1|\le 2\rho,\ {\frac 1 2}r\le |x''|\le 2r \bigr\}
\label{1-67}
\end{equation}
where condition (\ref{1-66}) is fulfilled does not exceed
$C\mu^{-1}h^{-3}r^2\rho^2 \gamma^{-1}\log \mu$ as well. Summation with respect to $\rho$ with $\rho \le C(\gamma/r)^{1/2}$ results in
$C\mu^{-1}h^{-3}r\log \mu$ and then the summation with respect to all $r, \gamma$ results in $\mu^{-1}h^{-3} |\log \mu|^2=o(\mu^{-1/2}h^{-3})$. So,
contribution to the remainder estimate of the part of the zone
$\cZ_\out$ where $r\ge \mu^{-1}\gamma^{-1}$, $\rho^2\le C\gamma$ and condition (\ref{1-66}) is fulfilled is $O(\mu^{-1/2}h^{-3})$.

The same arguments one can apply as $\rho^2 \ge C\gamma/r$ but in this case $\Delta =\epsilon \rho$ and condition (\ref{1-64}) is fulfilled automatically, there will be no factor $\log \mu$ and the summation of $C\mu^{-1}h^{-3}r^2\rho^2 \gamma^{-1}$
with respect to $\rho$, $r$, $\gamma$ results in $C\mu^{-1/2}h^{-3}$.

Therefore I arrive to

\begin{claim}\label{1-68}
Contribution to the remainder estimate of the subzone of $\cZ_\out$ where condition (\ref{1-66}) is fulfilled is $O(\mu^{-1/2}h^{-3})$.
\end{claim}

Finally, contribution to the remainder estimate of the subzone $\cY_{\gamma,\rho}$ with $\gamma \ge C\mu^{-1/2}$ where condition (\ref{1-66}) is violated does not exceed
\begin{equation*}
C\bigl( h^{-3}\rho \gamma\times \mu^{-2}\rho^{-1}\gamma^{-3}\times T_1^{-1}\bigr)= C\mu^{-1}\gamma^{-1}h^{-3}
\end{equation*}
since $T_1=\epsilon \mu^{-1}\gamma^{-1}$ and summation with respect to $\gamma$ results again in $O(\mu^{-1/2}h^{-3})$.

Combining with (\ref{1-68}) I obtain immediately

\begin{proposition} \label{prop-1-16} Under condition $(\ref{1-44})$
the contribution to the remainder estimate of the zone $\cZ_\out$ is $O\bigl(\mu^{-1/2}h^{-3}\bigr)$.
\end{proposition}

\subsection{Quantum dynamics. II. Inner zone}
\label{sect-1-4}

Now I need to analyze the contribution to the remainder estimate of the \emph{inner zone} $\cZ_\inn =\{|x_1|\le {\bar\gamma}_0 = C\mu^{-1/2}\}$ under condition (\ref{1-45}). Both classical and quantum dynamics are confined to this zone as $T\le \epsilon$ and it starts in $B(0,{\frac 1 2})$; propagation speed with respect to $x$ does not exceed $C_0$. Note that the correction procedure still works as long as one eliminates only unbalanced factors $Z_2$, $Z_2^\dag$.

\subsubsection{}
\label{sect-1-4-1} Contribution to the remainder estimate of subzone
$\{|Z_1|\le \rho, |x_1|\le \gamma \}$ does not exceed
$CT^{-1} h^{-3}\rho^2r^2\gamma = C\mu h^{-3} \gamma \rho^2$ with $T=T_0=\epsilon \mu^{-1}$; plugging $\gamma = c\mu^{-1/2}$ and
$\rho = c\mu^{-1/2}$ I get $O(\mu^{-1/2}h^{-3})$ and therefore
one needs to consider only the contribution to the remainder estimate of the subzone $\cZ_\inn \cap \{|Z_1|\ge C \mu^{-1/2}\}$; but then proposition \ref{1-12} remains true here and one can increase $T=T_0$ to $T=T_1= \epsilon \mu^{-1/2}$\,\footnote{\label{foot-9} As in footnote $^{\ref{foot-6}}$ one must take $T_1=\epsilon \min (\mu^{-1/2},\rho)$ but contribution to the remainder estimate of subzone $\cZ_\inn\cap \{|Z_1|\le \mu^{-1/2}\}$ is $O(\mu h^{-3}\rho^2\gamma)=O(\mu ^{-1/2}h^{-3})$.}.

Then the contribution to the remainder estimate of zone
$\{|Z_1|\le \rho, |x_1|\le \gamma , |x''|\le r\}$ does not exceed
$CT ^{-1} h^{-3}\rho^2r^2\gamma = C\mu ^{1/2} h^{-3} \gamma r^2\rho^2$ with $T=\epsilon \mu^{-1/2}$; in particular, the total contribution of the whole inner zone $\cZ_\inn$ is $O(h^{-3})$ and one needs to recover factor $\mu^{-1/2}$.

To recover this factor one needs just to increase $T=T_1=\epsilon\mu^{-1/2}$ to $T=T_2=\epsilon$. Before doing this just note that the contributions to the remainder estimates of subzones
$\cZ_\inn \cap\{|Z_1|\le c\mu^{-1/4}\}$,
$\cZ_\inn \cap\{|x''|\le c\mu^{-1/4}\}$ and
$\cZ_\inn \cap\{|Z_1|\cdot |x''| \le c\mu^{-1/4}|\log \mu|^{-1}\}$
are $O(\mu^{-1/2}h^{-3})$ and therefore

\begin{claim}\label{1-69} One needs to consider only the contribution to the remainder estimate of the subzone
$\cZ_\inn \cap \{|x''|\ge c_1 \mu^{-1/4},\; |Z_1| \ge c_1 \mu^{-1/4},\; |Z_1|\cdot |x''| \ge c_1\mu^{-1/4}|\log \mu|^{-1} \}$.
\end{claim}

I remind that the classical dynamics which starts at
$\{|x_1|\le \gamma, |x''|=r \ge C_1\mu^{-1/4}\}$ remains in
$\{|x_1|\le C\gamma, \bigl||x''|-r\bigr| \le C_1\mu^{-1/2}r^{-1}\}$
as $T=\epsilon $; actually it could be larger as $\rho r \ll 1$.

Repeating arguments of the proof of proposition \ref{prop-1-14} one can justify it easily for quantum propagation as well:

\begin{proposition}\label{prop-1-17} Under condition $(\ref{1-45})$ with small enough $\delta>0$ the quantum dynamics which starts at
$\{|x_1|\le \gamma, |x''|=r \ge C_1\mu^{-1/4}\}$ remains in
$\{|x_1|\le C\gamma, \bigl||x''|-r\bigr| \le C_1\mu^{-1/2}r^{-1}\}$
as $T=\epsilon $.
\end{proposition}

Now let us consider again function ${\frac 1 2}x_1^2$; then
\begin{equation}
{\frac 1 2}\{a , x_1^2\} =
\Re \Bigl(\{Z_1^\dag, x_1\} x_1Z_1+\{Z_2^\dag, x_1\}x_1 Z_2\Bigr)
\label{1-70}
\end{equation}
and thus should be corrected by
\begin{equation}
\mu^{-1}\Re i\Bigl( \{Z_1^\dag,x_1\} f_1^{-1}x_1Z_1 +\{Z_2^\dag,x_1\} f_2^{-1}x_1Z_2\Bigr)
\label{1-71}
\end{equation}
which is $O(\mu^{-1}r)$ as $r\ge \mu^{-1/2}$ and leads to the new error
\begin{align}
\mu^{-1}\Re i\Bigl(
&\bigl\{a,\{Z_1^\dag,x_1\} f_1^{-1}x_1 \bigr\} Z_1 +
\{Z_1^\dag,x_1\} f_1^{-1}x_1 \bigl( \{a, Z_1\} +i\mu Z_1\bigr)+\label{1-72}\\
&\bigl\{a,\{Z_2^\dag,x_1\} f_2^{-1}x_1 \bigr\} Z_2 +
\{Z_2^\dag,x_1\} f_2^{-1}x_1 \bigl( \{a, Z_2\} +i\mu Z_2\bigr)
\Bigr).\notag
\end{align}
This new error is $O(\mu^{-1})$ and thus taking $T=\epsilon r$ one can see that the oscillation of $x_1^2$ would be $O(\mu^{-1}r)$ which leads to

\begin{proposition}\label{prop-1-18} (i) The classical dynamics which starts at
$\{|x_1| = \gamma, |x''|=r \ge C_1\mu^{-1/4}\}$ with
$C_0 \mu^{-1/2}r^{1/2}\le \gamma \le C_0\mu^{-1/2} $ remains in $\{C^{-1}\gamma \le |x_1|\le C\gamma\}$
as $T=\epsilon r$.

\smallskip
\noindent
(ii) The classical dynamics which starts at
$\{|x_1| \le \gamma, |x''|=r \ge C_1\mu^{-1/4}\}$ with
$\gamma = C_0\mu^{-1/2}r^{1/2} $ remains in $\{ |x_1|\le C\gamma\}$
as $T=\epsilon r$.

\smallskip
\noindent
(iii) Under condition $(\ref{1-45})$ with small enough $\delta>0$ statements (i),(ii) remain true for a quantum dynamics as well.
\end{proposition}

So, one can see that the dynamics in the zones described in proposition \ref{prop-1-18}(i), (ii) are different with the dynamics in the former resembling the dynamics in $\cZ_\out$.
Then the arguments of the previous subsection work perfectly in the zone described in proposition \ref{prop-1-18}(i); tedious but easy\footnote{\label{foot-10} Because one needs to worry only about eliminating factor $\mu^{1/2}$ while in the previous subsection the offending factor was $\mu^{3/2}$. One can just take $\Delta = C\mu^{-1/2}$. Also one can assume that $|x_1|\ge {\hat\gamma} = c_1\mu^{-3/5}$; otherwise $|x''|\le r= c_2\mu ^{-1/5}$ and the contribution to the remainder estimate of this zone would be $O(\mu h^{-3}\gamma^2 r^2)= O(\mu^{-1/2}h^{-3})$.

But then $\rho \ge \mu^{-1}\gamma^{-1}$ in the zone of interest and one can finally upgrade $T_1=\epsilon\mu^{-1/2}$ to $T_1=\epsilon \mu^{-1}\gamma^{-1}$ and the crude estimate to $O(\mu h^{-3}\rho^2\gamma^2r^2)$.} details I leave to the reader:

\begin{proposition}\label{prop-1-19} Under condition $(\ref{1-45})$ with sufficiently small $\delta>0$ the contribution to the remainder estimate of the zone $\{ C_0\mu^{-1/2}|x''|^{1/2} \le |x_1| \le C_0\mu^{-1/2}\}$ does not exceed $C\mu^{-1/2}h^{-3}$.
\end{proposition}

Thus I am left with the zone
\begin{equation}
\bigl\{ C_1\mu^{-1/6}\le |x''|\le c, \quad |x_1|\le C_1\mu^{-1/2}|x''|^{1/2} \bigr\}.
\label{1-73}
\end{equation}

Again the dynamics in the zones
$\cY_{\gamma, \rho, r}$ with $\rho r \le \epsilon \mu \gamma^2$
and $\rho r \ge \epsilon \mu \gamma^2$
are very different with the former one more similar to the dynamics in the outer zone (as long as it remains there).

Similarly the propagation speed of the symbol $|Z_1|^2$ corrected by $O(\mu^{-1} \gamma^{-1}r |Z_1|+\mu^{-1})$ does not exceed $C\bigl(\mu^{-1}\gamma^{-2}r|Z_1|+\mu^{-1}\bigr)$ and therefore the propagation speed of $|Z_1|$ does not exceed $C\bigl(\mu^{-1}\gamma^{-2}r+\mu^{-1/2}\bigr)$
and thus magnitude of $|Z_1|\asymp \rho$ is preserved on the time interval $T=\epsilon\bigl(\mu^{-1}\gamma^{-2}r+\mu^{-1/2}\bigr)^{-1}\rho $ as long as $\gamma \ge C(\mu r \rho)^{1/2}$. As a result I arrive to

\begin{proposition}\label{prop-1-20} (i) The classical dynamics which starts at
$\{|x_1| = \gamma, |x''|=r \ge C_1\mu^{-1/4}, |Z_1|=\rho \}$ with
$C_0 \mu^{-1/2}r^{1/2}\rho^{1/2} \le \gamma \le C_0\mu^{-1/2} $ remains in $\{C^{-1}\rho \le |Z_1|\le C\rho \}$
as $T=\epsilon \min\Bigl( r, \mu \gamma r \rho ^{-1}, \mu^{-1}\rho^{-1}\Bigr)$.

\smallskip
\noindent
(ii) The classical dynamics which starts at
$\{|x_1| \le \gamma, |x''|=r \ge C_1\mu^{-1/4}, |Z_1|=\rho \}$, with $\gamma = C_0\mu^{-1/2}r^{1/2}\rho^{1/2} $ remains in
$\{C^{-1}\rho \le |Z_1|\le C\rho \}$
as $T=\epsilon \min\Bigl( r, \mu \gamma r \rho ^{-1}, \mu^{-1}\rho^{-1}\Bigr)$.

\smallskip
\noindent
(iii) Under condition $(\ref{1-45})$ with small enough $\delta>0$ statements (i),(ii) remain true for a quantum dynamics as well.
\end{proposition}

\subsubsection{}
\label{sect-1-4-2} Then again arguments of the previous subsection work perfectly in the zone described in proposition \ref{prop-1-20}(i); tedious but easy details (see footnote $^{\ref{foot-9}}$ with the obvious modifications.) I leave to the reader:

\begin{proposition}\label{prop-1-21} Under condition $(\ref{1-44})$ with sufficiently small $\delta>0$ the contribution to the remainder estimate of the zone
$\{ C_0\mu^{-1/2}|x''|^{1/2}|Z_1|^{1/2} \le |x_1| \le C_0\mu^{-1/2}\}$ does not exceed $C\mu^{-1/2}h^{-3}$.
\end{proposition}

Therefore I am left with the \emph{true inner zone}
\begin{equation}
\cZ_\inn^0 \Def\bigl\{ |x_1|\le C_1\mu^{-1/2}|x''|^{1/2}|Z_1|^{1/2} \bigr\}.
\label{1-74}
\end{equation}
Similarly to (\ref{1-69})

\begin{claim}\label{1-75} One needs to consider only the contribution to the remainder estimate of the subzone
$\cZ^0_\inn \cap \{|x''|\ge c_1 \mu^{-1/6},\; |Z_1| \ge c_1 \mu^{-1/6}, \;|Z_1|\cdot |x''| \ge c_1\mu^{-1/6}|\log \mu|^{-1} \}$.
\end{claim}

Let us consider classical dynamics in this zone first; I am interested in the time interval $T\le \epsilon r$. There
\begin{equation*}
\{a,Z_1\}= \{Z_1^\dag,Z_1\}Z_1 + \{Z_2^\dag,Z_1\}Z_2+\{Z_2,Z_1\}Z_2^\dag-\{V,Z_1\}
\end{equation*}
and then I correct $Z_1$ to
\begin{equation}
\zeta_1 = Z_1 + i\mu ^{-1}f_2^{-1} \{Z_2^\dag,Z_1\}Z_2 -i\mu ^{-1}f_2^{-1} \{Z_2,Z_1\}Z_2^\dag
\label{1-76}
\end{equation}
satisfying
\begin{multline*}
\{a,\zeta_1\}\equiv \{Z_1^\dag,Z_1\}Z_1 -if_2^{-1}\mu^{-1}\Bigl( \bigl\{Z_2,\{Z_2^\dag, Z_1\}\bigr\}- \bigl\{Z_2^\dag ,\{Z_2, Z_1\}\bigr\} \Bigr)|Z_2|^2-\{V,Z_1\}=\\
\{Z_1^\dag,Z_1\}Z_1 +if_2^{-1}\mu^{-1} \bigl\{Z_1, \{Z_2,Z_2^\dag\}\bigr\} |Z_2|^2-\{V,Z_1\}\equiv \\
\{Z_1^\dag,Z_1\}Z_1 - f_2^{-1}\{Z_1,f_2\}|Z_2|^2 -\{V,Z_1\}
\qquad\mod O(\mu^{-1/2}).
\end{multline*}
Since I do not mind to change the speed to one of the same magnitude and I consider only level 0, I can assume with no loss of the generality that
\begin{equation}
f_2=1;
\label{1-77}
\end{equation}
then finally $\{a,\zeta_1\}\equiv \{\zeta_1^\dag \zeta_1 -V,\zeta_1\}$ $\mod O(\mu^{-1/2})$.

Also corrected functions $x_j$ satisfy
$\{a,x'_j\}\equiv \{\zeta_1^\dag \zeta_1 -V,x'_j\}$
$\mod O(\mu^{-1/2})$ and I arrive to the following conclusion:

\begin{claim}\label{1-78} Under condition (\ref{1-77}) the classical dynamics in zone $\{|x_1|\le C\mu^{-1/2}\}$ in variables $x_j,Z_1,Z_2$ is described modulo $O(\mu^{-1}+\mu^{-1/2}T)$ by the solution to the short system
\begin{align}
&{\frac {d Z_1}{dt}}= {\frac 1 2}\{Z_1^\dag Z_1 -V,Z_1\}\label{1-79},\\
&{\frac {dx_j}{dt}}= {\frac 1 2}\{Z_1^\dag Z_1 -V,x_j\}\label{1-80},\\
&\{Z_1,Z_2\}=O(\mu^{-1}),\quad \{Z_1,Z_2^\dag\}=O(\mu^{-1}).\label{1-81}
\end{align}
\end{claim}

\vskip-20pt
Further, for each given initial point one can replace all functions modulo
$(x_1-{\frac 1 2} r^2)$, $(x_2-\theta r^2)$ so that I get two-dimensional system in variables $(r,\theta)$ described by the Hamiltonian ${\frac 1 2}(Z_1^\dag Z_1-V)$ with $Z_1,Z_1^\dag$ satisfying commutator relation
\begin{equation}
\{Z_1,Z_1^\dag\}=\mu \phi {\bar r} (r-{\bar r})+O(\mu^{-1}))
\label{1-82}
\end{equation}
which is exactly system studied extensively in \cite{IRO6} with magnetic field
$ \phi (r-{\bar r})$ and coupling constant ${\bar\mu}=\mu{\bar r}$. Condition
(\ref{1-74}) means exactly that I am in the inner zone
$\{|r-{\bar r}|\le C {\bar\mu}^{-1/2}|Z_1|^{1/2}\}$ for this system; since (in contrast to \cite{IRO6}) $|Z_1|$ is not disjoint from 0 anymore one needs to define inner zone in this way.

Then the drift speed with respect to $\theta r$ is
\begin{equation}
\asymp \Bigl( \eta - k^*\phi^{-1/2} \rho \Bigr) \rho r
\label{1-83}
\end{equation}
where constant $k^*\approx 0.66$ is defined in \cite{IRO6} and $\eta$ is the corrected symbol $\mu x_1^2/2$ which leads to its propagation speed given by (\ref{1-72}) multiplied by $\mu$.

One can see easily that the fourth term there
$\{Z_2^\dag,x_1\} f_2^{-1}x_1 \bigl( \{a, Z_2\} +i\mu Z_2\bigr)$ is $O(\mu^{-1/2})$. Further, the third term $ \bigl\{a,\{Z_2^\dag,x_1\} f_2^{-1}x_1 \bigr\} Z_2$ could be rewritten modulo $O(\mu^{-1/2})$ and terms with unbalanced $Z_2,Z_2^\dag$ as $|\{Z_2,x_1\}|^2 f_2^{-1}|Z_2|^2$ and disappears after applying $\Re i$.

So I am left with just two terms
\begin{equation}
\Re i\Bigl(
\bigl\{a,\{Z_1^\dag,x_1\} \phi^{-1} \bigr\} Z_1 +
\{Z_1^\dag,x_1\} \phi^{-1} \bigl( \{a, Z_1\} +i\mu Z_1\bigr)
\Bigr)\label{1-84}
\end{equation}
where $\phi= f_1/x_1$ and one can replace $a$ by $Z_1^\dag Z_1-V$ (other terms have unbalanced $Z_2,Z_2^\dag$).

\subsubsection{}
\label{sect-1-4-3} So I am \emph{exactly} in \cite{IRO6} situation. I remind that in \cite{IRO6} symbol $Z_1$ could be reduced to
\begin{equation}
Z_1=e^{i\beta (x)} \Bigl(\xi_1 + i\alpha (x) \bigl(\xi_2- V_2(x)\bigr)\Bigr)
\label{1-85}
\end{equation}
with real-valued $\alpha,\beta, V_2$ such that $\alpha=1+O(x_1)$,
$V_2=x_1^2/2 +O(x_1^3)$.

One can see easily that in this case that corrected $\mu x_1^2$ would be $\xi_2+O(\mu^{-1/2})$; I want to remind that in \cite{IRO6} exactly deviation of $\xi_2$ from $k^*\phi ^{-1/2}\rho$ was used as the measure of the drift with $\rho = V^{1/2}$. Note that the potential $V$ now should be replaced by $(V-E_2)$ with $E_2 =|Z_2^2|$.

Therefore due to the logarithmic uncertainty principle the violation of periodicity is observable after the first turn as
\begin{equation}
\Bigl( \eta - k^*\phi^{-1/2} \rho \Bigr) \rho r \mu^{-1/2} \times \Bigl( \eta - k^*\phi^{-1/2} \rho \Bigr) \ge Ch|\log h|;
\label{1-86}
\end{equation}
under assumption (\ref{1-45}) the latter condition is the automatic corollary of
\begin{equation}
\Delta=| \eta - k^*\phi^{-1/2} \rho |\ge C\mu^{-1/2}r^{-1/2}\rho^{-1/2}.
\label{1-87}
\end{equation}
Under this condition the drift speed for the ``short'' system is greater than $C\epsilon_1\mu^{-1/2}$ which is the larger than the error in the drift speed
$c \mu^{-1/2}$ (see (\ref{1-78}).

In this case $T_1=\epsilon r \Delta $ and the rest is easy but extra logarithmic factor would appear; however considering the right direction one can take
$T_1=\epsilon r \rho^{1/2} \Delta^{1/2}$ and following \cite{IRO6} with the standard justification on the quantum level one can prove easily that

\begin{claim}\label{1-88}
The contribution of the zone where condition (\ref{1-87}) is fulfilled to the remainder estimate is $O(\mu^{-1/2}h^{-3})$.
\end{claim}

On the other hand, contribution to the remainder estimate of zone where this condition is violated would not exceed
\begin{equation*}
C\Bigl(\int \mu^{-1/2}h^{-3}r^{-1/2}\rho ^{-1/2} rdr d\rho\Bigr)=O(\mu^{-1/2}h^{-3})
\end{equation*}
as well.

So, the main result of this section is proven:

\begin{proposition}\label{prop-1-22} Under condition $(\ref{1-45})$ with small enough $\delta>0$ the remainder estimate is $O(\mu^{-1/2}h^{-3})$ while the principal part is given by the standard Weyl formula (which is equal modulo the remainder estimate to magnetic Weyl formula).
\end{proposition}

This statement implies trivially

\begin{corollary}\label{cor-1-23} Under conditions $(\ref{0-2})$ and $(\ref{1-45})$ estimate $(\ref{0-3})$ holds.
\end{corollary}

\section{Strong magnetic field. Canonical forms}
\label{sect-2}

Now I want to consider the main case
\begin{equation}
h^{-\delta}\le \mu \le Ch^{-1}
\label{2-1}
\end{equation}
with an arbitrarily small exponent $\delta >0$. In this section I consider different canonical forms of the Magnetic Schr\"odinger Operator in question (depending on the zone).

\subsection{Precanonical Form}
\label{sect-2-1}
\subsubsection{}
\label{sect-2-1-1} Let
\begin{equation}
U_2(x,\xi)=f_2(x)^{-1/2}Z_2(x,\xi)
\label{2-2}
\end{equation}
Then
\begin{phantomequation}\label{2-3}\end{phantomequation}
\begin{equation}
-i\mu^{-1}\{U_2,U_2^\dag\}\equiv 2\qquad \mod \cO_k
\tag*{$(2.3)_k$}
\label{2-3-k}
\end{equation}
with $k=1$ where $\cO_k=\cO_k(\mu^{-1} Z_1, \mu^{-1} Z_1^\dag, \mu^{-1} Z_2, \mu^{-1} Z_2^\dag)$ is the space of the sum of polynomials containing monoms of order $k$ or higher with respect to $(\mu^{-1} Z_1, \mu^{-1} Z_1^\dag, \mu^{-1} Z_2, \mu^{-1} Z_2^\dag)$ with the coefficients, smoothly depending on $x$.

Note that correcting $U_2$ by $U_2\cO_k$ one trades $(\ref{2-3})_k$ by $(\ref{2-3})_{k+1}$ and therefore after an appropriate correction (\ref{2-3}) holds with arbitrarily large $k=M$.

So, correcting $U_2$:
\begin{equation}
U_2 \mapsto U_2 + \sum_{\substack{k+l+p+q\ge 2\\k\ge 1} } \alpha_{klpq}(x)U_2^kU_2^{\dag\,l}Z_1^pZ_1^{\dag\,q}\mu^{1-k-l-p-q}
\label{2-4}
\end{equation}
one one can achieve
\begin{equation}
-i\mu^{-1}\{U_2,U_2^\dag\}\equiv 2\qquad\mod O(\mu^{-M})
\label{2-5}
\end{equation}
with arbitrarily large $M$.

Note that $\{U_2,Z_1\}\equiv \{U_2^\dag ,Z_1\}\equiv 0$ modulo $\cO_1$ both originally and after correction (\ref{2-4}).

I claim that

\begin{proposition}
\label{prop-2-1}
Correcting
\begin{align}
&Z_1 \mapsto Z_1 + \sum_{\substack{k+l+p+q\ge 2\\k+l\ge 1 } }
\alpha'_{klpq}(x) U_2^kU_2^{\dag\,l}Z_1^pZ_1^{\dag\,q}\mu^{1-k-l-p-q},
\label{2-6}\\
&x_j \mapsto x_j + \sum_{k+l\ge 1 }
\alpha''_{jklpq}(x) U_2^kU_2^{\dag\,l}Z_1^pZ_1^{\dag\,q}\mu^{-k-l-p-q},\label{2-7}
\end{align}
one can arrange
\begin{align}
&\{U_2, x_j\}\equiv \{U_2,Z_1\}\equiv \{U_2,Z_1^\dag\}\equiv 0,\label{2-8}\\
&\{U_2^\dag, x_j\}\equiv \{U_2^\dag,Z_1\}\equiv \{U_2^\dag,Z_1^\dag\}\equiv 0\qquad \mod O(\mu^{-M}).
\label{2-9}
\end{align}
\end{proposition}

\begin{proof}
Really, correcting in the way described (but with the sum over $k\ge 1$ rather than $k+l\ge1$) I can achieve (\ref{2-9}). However, \emph{so far corrected symbol $x_j$ is not necessarily real-valued and corrected symbols $Z_1$ and $Z_1^\dag$ are not necessarily complex conjugate}.

Then due to (\ref{2-9}), (\ref{2-5}) and Poisson identity
$\{U_2^\dag, \{U_2,y\}\}\equiv 0$ with $y=x_j, Z_1,Z_1^\dag$. Therefore
\begin{align*}
&\{U_2,Z_1\}\equiv\sum_{l+p+q\ge 1} \beta_{lpq}(x)U_2^{\dag\,l}Z_1^pZ_1^{\dag\,q}\mu^{1-l-p-q},\\
&\{U_2,Z_1^\dag\}\equiv\sum_{l+p+q\ge 1} \beta'_{lpq}(x)U_2^{\dag\,l}Z_1^pZ_1^{\dag\,q}\mu^{1-l-p-q},\\
&\{U_2,x_j\}\equiv\sum_{l+p+q\ge 1} \beta''_{jlpq}(x)U_2^{\dag\,l}Z_1^pZ_1^{\dag\,q}\mu^{-l-p-q}
\end{align*}
where $x$, $Z_1$ and $Z_1^\dag$ are already corrected symbols. Since at this moment $x_j$ are not necessarily real-valued one needs to plug them formally into functions and use a Taylor decomposition with respect to $\mu^{-1}Z_1,\mu^{-1}Z_2,\mu^{-1}U_1,\mu^{-1}U_2$.

Then one can correct $x_j$, $Z_1$ and $Z_1^\dag$ according to formulae (\ref{2-7}),(\ref{2-8}) with summation $k=0$ so that (\ref{2-8}) holds.

Note that instead of $Z_1$ and $Z_1^\dag$ one could consider $\Re Z_1$ and $\Im Z_1$ and correct them deriving (\ref{2-8})-(\ref{2-9}). At this moment corrected symbols $x_j$, $\Re Z_1$ and $\Im Z_1$ are not necessarily real-valued. However
$U_2$ and $U_2^\dag$ are truly complex conjugate and then
(\ref{2-8})-(\ref{2-9}) for corrected symbols $x_j$, $\Re Z_1$, $\Im Z_1$ imply the same equalities for the real parts of them. Let us replace then (corrected) symbols $x_j$, $\Re Z_1$ and $\Im Z_1$ by their real parts. After this (corrected) symbols $x_j$ become real-valued and (corrected) symbols $Z_1=\Re Z_1+i\Im Z_1$ and $Z_1^\dag=\Re Z_1-i\Im Z_1$ become complex conjugate.
\end{proof}

\subsubsection{}
\label{sect-2-1-2} One can rewrite operator in question in these new variables as
\begin{multline}
2 A\equiv f_2(x) U_2^\dag U_2 + Z_1^\dag Z_1 -V(x)+\\
\\
\sum_{k+l+p+q\ge 3} b_{klpq}(x)U_2^kU_2^{\dag\,l}Z_1^pZ_1^{\dag\,q}\mu^{2-k-l-p-q}+
\sum_{k+l+p+q\ge 3} b'_{klpq}(x)U_2^kU_2^{\dag\,l}Z_1^pZ_1^{\dag\,q}\mu^{-k-l-p-q}
\label{2-10}
\end{multline}
modulo lower order terms where the first line is the ``main part'' and the second line is the ``perturbation''; here the first part of the ``perturbation'' comes from $f_2U_2^\dag U_2 +Z_1^\dag Z_1$ and the second part comes from $-V$.

\begin{proposition}\label{prop-2-2} For arbitrarily large $M$ there exist $\epsilon=\epsilon (M)>0$ and a real valued symbol
\begin{equation}
\cL = \sum_{k+l \ge 1} \ell_{kpqm}(x)U_2^kU_2^{\dag\,l}Z_1^pZ_1^{\dag\,q}\mu^{-1-k-l-p-q-m}
\label{2-11}
\end{equation}
such that in the strip $\{|x_1|<\epsilon\}$ the following equality holds (so far only for the principal symbols of operators in question):
\begin{multline}
A^\#\Def e^{-i\mu^{-1}h^{-1}\cL} Ae^{i\mu^{-1}h^{-1}\cL}\equiv
{\frac1 2}\Bigl(f_2(x^\#) U_2^\# U_2^{\#\dag}+ Z_1^\# Z_1^{\#\dag} -V(x^\#)\Bigr)+\\
\sum_{2k+p+q+2m\ge 3}^{} B_{kpqm}(x^\#)\bigl(U_2^\#U_2^{\#\dag}\bigr)^k Z_1^{\#\,p} Z_1^{\#\dag\,q}\mu^{2-2k-p-q-2m}
\label{2-12}
\end{multline}
where $e^{-i\mu^{-1}h^{-1}\cL}$, $e^{i\mu^{-1}h^{-1}\cL}$ are (formal) $\mu^{-1}h$-FIOs and
\begin{multline}
x_j^\# = e^{-i\mu^{-1}h^{-1}\cL}x_je^{i\mu^{-1}h^{-1}\cL},\quad
U_2^\# = e^{-i\mu^{-1}h^{-1}\cL}U_2e^{i\mu^{-1}h^{-1}\cL},\\
Z_1^\# = e^{-i\mu^{-1}h^{-1}\cL}Z_1 e^{i\mu^{-1}h^{-1}\cL}
\label{2-13}
\end{multline}
are linked to the original symbols $x_j$, $U_2$, $Z_1$ by formulae similar to $(\ref{2-4}),(\ref{2-6}),(\ref{2-7})$ but containing also factors $\mu^{-1}$
exactly on the same role as $U_2$ or $U_2^\dag$ so that $k+l$ in the condition of summation is replaced by $k+l+m$\,\footnote{\label{foot-11} More precisely
\begin{align}
&U_2 \mapsto U_2^\#=U_2 + \sum_{k+l+p+q+m\ge 2 } \alpha_{klpqm}(x)U_2^kU_2^{\dag\,l}Z_1^pZ_1^{\dag\,q}\mu^{1-k-l-p-q-m},
\tag*{$(2.4)^*$}\label{2-4-*}\\
&Z_1 \mapsto Z_1^\#=Z_1 + \sum_{k+l+p+q+m\ge 2 }
\alpha'_{klpqm}(x) U_2^kU_2^{\dag\,l}Z_1^pZ_1^{\dag\,q}\mu^{1-k-l-p-q-m},
\tag*{$(2.6)^*$}\label{2-6-*}\\
&x_j \mapsto x^\#_j=x_j + \sum_{k+l+p+q+m\ge 1 }
\alpha''_{jklpqm}(x) U_2^kU_2^{\dag\,l}Z_1^pZ_1^{\dag\,q}\mu^{-k-l-p-q-m}.
\tag*{$(2.7)^*$}\label{2-7-*}
\end{align}
}.

All equalities here hold modulo $O(\mu^{-M})$.
\end{proposition}

\begin{proof} Proof is standard (see \cite{IRO5}). The problem of resonance is avoided because I remove only terms $b_* (x)U_2^kU_2^{\dag\,l}Z_1^pZ_1^{\dag\,q}$ with $k\ne l$ which leads to the denominator
\begin{equation*}
\mu \bigl((k-l)f_2 + (p-q)f_1\bigr)
\end{equation*}
in the correction term; this denominator is of magnitude $\mu (k-l)$ since $k\ne l$, $f_2\asymp 1$ and $f_1\asymp |x_1|\le \epsilon (M)$.
\end{proof}

Now one can upgrade all the above arguments to operators so that equalities (\ref{2-5}),(\ref{2-8}) and (\ref{2-9}) hold for commutators rather than for Poisson brackets and (\ref{2-12}) is replaced by
\begin{multline}
A^\#\Def e^{-i\mu^{-1}h^{-1}\cL} Ae^{i\mu^{-1}h^{-1}\cL}\equiv
{\frac 1 2}\Bigl(f_2^\w U_2^*\bullet U_2 + Z_1^*\bullet Z_1 -V^\w\Bigr)+\\
\sum_{2k+p+q+2m+2s\ge 3}^{} B_{kpqms}^\w (U_2\bullet U_2^*)^l Z_1^pZ_1^{*\,q}\mu^{2-2k-2m-p-q-s}h^s
\tag*{$(2.12)^*$}
\label{2-12-*}
\end{multline}
where $a^\w$ is the Weyl quantization of symbol $a$ and one should remember that the corrected symbol $x$ is not just a function of coordinates but symbols and therefore the functions of $x$ are replaced by the symbols of $\mu^{-1}h$-PDOs. Here and below $\bullet$ means the symmetrized product: $K\bullet L={\frac 1 2}(KL+LK)$. So, I arrive to

\begin{proposition}\label{prop-2-3} For arbitrarily large $M$ and arbitrarily small $\delta>0$ there exist $\epsilon=\epsilon (M,\delta)>0$ and a real valued symbol
\begin{equation}
\cL = \sum_{k+l+p+q+m+2s\ge 1} \ell_{kpqms}(x)U_2^kU_2^{\dag\,l}Z_1^pZ_1^{\dag\,q}\mu^{-1-k-l-p-q-m-s}h^s
\tag*{$(2.11)^*$}\label{2-11-*}
\end{equation}
such that under condition $(\ref{2-1})$ in the strip $\{|x_1|<\epsilon\}$ formula \ref{2-12-*} holds where $e^{-i\mu^{-1}h^{-1}\cL}$, $e^{i\mu^{-1}h^{-1}\cL}$ are $\mu^{-1}h$-FIOs and
$x_j^\#,U_2^\#,Z_1^\#$ are full symbols still defined by $(\ref{2-13})$ and they
are linked to the original symbols $x_j$, $Z_1$, $U_2$ by formulae similar to \ref{2-4-*},\ref{2-6-*},\ref{2-7-*} but containing also factors $(\mu^{-1}h)$
of the double value of $U_2$ or $U_2^\dag$ or $\mu^{-1}$ so that $k+l+m$ in the condition of summation is replaced by $k+l+m+2s$\,\footnote{\label{foot-12} More precisely
\begin{align}
&U_2 \mapsto U_2^\#=U_2 + \sum_{k+l+p+q+m+2s\ge 2 } \alpha_{klpqms}(x)U_2^kU_2^{\dag\,l}Z_1^pZ_1^{\dag\,q}\mu^{1-k-l-p-q-m-s}h^s,
\tag*{$(2.4)^{**}$}\label{2-4-**}\\
&Z_1 \mapsto Z_1^\#=Z_1 + \sum_{k+l+p+q+m+2s\ge 2 }
\alpha'_{klpqmm}(x) U_2^kU_2^{\dag\,l}Z_1^pZ_1^{\dag\,q}\mu^{1-k-l-p-q-m-s}h^s,
\tag*{$(2.6)^{**}$}\label{2-6-**}\\
&x_j \mapsto x^\#_j=x_j + \sum_{k+l+m+2s\ge 1 }
\alpha''_{jklpqms}(x) U_2^kU_2^{\dag\,l}Z_1^pZ_1^{\dag\,q}\mu^{-k-l-p-q-m-s}h^s,
\tag*{$(2.7)^{**}$}\label{2-7-**}
\end{align}
}.

All equalities here and below hold modulo $O(h^M)$.
\end{proposition}

\subsubsection{}
\label{sect-2-1-3}
Note that due to (\ref{2-5}) for $U_2^\#$ there exists $\mu^{-1}h$-FIO $\cF$ such that
\begin{equation}
\cF^*\cF\equiv \cF\cF^*\equiv I\qquad\text{and}\qquad
\cF^*U_2^\#\cF \equiv (hD_2 -i\mu x_2).
\label{2-14}
\end{equation}
Then
\begin{equation}
\cF^*(U_2^\#\bullet U_2^{\#\dag})\cF \equiv \bigl(h^2D_2^2+\mu^2 x_2^2\bigr)
\label{2-15}
\end{equation}
and also due to (\ref{2-8}),(\ref{2-9})
$\mu^{-1}\cF^*Z_1^\#\cF$, $\cF^*x_j^\#\cF$ and $\cF^* B_{kpqms}\cF$ are $\mu^{-1}h$-PDOs with symbols depending only on $x',\xi'$ ($x'=(x_1,x_3,x_4)$).
Let us redefine $Z_1^\# \redef \mu^{-1}\cF^*Z_1^\#\cF$,
$B^\w_{kpqms} \redef \cF^* B^\w_{kpqms}\cF$ and $A^\# \redef \cF^* A^\# \cF$; then also automatically
$ Z_1^{\# \dag} \redef\mu^{-1}\cF^*Z_1^{\#\dag}\cF$ is an adjoint operator. I remind that $f_2^\w$ and $B^\w_{kpqms}$ are operators with the real-valued symbols and $A^\#$ is again defined by \ref{2-12-*}.

Now one can decompose $u$ into series
\begin{gather}
u(x,y,t)=\sum_{n,n'\in \bZ^+} \mu h^{-1} u_{nn'}(x',y',t)\,
\Upsilon_n( x_2)
\Upsilon_{n'}( y_2), \quad
\Upsilon_n( x_2)\Def \upsilon_n\bigl(\mu^{1/2}h^{-1/2} x_2\bigr)
\label{2-16}
\end{gather}
where $\upsilon_n$ are (real-valued and orthonormal) Hermite functions and then replace
operator $\bigl(h^2D_2^2+\mu^2x_2^2\bigr)$ by $(2n+1)\mu h$, thus reducing
operator $A^\#$ to the family of 3-dimensional $\mu^{-1}h$-PDOs
\begin{multline}
\cA_n \Def {\color{blue}
{\frac 1 2}\Bigl(f_2^\# (2n+1) \mu h + (Z_1^\# \bullet Z_1^{\#\dag}) -V^\# \Bigr)}+\\
\sum_{2k+p+q+2m+2s\ge 1}^{} B^w_{kpqms}\times \bigl((2n+1)\mu h\bigr)^k
Z_1^{\#\,p}Z_1^{\#\dag\,q}\mu ^{2-2k-2m-p-q-s} h ^s;
\label{2-17}
\end{multline}
again the first line is the \emph{main part} of operator and $f^\#$, $V^\#$ are
transformed $f_2$ and $V$.

\begin{remark}\label{rem-2-4}
(i) In this and below formulae one actually needs to consider
\begin{equation}
n\le C_0/(\mu h)+C|\log h|.
\label{2-18}
\end{equation}
Really, I am interested in the domain in the phase space where $\{a\le c\}$ or, equivalently, $\{|U_2|\le c, |Z_1|\le c\}$. Note that
condition $\{|U_2|\le \varrho\}$ does not contradict to the logarithmic uncertainty principle as long as $\varrho^2 \ge C\mu h |\log h|$ and thus one can take an upper bound for $a$ as $\varrho^2 = c+ C\mu h |\log h|$ which implies (\ref{2-18}). This estimate perfectly suits my purposes as long as
\begin{equation}
h^{-\delta} \le \mu \le \epsilon h^{-1}|\log h|^{-1}.
\label{2-19}
\end{equation}
\noindent
(ii) From the operator point of view which one needs to apply only as
\begin{equation}
\epsilon h^{-1}|\log h|^{-1}\le \mu \le \epsilon h^{-1}
\label{2-20}
\end{equation}
in the final analysis I need to consider only $n\le C_0/(\mu h)$.
\end{remark}

\begin{remark}\label{rem-2-5}
(i) Obviously $f_j^\#$, $V^\#$ and $\mu^{-1}Z_1^\#$ are $\mu^{-1}h$-PDOs with symbols $f_j\circ \Psi_2$, $V\circ \Psi_2$ and $\mu ^{-1}Z_1\circ \Psi_2$
respectively where $\Psi_2:\bR^6 \ni (x',\xi')\to x\in \bR^4$ is some smooth map. Actually I am interested in the symbols of $f_2^\#$ and $V^\#$ only modulo
symbols of $\mu^{-1}Z^\#_1$, $\mu^{-1}Z_1^{\#\dag}$ and therefore actually I am interested only in the map $\Sigma \to \bR^4$. However at this moment I do not have a natural parametrization of $\Sigma$.

\smallskip
\noindent
(ii) Note that $f_1^\#$ is not necessarily $x_1$ anymore. Still since
$\nabla \Re Z_1$, $\nabla \Im Z_1$, $\nabla \Re Z_2$, $\nabla \Im Z_2$ and $\nabla f_1$ were linearly independent\footnote{\label{foot-13} All the linear independence statements are uniform.} this linear independence statement is true for $\nabla \Re Z_1^\#$, $\nabla \Im Z_1^\#$ and $\nabla f_1^\#$\,\footnote{\label{foot-14} Where in the former case it was $\nabla_{x,\xi}$ and in the latter one $\nabla_{x',\xi'}$; I remind that $x'=(x_1.x_3,x_4)$.} as well. Therefore after $\mu^{-1}h$-FIO transformation (in $(x',\mu^{-1}hD')$) one can achieve $f_1^\#=x_1$.
\end{remark}

So, I arrive to

\begin{proposition}\label{prop-2-6} Let condition $(\ref{2-1})$ be fulfilled. Then by means of $\mu^{-1}h$-transform, decomposition $(\ref{2-16})$ and one more $\mu^{-1}h$-transform one can reduce operator to the family of 3-dimensional operators $(\ref{2-17})$ with $n$ satisfying $(\ref{2-18})$ and with
\begin{equation}
i \mu^{-1}\{Z_1^\#,Z_1^{\#\dag}\}=f_1^\#=x_1\qquad \mod (Z_1,Z_1^\dag).
\label{2-21}
\end{equation}
\end{proposition}

\subsection{Canonical form away from $\Lambda$}
\label{sect-2-2}

\subsubsection{}
\label{sect-2-2-1} Let us consider first the canonical form in the subdomain, disjoint from $\Lambda$ which means that
\begin{equation}
|\{Z_1, f_1\}|\ge \epsilon_0.
\label{2-22}
\end{equation}
Let us reduce (\ref{2-17}) to a canonical form (multiplied by an elliptic operator). Note that the multiplication of $Z_1^\#$ by symbol $\beta$ implies the multiplication of $Z_1^{\#\dag}$ by $\beta^\dag$ and therefore the multiplication of $f_1^\#$ by $|\beta|^2$ modulo $\cO_1$ and finally the multiplication of
$\{Z_1^\#, f_1^\#\}$ by $|\beta|^2\beta$ modulo ${\tilde\cO}_1$ as well where here and below ${\tilde{\cO}_m}$ means the sum of monoms of the type
$\alpha _{ljk}\mu ^{l-m} x_1^l Z_1^{\#\,j}Z_1^{\#\dag\,k}$ with $m\ge j+k+l$ and smooth coefficients $\alpha_{ljk}$. Also to maintain (\ref{2-12}) and (\ref{2-17}) one needs to multiply operators $A$ and $\cA_n$ (and thus $f_2^\#,V^\#$ and all the perturbation terms) by $|\beta|^2$.

Then picking up
\begin{equation*}
\beta = \{Z_1^\#, f_1^\#\}^\dag \cdot |\{Z_1^\#, f_1^\#\}|^{-4/3}
\end{equation*}
one can achieve
\begin{equation}
\{Z_1^\#, f_1^\#\}=1+\alpha
\label{2-23}
\end{equation}
with $\alpha \in {\tilde\cO}_m$ and $m=1$.
Consequently, multiplying $Z_1^\#$
by $(1+\beta)$ with $\alpha \in {\tilde\cO}_m$ one changes
$\{Z_1^\#,f_1^\#\}$ by $(2\beta +\beta^\dag)\{Z_1^\#,f_1\}$ modulo ${\tilde{\cO}_{m+1}}$ and picking
$\beta ={\frac 1 3}\alpha^\dag - {\frac 2 3}\alpha$
one can achieve (\ref{2-23}) with $\alpha \in {\tilde\cO}_M$ with arbitrarily large $M$ (and after asymptotic summation with $M=\infty$). Unfortunately elements of ${\tilde\cO}_M$ (in contrast to elements of $\cO_M$) are not necessarily negligible unless $|x_1|\le h^\delta$.

However, it is easy to see that multiplying $Z_1^\#$ by $(1+\beta)$ with
$\beta \in \cO_m \cap {\tilde\cO}_\infty$ one changes $\{Z_1^\#,f_1^\#\}$ by
$\{Z_1^\#,\phi\}$ with
\begin{equation*}
\phi= f_1(\beta +\beta^\dag) +i\mu^{-1}\{Z_1^\#,\beta^\dag\}-
i\mu^{-1} \{Z_1^{\#\dag},\beta\}
\end{equation*}
modulo $\cO_m\cap{\tilde\cO}_\infty$. Obviously there exists
$\beta \in \cO_1\cap {\tilde\cO}_\infty$ eliminating error $\alpha$ in (\ref{2-23}) modulo $ \cO_1\cap {\tilde\cO}_\infty$. One can see easily that the only restriction to $\phi \in \cO_{m-1}\cap{\tilde\cO}_\infty$ is that it must be real valued.

Therefore step by step one can make
$\Re \alpha \in \cO_M\cap {\tilde\cO}_\infty$: $\{\Re Z_1, f_1\}\equiv 1$
$\mod \cO_M\cap {\tilde\cO}_\infty$ with arbitrarily large $M$. But then
setting
$f_1^\#= i\mu^{-1}\{Z_1^\#,Z_1^{\#\dag}\}= 2\mu^{-1}\{\Re Z_1^\#,\Im Z_1^\#\}$
one gets
\begin{equation*}
\bigl\{\Re Z_1^\#,\{\Im Z_1^\#,f_1\}\bigr\}= \bigl\{\Im Z_1^\#,\{\Re Z_1^\#,f_1\}\bigr\}
\in \cO_{M-1}\cap {\tilde\cO}_\infty
\end{equation*}
and since
$\{\Im Z_1^\#,f_1\}\in \cO_{M-1}$ as $f_1=0$ and $\{\Re Z_1^\#,f_1^\#\}\sim 1$ \ I derive that $\{\Im Z_1^\#,f_1\}\in \cO_{M-1}\cap {\tilde\cO}_\infty$. Then
\begin{equation}
-i\mu^{-1}\{ Z_1^\#, Z_1^{\#\dag} \} =2 f_1^\#,\qquad \{Z_1^\#,f_1^\#\}\equiv 1 \mod O(\mu^{-M})
\label{2-24}
\end{equation}
and then after an appropriate $\mu^{-1}h$-FIO transformation
\begin{equation}
Z_1^\#\equiv hD_1+i\bigl(hD_3 -{\frac 1 2}\mu x_1^2\bigr)
\label{2-25}
\end{equation}
(again I am using a linear independence of $\nabla \Re Z_1^\#$,
$\nabla \Im Z_1^\#$ and $\nabla f_1$).
However, one can replace operators $\beta(x',\mu^{-1}hD')$ by their decompositions into powers of $\mu^{-1}hD_1={\frac 1 2}(Z_1^\#+Z_1^{\#\dag})$ and $\mu^{-1}hD_3-{\frac 1 2}x_1^2 = {\frac 1 {2i}}\mu^{-1}(Z_1^\#-Z_1^{\#\dag})$
thus arriving to decomposition (\ref{2-26}) below with all the operators of the $\beta (x_1,x'', \mu^{-1}hD_4 ; \mu^{-1}h)$ and the following statement is proven:

\begin{proposition}\label{prop-2-7} Let conditions $(\ref{2-1})$ and $(\ref{2-22})$ be fulfilled. Then by means of $\mu^{-1}h$-transform, decomposition $(\ref{2-16})$ and one more $\mu^{-1}h$-transform one can reduce original operator $A$ to the family of 3-dimensional operators
\begin{multline}
\cA_n= {\frac 1 2}\sigma \Bigl({\color{blue}
(Z_1^\# \bullet Z_1^{\#\dag}) -W_n^\# }+\\
\sum_{2k+p+q+2m+2s\ge 3}^{} B^\w_{kpqms}\times \bigl((2n+1)\mu h\bigr)^k
Z_1^{\#\,p}Z_1^{\#\dag\,q}\,\mu ^{2-2k-2m-p-q-s} h ^s\Bigr)
\label{2-26}
\end{multline}
where symbol $ \sigma=|\{Z_1 ,f_1\}|^{2/3}$ is bounded and disjoint from 0, $Z_1^\#$ given by $(\ref{2-25})$ and

\begin{claim}\label{2-27} $\W_n^\#$ is an operator with the symbol
$\sigma^{-1}\bigl(V- f_2 (2n+1)\mu h\bigr) \circ \Psi$, $\Psi$ is the smooth diffeomorphism $\bR^4\to \bR^4$, with $D\Psi$ transforming $\Span \bigl( (0,1,0;0), (0,0,0;1)\bigr)$ into $\bK_2$,\newline
$\Span \bigl( ( 1, 0 , 0, 0), (0,0,1,0)\bigr)$ into $\bK_1$
and, in particular, $(0,0,1;0)$ into some element of $\bK_1 $ of $(0,*,*,*)$ form. Also $|\det D\Psi ^{-1}|= f_2$.\end{claim}
\end{proposition}
Again the first line in $(\ref{2-26})$ is the main part.

\subsubsection{}
\label{sect-2-2-2} Now let us consider domain where condition (\ref{2-22}) is violated. Note first that $ \bigl\{(x,\xi): f_1=Z_1=\{Z_1,f_1\}=0\bigr\}$ was an involutive manifold of codimension 5 in $\bR^4\times \bR^4$ with

\begin{claim}\label{2-28}
$\nabla f_1,\ \nabla \Re Z_1,\ \nabla \Im Z_1,\ \nabla \{\Re Z_1,f_1\}$ and $\nabla \{\Im Z_1,f_1\}$ linearly independent on $\Lambda'$\,$^{\ref{foot-13}}$.
\end{claim}
Also note that
\begin{equation*}
\{Z_1,x_1\}= x_3+ix_4, \quad \{Z_1,x_3\}=x_3+ix_4,\quad \{Z_1,x_4\}=x_4-ix_3
\qquad\text{as\ } x_1=0
\end{equation*}
and
\begin{equation}
\{\alpha^\#,\beta^\#\}=\{\alpha,\beta\}+
i\mu^{-1}f_2^{-1}\Bigl(\{Z_2^\dag,\alpha\}\{Z_2,\beta\}-
\{Z_2,\alpha\}\{Z_2^\dag,\beta\}\Bigr)+O(\mu^{-2}).
\label{2-29}
\end{equation}
However $\{Z_2,x_j\}$ are actually arbitrary and therefore while we know that
\begin{equation}
\Lambda^\#=\bigl\{f_1^\#=Z_1^\#=\{Z_1^\#,f_1^\#\}=0\bigr\}
\label{2-30}
\end{equation}
is a manifold of codimension 5 with
\begin{claim}\label{2-31}
$\nabla f^\#_1,\ \nabla \Re Z^\#_1,\ \nabla \Im Z^\#_1,\ \nabla \{\Re Z^\#_1,f_1\}$ and $\nabla \{\Im Z^\#_1,f^\#_1\}$ linearly independent on $\Lambda^\#$\,$^{\ref{foot-13},\ref{foot-14}}$;
\end{claim}
but its symplectic structure is not fixed and there is no canonical form uniformly near $\Lambda^\#$.

Instead let us consider point ${\bar z}$ ($=(0,0)\in \bR^6$ for simplicity of notations) with
\begin{equation}
|f_1|\le \gamma= \epsilon r^2,\quad |\{Z_1,f_1\}|=r \qquad \text{with\ } r\ge \mu^{\delta-1/3},
\label{2-32}
\end{equation}
and reduce an operator to a canonical form in its $\epsilon(\gamma, r)$-vicinity with respect to $(x_1;x_2,x_3,x_4)$, leaving domain $\{|f_1|\ge \epsilon |\{Z_1,f_1\}|^2\}$ for a later.

Then after reduction of the previous subsection condition (\ref{2-32}) remains valid for $Z_1^\#$ and $f_1^\#$; therefore now multiplication by a symbol
$\beta = \{Z_1^\#,f_1^\#\}^\dag\cdot |\{Z_1^\#,f_1^\#\}|^{-4/3}r^{1/3}$ provides a modified equality (\ref{2-23}), namely,
\begin{equation}
\{Z_1^\#, f_1^\#\}=r+\alpha.
\label{2-33}
\end{equation}
However, $\beta$ is an uniformly smooth symbol; more precisely $\beta={\hat\beta} (r^{-1}x', r^{-1} \xi')$ with uniformly smooth symbol ${\hat\beta}$, which means that after original rescaling $x'\mapsto rx'$\,\footnote{\label{foot-15} Since symplectic structure of $\Lambda^\#$ is not defined one cannot confine this non-smoothness to a couple of coordinates.}
$\beta={\hat\beta} (x, \mu^{-1}hr^{-2}D')$ and $x_1^{-1}\alpha$ is of the same type.

Continuing as before\footnote{\label{foot-16} With effective semiclassical parameter $\mu^{-1}hr^{-2}$ rather than $\mu^{-1}h$.} one can achieve a modified equality (\ref{2-24}), namely,
\begin{equation}
-i\mu^{-1}r^{-1}\{ Z_1^\#, Z_1^{\#\dag} \} =2f_1^\#,\qquad \{Z_1^\#,f_1^\#\}\equiv r \mod O(\mu^{-s})
\label{2-34}
\end{equation}
where now $\mu^{-1}r^{-1}Z_1$ is $\mu^{-1}hr^{-2}$-PDO
and then modified equality (\ref{2-25})
\begin{equation}
Z_1^\#\equiv \bigl( r h_* D_1+ h_* D_4\bigr) +i\bigl(h_*D_3 -{\frac 1 2}\mu r x_1^2\bigr)
\label{2-35}
\end{equation}
where extra terms $h_*D_4$ and $h_*D_3$ appear because (original) $\nabla \Re Z_1$ and $\nabla \Im Z_1$ are linearly independent, $h_*=hr^{-1}$.

Then after another rescaling $x_1\mapsto r x_1$, $D_1\mapsto r^{-1}D_1$ (\ref{2-35}) becomes
\begin{equation*}
Z_1^\#\equiv
\bigl( h_* D_1+ h_* D_4\bigr)+i\bigl(h_*D_3 -{\frac 1 2}\mu r^3 x_1^2\bigr)
\end{equation*}
and after transformation $(x_1,x_3,x_4; D_1,D_3,D_4)\mapsto (x_1,x_3, x_4+x_1;D_1-D_4,D_3,D_4)$ becomes
\begin{equation}
Z_1^\#= h_*D_1 +i\bigl(h_*D_3 -{\frac 1 2}\mu_* x_1^2\bigr),\qquad h_*=hr^{-1},\quad \mu_*=\mu r^3
\label{2-36}
\end{equation}
which is exactly equality (\ref{2-25}) with $h$ and $\mu$ replaced by $h_*$ and $\mu_* $ respectively.

Meanwhile, in the original operators of type
$\beta (r^{-1}x',\mu^{-1}hr^{-1} D';\mu^{-1}hr^{-2})$ (with the smooth symbol $\beta$) were transformed subsequently into
$\beta (x', \mu^{-1}h r^{-1}D' ;\mu^{-1}hr^{-2})$, then into\newline
$\beta (rx_1,x'', \mu^{-1}hr^{-3}D_1,\mu^{-1}hr^{-2}D''; \mu^{-1}hr^{-2})$
and finally into
\begin{equation*}
\beta (rx_1,x'', \mu^{-1}hr^{-3}D_1,\mu^{-1}hr^{-2}D_3, \mu^{-1}hr^{-2}D_4 ; \mu^{-1}hr^{-3}).
\end{equation*}
So, I arrive to (\ref{2-17})-like decomposition with $Z_1^\dag$ given by (\ref{2-36}) and $\sigma$, $W_n$ and $B^\w_{kpqms}$ of the above type. However, one can replace such operators by their decompositions into powers of $\mu^{-1}hr^{-3}D_1={\frac 1 2}\mu^{-1}r^{-2}(Z_1^\#+Z_1^{\#\dag})$ and $\mu^{-1}hr^{-2}D_3-{\frac 1 2}r^2x_1^2 = {\frac 1 {2i}}\mu^{-1}r^{-1}(Z_1^\#-Z_1^{\#\dag})$
thus arriving to decomposition (\ref{2-38}) below with all the operators of the simplified type
\begin{equation}
\beta (rx_1,x'', \mu^{-1}hr^{-2}D_4 ; \mu^{-1}hr^{-3});
\label{2-37}
\end{equation}
factors $\mu^{-1}hr^{-4}$ appear because one needs to commute $Z_1^\#,Z_1^\#$ between themselves and with operators of (\ref{2-37})-type.

So, I finally arrive to

\begin{proposition}\label{prop-2-8} Let conditions $(\ref{2-1})$ and $(\ref{2-32})$ be fulfilled. Then by means of $\mu^{-1}h$-FIO transform, decomposition $(\ref{2-14})$ and the series of rescalings and $\mu^{-1}h r^{-2}$-FIO transform one can reduce described above original operator $A$ to the family of 3-dimensional operators
\begin{multline}
\cA_n\Def{\frac 1 2}\sigma^\# \Bigl({\color{blue}
(Z_1^\# \bullet Z_1^{\#\dag}) -W_n^\# }+\\
\sum_{2k+p+q+2m+2s\ge 3}^{} B^\w_{kpqms}\times \bigl((2n+1)\mu h\bigr)^k
Z_1^{\#\,p}Z_1^{\#\dag\,q}\,\mu ^{2-2k-2m-p-q-s} h ^s r^{-2p-2q-4m-4s}\Bigr)
\label{2-38}
\end{multline}
where $\sigma^\#$ is transformed operator with symbol $\sigma = |\{Z_1 ,f_1\}|^{2/3}r^{-2/3}$ which is bounded and disjoint from $0$ and $Z_1^\#$ given by $(\ref{2-35})$.
\end{proposition}
Again the first line in $(\ref{2-38})$ is the main part.

\begin{remark}\label{rem-2-9}
(i) Here operators $\sigma^\#$, $W_n^\#$ and $B^\w_{kpqms}$ are of the (\ref{2-37})-type
with uniformly smooth symbol $\beta$. Therefore these symbols are quantizable provided under condition (\ref{2-32}) for sure;

\smallskip
\noindent
(ii) I leave for a section 3 the the calculation of the symbol of $W$ and the discussion of the corresponding diffeomorphism.
\end{remark}
\begin{figure}[h!]
\centering
\caption{\label{fig-1} Reduction was done in the colored zone.}
\end{figure}

\subsubsection{}\label{sect-2-2-3}
Let us look what happens with the cut-off symbol $\psi$ as a result of these transformations. After transformations of subsubsections \ref{sect-2-1-1}, \ref{sect-2-1-2} after decomposition (\ref{2-16}) it becomes a matrix operator
\begin{equation*}
\sum _{k,l} \phi_{kl} (x',\mu^{-1}hD') \mu^{-k-l} U_2^k U_2^{\dag l}
\end{equation*}
with basis-shifting (in $L^2(\bR)$) operators $(2(n+1)\mu h)^{-1/2}U_2^\dag$ and
$(2n\mu h)^{-1/2}U_2$.

Then after transformations of subsubsection \ref{sect-2-2-2} these matrix operator becomes
\begin{equation}
\sum_{k,l,p,q,m,s}^{} \phi^\w_{klpqms}
U_2^k U_2^{\dag l}
Z_1^{\#\,p}Z_1^{\#\dag\,q}\,\mu ^{-k-l-m-p-q-s} h ^s r^{-2p-2q-2m-4s}
\label{2-39}
\end{equation}
with the scalar main part of the form (\ref{2-37}).

\subsection{Canonical form in the strictly outer zone}
\label{sect-2-3}

\subsubsection{}
\label{sect-2-3-1}
Now one can properly redefine the first part of the \emph{outer zone\/}
\begin{equation}
\cZ_{\out,I} = \Bigl\{ C\mu^{-1/2} \dist (x,\Lambda)^{1/2} \le \dist(x,\Sigma)\le \epsilon \dist (x,\Lambda)^2\Bigr\}
\label{2-40}
\end{equation}
which in the notations $\dist(x,\Sigma)=|x_1|$ and
$\dist (x,\Lambda)\asymp r$ is $\bigl\{\mu^{-1/2}r^{1/2}\le |x_1|\le \epsilon r^2\bigr\}$. Note that after rescaling of the previous section $x_1 \mapsto x_1r^{-2}$, $\mu \mapsto \mu_*=\mu r^3$ it becomes $\bigl\{|x_1|\ge C\mu_*^{-1/2}\bigr\}$.

In this subsubsection I restrict myself to the part of \emph{strictly outer zone\/}
\begin{equation}
\cZ^*_{\out,I} = \Bigl\{ C\mu^{3\delta-1/2} \dist (x,\Lambda)^{1/2} \le \dist(x,\Sigma)\le \epsilon \dist (x,\Lambda)^2\Bigr\}
\label{2-41}
\end{equation}
which in the same notations is $\bigl\{\mu^{-1/2}r^{1/2}\le |x_1|\le \epsilon r^2\bigr\}$ and after becomes $\bigl\{|x_1|\ge C\mu_*^{3\delta-1/2}\bigr\}$ (may be with a different $\delta>0$) and provides a nice quantization of all symbols below.

Then since $\mu_*\ge h^{-\delta}$ one can apply the canonical form of the Schr\"odinger operator with the strong magnetic field. Let us consider an operator obtained after reduction of the previous subsection; there
$|x_1|\asymp \gamma_*=\gamma r^{-2}$ if originally $|x_1|\asymp \gamma$.
Notation $\gamma=r^2$ was temporary.

Applying rescaling $(x_1,x_3)\mapsto (x_1\gamma_*,x_3\gamma_*)$ one arrives to transformed (\ref{2-36}) of the same type but with $h_*$ and $\mu_*$ replaced by $h_*\gamma_*^{-1}=hr\gamma^{-1}$ and $\mu_*\gamma_*^2=\mu \gamma^2r^{-1}$ respectively:
\begin{equation}
Z_1^\# =
\mu \gamma^2 r^{-1}\Bigl( \hbar D_1 + i \bigl( \hbar D_3 +
{\frac 1 4} x_1^2\bigr)\Bigr)
\label{2-42}
\end{equation}
with $|x_1|\asymp 1$ and
\begin{equation}
\hbar =\mu^{-1}hr^2\gamma^{-3}\le c\mu ^{1/2}h^{1+3\delta} r^{1/2}\ll c h^{1/2}.
\label{2-43}
\end{equation}
Also $\mu^{-1}hr^{-2}D_3 \mapsto \mu^{-1}h\gamma^{-1}D_3$.

Let us consider the part where $x_1>0$; another part is separated and treated exactly in the same manner.
Introducing
\begin{equation}
U_1= x_1^{-1/2}\Bigl( \hbar D_1 + i \bigl( \hbar D_3 +
{\frac 1 4} x_1^2\bigr)\Bigr)
\label{2-44}
\end{equation}
one can see easily that
\begin{align}
&|U_1|\le c\nu^{-1},\qquad \nu \Def \mu\gamma^2r^{-1}\ge h^{-\delta},\label{2-45}\\
&\{U_1,U_1^\dag\}\equiv 1\quad \mod U_1,U_1^\dag\label{2-46}
\end{align}
where (\ref{2-45}) follows from $|Z^\#_1|\le c$ and (\ref{2-46}) is understood in the sense of $\hbar$-symbols.

Therefore there exists $\hbar$-FIO transformation $\cF'$ of
$(x_1,x_3; \hbar D_1,\hbar D_3)$ thus not affecting $(x_4, \mu^{-1}hr^{-2} D_4)$ so that
\begin{equation}
U_1\mapsto U_1^\# \equiv \hbar D_1 + i x_1 \qquad \mod \cO_2 (U_1,U_1^\dag);
\label{2-47}
\end{equation}
then operator (\ref{2-38}) is transformed into
\begin{multline}
{\frac 1 2}\sigma \sigma_1 \Bigl({\color{blue}
\nu^2 U_1^\#\bullet U_1^{\#\dag} -W_n^{\#\#} }+\\
\sum_{2k+p+q+2m+2s+2l\ge 3}^{} B^\w_{kpqmsl}\times \bigl((2n+1)\mu h\bigr)^k
\bigl(\nu U_1^\# \bigr)^p \bigl(\nu U_1^{\#\dag}\bigr)^q \,\mu ^{2-2k-2m-p-q-s} h ^s r^{-2p-2q-4m-4s} \hbar ^l\Bigr)
\label{2-48}
\end{multline}
with $W_n^{\#\#}$ the result of transformation of $x_1^{-1}W_n^\#$.

Here all the operators $\sigma$, $\sigma_1$, $W_n^{\#\#}$ and $B^\w_{kpqmsl}$ of the type
\begin{equation*}
\beta (x_1, x_3, x_4; \hbar D_1, \hbar D_3, \mu^{-1}hr^{-2}D_4).
\end{equation*}
Note that $\hbar \nu =h r\gamma^{-1}\ll \mu^{1/2} h r^{1/2} \le h^{1/2}$ in the zone in question. One can decompose all such operators into power series with respect to $x_1$, $\hbar D_1$ and rewrite (\ref{2-48}) with all operators of the type
\begin{equation}
\beta ( x_3, x_4; \hbar D_3, \mu^{-1}hr^{-2}D_4),\qquad \hbar= \mu^{-1}hr^2\gamma^{-3}
\label{2-49}
\end{equation}
and with $\bigl(\nu U_1^\# \bigr)^p$, $\bigl(\nu U_1^{\#\dag}\bigr)^q$ replaced by $\bigl(\nu U_1^\# \bigr)^{p+p'}\nu^{-p'}$, $\bigl(\nu U_1^{\#\dag}\bigr)^{q+q'}\nu^{-q'}$ respectively since $\nu \ge h^{-\delta}$.

Repeating arguments of the previous subsections one can transform this operator into
\begin{align}
&{\frac 1 2}\sigma \sigma_1 \Bigl({\color{blue}
\nu^2 U_1^\#\bullet U_1^{\#\dag} -W_n^{\#\#} }+\label{2-50}\\
&\smashoperator{\sum_{2q+2l+2j\ge 3}^{} }\qquad B^\w_{qjl}\times 
\bigl(\nu^2 U_1^\#\bullet U_1^{\#\dag}\bigr)^q \, \nu^{2-2q-2j}\hbar ^l+\notag\\
&\smashoperator{\sum_{2k+2q+2m+2s+2l\ge 3}^{}}\qquad B^\w_{kqjmsl}\times \bigl((2n+1)\mu h\bigr)^k
\bigl(\nu^2 U_1^\#\bullet U_1^{\#\dag}\bigr)^q \,\mu ^{2-2k-2m-2q-s} h ^s r^{-4q-4m-4s} \nu^{-2j} \hbar^l\Bigr)\notag
\end{align}
where the second line comes from decomposing of $W_n^{\#\#}$ and because it is rescaled before
\begin{equation}
|B_{qjl}|\le C\gamma r^{-1},
\label{2-51}
\end{equation}
Now one can apply decomposition with respect to
$\hbar^{-1/4} \upsilon_p (\hbar ^{-1/2}x_1)$:
\begin{equation}
u_{nn'}(x',y',t)=\sum_{p,p'\in \bZ^+} \hbar^{-1} u_{npn'p'}(x'',y'',t)\,
\upsilon_p\bigl(\hbar^{-1/2} x_1)
\upsilon_{p'}\bigl(\hbar^{-1/2} y_1)
\label{2-52}
\end{equation}
with $x''=(x_3,x_4)$. Then the operator (\ref{2-50}) becomes
\begin{align}
&{\frac 1 2}\sigma \sigma_1 \Bigl({\color{blue}
(2p+1)\mu h \gamma -W_n^{\#\#} }+\label{2-53}\\
&\smashoperator{\sum_{2q+2l+2j\ge 3}^{} }\qquad B^\w_{qjl}\times 
\bigl((2p+1)\mu h\gamma\bigr)^q \, \nu^{2-2q-2j}\hbar ^l+\notag\\
&\smashoperator{\sum_{2k+2q+2m+2s+2l\ge 3}^{}}\qquad B^\w_{kqjmsl}\times \bigl((2n+1)\mu h\bigr)^k
\bigl( (2p+1) \mu h\gamma \bigr)^q \,\mu ^{2-2k-2m-2q-s} h ^s r^{-4q-4m-4s} \nu^{-2j} \hbar^l\Bigr).\notag
\end{align}
since $\nu^2 \hbar = \mu h \gamma$ and $2k+2q+2m+2s+2j+2l\ge 3$ is equivalent to $k+q+m+s+j+l\ge 2$.

In the decomposition (\ref{2-53}) the role of $m$, $s$, $j$ and $l$ is just to bound a magnitude and indicate dependence on $\mu$ and $h$ since dependence on $r,\gamma$ (but not magnitude) is not important. However $\nu =\mu \cdot \gamma^2r^{-1}\le \mu r^2$, $\hbar = \mu^{-1}h\cdot r^2\gamma^{-3}\ge \mu^{-1}hr^{-4}$; therefore one can join terms with $s>0$ or $m>0$ to the terms with $s=0$ and $m=0$ without changing $s+l$, $m+j$ thus getting the simplified correction term
\begin{equation*}
\sum_{k+q+j+l\ge 2}^{} B^\w_{kqjmsl}\times \bigl((2n+1)\mu h\bigr)^k
\bigl( (2p+1) \mu h\gamma \bigr)^q \,\mu ^{2-2k-2q} h ^s r^{-4q} \nu^{-2j} \hbar^l
\end{equation*}

Therefore I arrive to the following

\begin{proposition}\label{prop-2-10} Let condition $(\ref{2-1})$ be fulfilled. Let us consider zone $\cZ^*_{\out,I}$, described by (\ref{2-40}).

Then by means of $\mu^{-1}hr^{-2}$-FIO transforms, decomposition $(\ref{2-16})$, the ``special'' $\hbar$-FIO transform and decomposition $(\ref{2-52})$ one can reduce the original operator $A$
to the family of 2-dimensional operators
\begin{multline}
\cA_{np}\Def {\frac 1 2}\sigma \sigma_1 \Bigl({\color{blue}
H_{pn} }+\\
\sum_{k+q+j+l\ge 2}^{} B^\w_{kqjmsl}\times \bigl((2n+1)\mu h\bigr)^k
\bigl( (2p+1) \mu h\gamma \bigr)^q \,\mu ^{2-2k-2q} h ^s r^{-4q} \nu^{-2j} \hbar^l\Bigr)
\label{2-54}
\end{multline}
with $n\le C_0/(\mu h)$, $p\le C_0/(\mu h\gamma)$ where all operators are of $(\ref{2-49})$-type and $\sigma,~\sigma_1$ are bounded and disjoint from $0$.
\end{proposition}

\begin{remark}\label{rem-2-11} (i) While in the canonical form (\ref{2-38}) only one cyclotron movement was separated, in the canonical form (\ref{2-54}) both of them are separated leaving only pure drift movement which will be studied in the next section;

\smallskip
\noindent
(ii) I leave for a section 3 the discussion of the actual elements where reduction was done, corresponding diffeomorphism and the calculation of the symbol of $W$;

\smallskip
\noindent
(iii) Cut-off symbol $\psi$ which after previous transformations became (\ref{2-39}) now becomes
\begin{equation}
\sum_{k,i, p,q, s,l}^{} \phi_{kipqsl}\times U_2^i U_2^{\dag\,k}
U_1^{\#\,p} U_1^{\#\dag\,q} \,\mu ^{-i-k-2m-p-q-s} h ^s r^{-2p-2q-4m-4s} \nu^{-2j-i-k} \hbar^l
\label{2-55}
\end{equation}
with (\ref{2-49})-type operators $\phi_{kipqsl}$.
\end{remark}

\subsubsection{}
\label{sect-2-3-2}
Now my goal is to achieve a similar reduction in the second part of the \emph{outer zone}
\begin{equation}
\cZ^*_{\out,II} = \Bigl\{ C\max\bigl( \dist (x,\Lambda)^2, \mu^{4\delta-2/3}\bigr) \le \dist(x,\Sigma)\le \epsilon \Bigr\}.
\label{2-56}
\end{equation}
Then the previous subsection \ref{sect-2-2} becomes irrelevant; however I can start from (\ref{2-17}). Again with no loss of the generality one can assume that $f_1^\#=x_1$. Let us consider a point $\bz=(\bx',\bxi') \in \bR^6$ and an element (\ref{2-55}) with $\gamma=|{\bar x}'_1|$ and $r=\gamma^{1/2}$ (exactly as in the previous subsection it was defined the other way. Without any loss of the generality one can assume that $\bz = (\gamma, 0,\dots,0)$.

Consider Taylor decomposition (with a remainder) of $\mu^{-1}Z_1^\#$ with respect to $x'$, $\mu^{-1}hD'$ and $r$ (in $r$-vicinity of $0$ all of them are $O(r)$); then due to $\{Z_1^\#,Z_1^{\#\dag}\}=O(r)$, $\{Z_1,x_1\}=O(r)$ with no loss of the generality one can assume that
$\mu ^{-1}Z_1^\# = \mu^{-1} hD_3+i\mu^{-1}hD_4 +O(r^2)$.
Note that $\mu^{-1}hD_1$ could be either as square (or higher degree) or with a linear or quadratic factor with respect to other variables. Also note that the quadratic terms in $O(r^2)$ containing exclusively $x_3,x_4$ must commute with $D_3+iD_4$ (otherwise $\{Z_1^\#,Z_1^{\#\dag}\}\sim 2x_1$ would not be possible) and therefore must be of the form $c(x_3+ix_4)^2$. However one could remove such terms by $\mu^{-1}h$-PDO affecting only $x'',\mu^{-1}hD''$ and leaving $x_1$, $\mu^{-1}hD_1$ intact. So, terms in the decomposition must contain either one of the factors $x_1$, $\mu^{-1}hD_3$, $\mu^{-1}hD_4$, or two factors $\mu^{-1}hD_1$, or three factors $x_3,x_4$; also $x_1$ must be accompanied by $O(r)$-type factor.

Then scaling
$x''\mapsto rx'$, $x_1\mapsto r^2 x_1$, $D''\mapsto r^{-1}D''$, $D_1\mapsto r^{-2}D_1$ and dividing by $r^3$ one gets completely legitimate $\mu^{-1}h r^{-4}$-PDO where each factor $\mu^{-1}hD''$ can accommodate division by $r^3$ and $\mu^{-1}hD_1$ can accommodate division by $r^2$ but there either at least two such factors, or there is an extra $O(r)$ type factor. This operator, with
the semiclassical parameter $\hbar=\mu^{-1}h\gamma^{-2}$ satisfies again $\{Z_1^\#,Z_1^{\#\dag}\}\sim 2x_1$ but must be considered near point $(1,0,\dots,0)$ and then by the standard arguments one can reduce operator to the canonical form (\ref{2-54}) with all the operators of the type
\begin{equation}
\beta (x_3, x_4; \hbar D_3, \hbar D_4),\qquad \hbar= \mu^{-1}h\gamma^{-2}.
\label{2-57}
\end{equation}
Note that as $r=\gamma^{1/2}$ this is consistent with the definition of $\hbar$ in the zone $\cZ_{\out,I}$ but form (\ref{2-57}) is more general then (\ref{2-49}).

So, I arrive to

\begin{proposition}\label{prop-2-12} Let condition $(\ref{2-1})$ be fulfilled.
Then one can reduce the original operator $A$ to the family $(\ref{2-54})$
of 2-dimensional operators with $n\le C_0/(\mu h)$, $p\le C_0/(\mu h\gamma)$ where all operators are of $(\ref{2-49})$-type and $\sigma,~\sigma_1$ are bounded and disjoint from $0$.
\end{proposition}

\begin{remark}\label{rem-2-13}
(i) All statements (i)-(iii) of remark \ref{rem-2-11} remain valid.
\end{remark}

\subsubsection{}
\label{sect-2-3-3}
I would like to finish this section by the following

\begin{remark}\label{rem-2-14} In $2D$-case canonical form would be \cite{IRO6}
\begin{multline}
\cA_p=\sigma \sigma_1 \Bigl({\color{blue}
(2p+1) \mu h \gamma -W ^{\#\#} }+\\
\sum_{2l+2j+2t\ge 3}^{} B^\w_{ljt}\times
\bigl((2p+1)\mu h \gamma \bigr)^l\times
\mu ^{2-2l-2j-t} h^t \gamma^{-4j-3t }\Bigr)
\label{2-58}
\end{multline}
with
\begin{equation}
W^{\#\#}_n=\beta ( \gamma^{-1} x_3, \mu^{-1}h \gamma^{-2} D_3).
\label{2-59}
\end{equation}
There instead of 2-parameter family of 2-dimensional PDOs I had only 1-parametric family of 1-dimensional PDOs.
\end{remark}

\section{Estimates}\label{sect-3}

In this section I derive remainder estimates with the main part in the standard but rather implicit form which is the sum of expressions (\ref{0-12})
\begin{equation}
h^{-1}\sum_{\iota}\int_{-\infty}^0 \Bigl( F_{t\to h^{-1}\tau} {\bar\chi}_{T_\iota}(t) \Gamma uQ_{\iota y}^t\Bigr)\,d\tau
\label{3-1}
\end{equation}
with $Q_\iota$ partition of unity and an appropriate $T_\iota$ where as in my previous papers ${\bar\chi}_T(t)={\bar\chi}(t/T)$ and ${\bar\chi}$ is supported in $[-1,1]$ and equal 1 at $[-{\frac 1 2},{\frac 1 2}]$, while the remainder estimate could depend on non-degeneracy condition.

I am going to consider different zones and apply different approaches and use different canonical forms in each of them; I start from the strictly outer zone.
\begin{figure}[h!]
\centering
%
\includegraphics{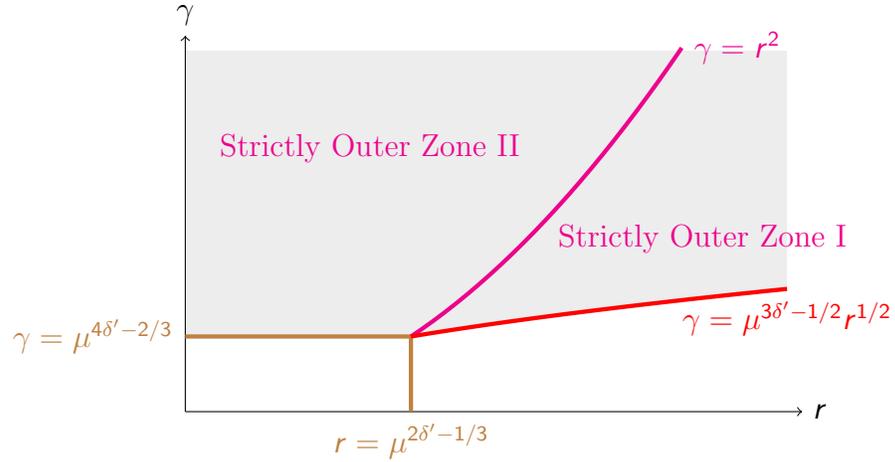}
\caption{\label{fig-2} My first target is strictly outer zone.}
\end{figure}

\subsection{Estimates in the strictly outer zone. I}
\label{sect-3-1}

I will start from the most massive \emph{Strictly Outer Zone I} $\cZ^*_{\out,I}$ described in (\ref{2-40}) with an arbitrarily small exponent $\delta_0>0$ and sufficiently small constant
$\epsilon >0$.

\subsubsection{}
\label{sect-3-1-1}
After rescalings and transformations described in the previous section the original operator \emph{somewhere} there was reduced to the family of two-dimensional $(\hbar_3, \hbar_4)$-PDOs $\cA_{pn}(x'',\hbar_3 D_3,\hbar_4 D_4)$ given by (\ref{2-54}) with
\begin{equation}
\hbar_3 =\mu ^{-1}h r^2\gamma^{-3}\ll \mu^{1/2}h,\qquad
\hbar_4 =\mu^{-1}h r^{-2} \ll \mu^{-1/3}h;
\label{3-2}
\end{equation}
one needs to consider only magnetic numbers
\begin{equation}
n\le C_0(\mu h)^{-1},\qquad p \le C_0(\mu h \gamma )^{-1}\qquad \qquad \text{as\ \ }\dist(x,\Sigma)\asymp \gamma;
\label{3-3}
\end{equation}
these two restrictions mean that energy levels do not exceed $\epsilon$.

These final operators are defined in $B(0,1)$. Therefore
\begin{claim}\label{3-4}
A contribution of each {\sf (pair,element)}\footnote{\label{foot-17} Where here and below a pair is $(p,n)$ and an element means the final element $B(0,1)$.} to to the asymptotics is $O\bigl(\hbar_3^{-1}\hbar_4^{-1}\bigr)=\linebreak O\bigl(\mu^2h^{-2}\gamma^3\bigr)$\,\footnote{\label{foot-18} And for many pairs it is actually of this
magnitude.}.
\end{claim}
However, one can see easily that this element before the final reduction was the ``slab'' $B(0,1)\cap \{|x_1- {\bar x_1}| + |x_3-{\bar x}_3|\le \epsilon \gamma_*\}$ with $\gamma_*=\gamma r^{-2}$, $|{\bar x_1}|\asymp\gamma_*$ where here $B(0,1)$ is an image of the original element $B(y', r)\cap \cZ_{\out, I}$. Therefore the total contribution of the given pair $(p,n)$ over subzone
\begin{equation*}
\cZ_{(r,\gamma)}=
\bigr\{ \gamma \le \dist(x,\Sigma) \le 2\gamma, r\le \dist (x,\Lambda )\le 2r\bigr\}
\end{equation*}
with $r\ge \epsilon \gamma^{1/2}$ to the asymptotics is
$O\bigl( \mu^2h^{-2}\gamma^3 \times r^2\gamma ^{-1}\bigr) =
O(\mu ^2h^{-2}\gamma ^2 r^2\bigr)$ and the total contribution of subzone $\cZ_{(r,\gamma)}$ to the asymptotics is
$O\bigl(\mu^2h^{-2}\gamma^2r^2 \times (\mu h\gamma )^{-1} \times (\mu h )^{-1} \bigr) = O(h^{-4}\gamma r^2\bigr)$ where the second and the third factors are numbers of permitted indices $p$ and $n$ respectively.

Then the total contribution of zone $\cZ^*_{\out,I}$ to asymptotics is $O(h^{-4})$\,\footnote{\label{foot-19} And is actually of this magnitude under condition (\ref{0-2}) as $\mu \le \epsilon h^{-1}$. These results are consistent with what one can get without reduction.}.

The same approach works in $\cZ^*_{\out,II}$: here one has the family of two-dimensional $\hbar_3$-PDO $\cA_{pn}(x'',\hbar_3 D_3,\hbar_3 D_4)$ with
$\hbar_3=\mu^{-1}h\gamma^{-2}$. Therefore a contribution of each {\sf(pair,element)} to to the asymptotics is $O(\hbar_3^{-2})=O(\mu^2h^{-2}\gamma^4)$. However there was nondiscriminatory rescaling $(x_3,x_4,\xi_3,\xi_4)\mapsto (r^{-1}x_3,r^{-1}x_4,r^{-1}\xi_3,r^{-1}\xi_4)$ while the number of such balls should be $O(r^{-2})$ since $\codim_\Sigma \Lambda=2$.

Therefore the total contribution of the given pair $(p,n)$ over subzone
\begin{equation*}
{\bar\cZ}_{(r,\gamma)}=
\bigr\{ \gamma \le \dist(x,\Sigma) \le 2\gamma, \dist (x,\Lambda )\le 2r\bigr\}
\end{equation*}
with $r=\gamma^{1/2}$ to the asymptotics is $O(\mu^2h^{-2}\gamma^3)$, the total contribution of subzone ${\bar\cZ}_{(r,\gamma)}$ to the asymptotics is
$O(h^{-4}\gamma^2)$ and the total contribution of zone $\cZ^*_{\out,I}$ to asymptotics is $O(h^{-4})$.

Now I need to analyze the contribution of such elements to the remainder estimate.

\subsubsection
\label{sect-3-1-2}
Let us consider some index pair $(p,n)$ and some final $B(0,1)$ element. Assume first that after all these transformations operator $\cA_{pn}$ is \emph{regular\/} on this {\sf(pair,element)} in the sense of (\ref{3-5}),(\ref{3-6}) and (\ref{3-8})
\begin{align}
&\eta_\alpha|\nabla_{x'',\xi''}^\alpha \cA_{pn}|\le c_\alpha\ell \qquad \forall \alpha:0<|\alpha|\le K,\label{3-5}\\
\intertext{with $\eta_\alpha=1$ as $\alpha =(*,0,*,0)$ and $\eta_\alpha=\gamma r^{-2}$ otherwise and $\ell$ satisfying}
&\ell\ge h^{-K_0}.\label{3-6}
\end{align}
Then the propagation speed with respect to $(x_3,\hbar_3D_3)$ is
$O(\ell \hbar_3/h )= O( \ell \mu ^{-1}r^2\gamma^{-3})$
and the propagation speed with respect to $( x_4,\hbar_4D_4)$
is $O(\ell \hbar_4/h\times r^2/\gamma )=O(\ell \hbar_3/h )$ as well.
Therefore the dynamics starting in $B(0,{\frac 1 2 })$ is retained in $B(0,1)$ as $|t|\le T_1$ with
\begin{equation}
T_1= \epsilon h\ell^{-1}\hbar_3^{-1} \asymp \mu \gamma^3 r^{-2}\ell^{-1};
\label{3-7}
\end{equation}
I remind that $\hbar_3\gg \hbar_4$.
In the same time under respective assumptions
\begin{equation}
|\nabla_{x_3,\xi_3} \cA_{pn}|+\zeta |\nabla _{x_4,\xi_4}\cA_{pn}| \asymp \ell,\quad \zeta=\gamma r^{-2}
\label{3-8}
\end{equation}
one can pick up $T_0$ from the uncertainty principle\,\footnote{\label{foot-20} I don't need the logarithmic uncertainty principle in this place and use the simpler version of the microlocal uncertainty principle.}
$\ell \hbar_j h^{-1} T_0\ge \hbar_j h^{-\delta'}$\,\footnote{\label{foot-21} Alternatively one could scale $x_4\mapsto (\hbar_4/\hbar_3)^{1/2} x_4$, $\xi_4\mapsto (\hbar_4/\hbar_3)^{1/2} \xi_1$, $\hbar_4\mapsto \hbar_3$.} or equivalently
\begin{equation}
T_0= C h^{1-\delta'} \ell^{-1}.
\label{3-9}
\end{equation}
with an arbitrarily small exponent $\delta'>0$.

Then the contribution of this {\sf(pair,element)} to the remainder estimate does not exceed
\begin{equation}
C\hbar_3^{-1}\hbar_4^{-1} \times T_0T_1^{-1}\asymp \hbar_4^{-1}h^{-\delta'} =
\mu h^{-1-\delta'} r^2
\label{3-10}
\end{equation}
where $C\hbar_3^{-1}\hbar_4^{-1} \times T_0$ estimates $\sup_{|\tau|\le \epsilon} |F_{t\to h^{-1}\tau} {\bar\chi}_T(t)\Gamma \bigl(u\psi\bigr)|$ and factor $T_1^{-1}$ is due to the standard Tauberian approach.

So, in comparison with the estimate of the contribution to the asymptotics (\ref{3-4}) the estimate of the contribution to the asymptotics picked an extra factor $\hbar_3 h^{-\delta'}$ with an arbitrarily small exponent $\delta'>0$; actually one can reduce $\delta'$ to $0$ but it is not needed here.

Let us pick up {\sf($p$,element)}, then break this element into $\zeta $-subelements with respect to $(x_4,\xi_4)$; then on each such subelement operator $\cA_{pn}$ is elliptic\footnote{\label{foot-22} See next subsubsection for justification; to prove (\ref{3-14}) and (\ref{3-15}) which are main results of this subsubsection one does not need these arguments.}
$|\cA_{pn}|\ge \epsilon \ell$
for all indices $n$ but $C\bigl(\ell/(\mu h) + 1\bigr)$. Therefore the total contribution to the remainder the pair {\sf($p$,element)} combination such that conditions (\ref{3-5}),(\ref{3-6}) and (\ref{3-8}) are fulfilled with $\ell \le {\bar\ell}$ does not exceed
\begin{equation}
C\mu h^{-1-\delta'} r^2\times \bigl({\bar\ell}(\mu h)^{-1}+1\bigr)=
C\bigl({\bar\ell}+\mu h\bigr) h^{-2-\delta'} r^2.
\label{3-11}
\end{equation}
Then summation over $\cZ_{(r,\gamma)}$ results in
\begin{equation}
Ch^{-4}\gamma r^2 \times \mu^{-1}hr^2\gamma^{-3}\times h^{-\delta'}({\bar\ell}+\mu h)=
C\bigl({\bar\ell}+\mu h\bigr)\mu^{-1}h^{-3-\delta'}r^4\gamma^{-2}
\label{3-12}
\end{equation}
and therefore
\begin{claim}\label{3-13} The total contribution to the remainder of all {\sf(pair,element)} combinations residing in $\cZ_{(r,\gamma)}\subset \cZ^*_{\out, I}$ and satisfying conditions (\ref{3-5}),(\ref{3-6}) and (\ref{3-8}) with $\ell\le{\bar\ell}$ does not exceed (\ref{3-12}).
\end{claim}
It implies

\begin{claim}\label{3-14} The total contribution to the remainder of all {\sf(pair,element)} combinations residing in the \emph{far outer zone\/}
$\cZ^*_{\out, I}\cap \bigl\{\dist(x,\Sigma)\ge \mu^{-1/4}h^{-\delta'}\dist (x,\Lambda)^2\bigr\}$ and satisfying conditions (\ref{3-5}),(\ref{3-6}) and (\ref{3-8}) does not exceed $C\mu^{-1/2}h^{-3}$.
\end{claim}

The same arguments work for ${\bar\cZ}_{(r,\gamma)}$ with $r=\gamma^{1/2}$
resulting in the estimate $C\mu^{-1}h^{-3-\delta'}$ of its contribution to the remainder; summation over $\cZ^*_{\out,II}$ results in expression $o(\mu^{-1/2}h^{-3})$. Therefore

\begin{claim}\label{3-15} The total contribution to the remainder of all {\sf(pair,element)} combinations residing in
$\cZ^*_{\out, II}$ and satisfying conditions (\ref{3-5}),(\ref{3-6}) and (\ref{3-8}) with $\zeta=1$\,\footnote{\label{foot-23} One should pick up here $\zeta_1=1$ because $\hbar_3=\hbar_4$.} does not exceed $C\mu^{-1/2}h^{-3}$.
\end{claim}
I consider the exceptional regular elements (not covered by (\ref{3-14}) in the next subsubsection. The rest of the analysis of this subsection is devoted to detecting and analysis of the irregular elements.

\subsubsection{}
\label{sect-3-1-3}
Note that
\begin{equation}
H_{pn}= \Bigl(\phi_1^{-1}\bigl( V- (2n+1)\mu h f_2 - (2p+1)\mu h f_1\bigr)\Bigr)\circ \Psi
\label{3-16}
\end{equation}
with uniformly smooth map $\Psi$ which does not depend on $V$; this exonerates elliptic arguments of the previous subsubsection and in particular (\ref{3-13}) since perturbation part satisfies $|\partial_n \cB_{pn}|\ll \mu h$; one can consider $n$ and $p$ as continuous parameters.

\begin{proposition}\label{prop-3-1} (i) After reduction in $\cZ^*_{\out, I}$ symbols $W\circ \Psi$, $f_1\circ \Psi$, $f_2\circ \Psi$ satisfy conditions $(\ref{3-5}),(\ref{3-6})$ with $\ell= \gamma r^{-1}$ and indicated above $\eta_\alpha$ and $\zeta$;
furthermore
\begin{equation}
|\nabla_{x_3,\xi_3} f_1\circ \Psi|\asymp \gamma;
\label{3-17}
\end{equation}
(ii) After reduction in $\cZ^*_{\out,II}$ symbols $W\circ \Psi$, $f_1\circ \Psi$, $f_2\circ \Psi$ satisfy conditions $(\ref{3-5}),(\ref{3-6})$ with $\ell= \gamma r^{-1}$;
furthermore
\begin{equation}
|\nabla_{x'',\xi''} f_1\circ \Psi|\asymp \gamma.
\label{3-18}
\end{equation}
\end{proposition}

\begin{proof} Proof follows easily from the reductions. More precise results will be needed and proven later.
\end{proof}

Therefore

\begin{claim}\label{3-19} As $p \ge C_0/(\mu hr)$ symbols $H_{pq}$ satisfy conditions (\ref{3-5}),(\ref{3-6}) and (\ref{3-8} with $\ell= \mu h \gamma p \ge C_0\gamma/r$.
\end{claim}

On the other hand one can prove easily from the reduction that

\begin{claim}\label{3-20} Perturbation symbols $\cB_{pn}$ satisfy conditions (\ref{3-5}) and (\ref{3-6}) with $\ell$ replaced by $\ell'=h^2\gamma^{-2}r^2 (p+1)^2 + \mu^{-1}h\gamma^{-3}r^2 \ll \mu h \gamma (p+1)$.
\end{claim}

Note that if I restrict myself only by $p\le C_0/(\mu hr)$ (condition to be explained later) instead of $p\le C_0/(\mu h\gamma)$, transition from (\ref{3-11}) to (\ref{3-12}) picks up an extra factor $\gamma/r$ and therefore

\begin{claim}\label{3-21} The total contribution to the remainder of all {\sf(pair,element)} combinations with $p\le C_0/(\mu h r)$, residing in
$\cZ^*_{\out, I}$ and satisfying conditions (\ref{3-5}),(\ref{3-6}) and (\ref{3-8}) does not exceed $C\mu^{-1/2}h^{-3}$.
\end{claim}

Now I want to cover indices $p\ge C_0/(\mu h r)$ and elements in $\cZ^*_{\out,I}$ which are not in the far outer zone and prove that their total contribution to the remainder also does not exceed $C\mu^{-1/2}h^{-3}$. To achieve this I plan gain a factor $\gamma r^{-2}$ in the estimates of the previous subsubsection, starting from (\ref{3-10}) by increasing $T_1$ given by (\ref{3-7}) to
\begin{equation}
T_2 = \epsilon r\mu_* \gamma_*^2=\epsilon \mu \gamma^2 .
\label{3-22}
\end{equation}
Note that (\ref{3-7}) appeared because it is a time for which propagation is confined to $B(0,1)$ in the final reduction, for larger $t$ shift with respect to $(x_3,\xi_3)$ would be too large while shift with respect to $(x_4,\xi_4)$ would be still less than $\epsilon$ as $|t|\le T_2$.

However considering propagation after intermediate reduction (i.e. reduction of subsection \ref{sect-2-2}) the shift with respect to $x_3$ would be of magnitude $\mu_*^{-1}\gamma_*^{-2} (\mu h \gamma p) |t|$ which would be observable but less than $1$ as $T_1\le |t|\le T_2$; this propagation just moves to the adjacent ``slabs'' in the intermediate element $B(0,1)$ corresponding to the different final elements.

Then the contribution of such {\sf(pair,element)} to the remainder estimate is given by a modified expression (\ref{3-10})
\begin{equation}
C\hbar_3^{-1}\hbar_4^{-1}\times T_0T_2^{-1} \asymp \mu h^{-1-\delta'} \gamma
\ell^{-1} \asymp \mu h^{-1-\delta'} \gamma
(\mu h \gamma p)^{-1}
\label{3-23}
\end{equation}
which leads to modified expression (\ref{3-12})
\begin{equation}
\mu h^{-1-\delta'} \gamma r^{-2} \times (\mu h)^{-1}\times (\mu h\gamma)^{-1}|\log h| \times r^2\gamma^{-1}
\label{3-24}
\end{equation}
where the second, the third and the fourth factors estimate respectively the
number of indices $n$, the sum of $(\mu h \gamma p)^{-1}$ with respect to $p$ and the number of the final elements with given $(r,\gamma)$.

Then the sum with respect to $\gamma$ could be estimated by by the same expression with $\gamma= \mu^{-1/2-\delta}r^{1/2}$, which is $C\mu^{-1/2-\delta}h^{-3-\delta'}r^{3/2}|\log h|$;
finally summation with respect to $r$ results in $O(\mu^{-1/2}h^{-3})$ as $\delta'>0$ is small enough..

Therefore I arrive to

\begin{claim}\label{3-25} The total contribution to the remainder of all {\sf(pair,element)} combinations with $p\ge C_0/(\mu h r)$, residing in
$\cZ^*_{\out, I}$ does not exceed $C\mu^{-1/2}h^{-3}$.
\end{claim}
which together with (\ref{3-21} concludes the analysis of regular elements; their total contribution to the remainder is $O(\mu^{-1/2}h^{-3})$.

\subsubsection
\label{sect-3-1-4}
Therefore in what follows I am interested only in non-regular elements, with $p\le C_0/(\mu h r)$; then $H_{pn}$ with $n$ violating ellipticity condition satisfies (\ref{3-5}),(\ref{3-6}) with $\ell =\gamma r^{-1}$ in $\cZ^*_{\out, I}$ and with $\ell =\gamma r^{-1}$, $\eta_\alpha=(1,1,1,1)$ in $\cZ^*_{\out, II}$.

First consider elements in $\cZ^*_{\out, I}$ with fixed $p\le C_0/(\mu h r)$ and $n$. After rescaling
\begin{equation}
x_4\mapsto \gamma r^{-2}x_4,\quad\xi_4\mapsto \gamma r^{-2}\xi_4\,\quad
\hbar_4\mapsto \hbar'_4\Def \mu^{-1}h\gamma^{-2}r^2
\label{3-26}
\end{equation}
symbol $\gamma^{-1}r H_{pn}$ becomes uniformly smooth and therefore one can introduce a scaling function
\begin{equation}
\ell= \ell_{pn}(x'',\xi'')=\epsilon r \gamma^{-1} |\nabla_{x'',\xi''}\cA_{pn}|+ {\bar\ell}_0, \qquad {\bar\ell}_0= C\hbar_3^{1/2}h^{-\delta''};
\label{3-27}
\end{equation}
obviously $|\nabla\ell|\le {\frac 1 2}$.
Let us introduce $\ell$-admissible partition of $B(0,1)$.
Then the symbols are still quantizable after rescaling $(x'',\xi'')\mapsto \ell^{-1}(x'',\xi'')$, $\hbar_3\mapsto \ell^{-2}\hbar_3$, $\hbar'_4\mapsto \ell^{-2}\hbar'_4$.

Further, as $\ell \ge 2{\bar\ell}_0$, the propagation speed does not exceed \begin{equation}
v\asymp \mu^{-1}\gamma^{-3}r^2 |\nabla \cA_{pn}|\asymp\mu^{-1}r\gamma^{-2}\ell
\label{3-28}
\end{equation}
in the corresponding partition subelement\footnote{\label{foot-24} I call this partition elements ``subelements'' because actually $B(0,1)$ itself is a rescaled partition element by itself.} (in the coordinates before rescaling of the subelements) and therefore the dynamics started in $B(z,\ell(z))$ remains in $B(z,2\ell(z))$ as $|t|\le T_1$ with
\begin{equation}
T_1 =\epsilon \ell v^{-1} =\epsilon \mu r^{-1} \gamma^2.
\label{3-29}
\end{equation}
On the other hand, the shift $v|t|$ is observable as $|t|\ge T_0$ with $T_0$ defined from the uncertainty principle
\begin{equation}
vT_0 \times \ell \ge \hbar_3 h^{-\delta'}.
\label{3-30}
\end{equation}
Really, if $ |\nabla_{x_3,\xi_3}\cA_{pn}|\asymp \ell\gamma$ then the propagation speed with respect to $(x_3,\xi_3)$ is exactly of magnitude $v$ and the uncertainty principle is due to (\ref{3-30}). Otherwise $|\nabla_{x_4,\xi_4}\cA_{pn}|\asymp \ell\gamma$ and the propagation speed with respect to $(x_4,\xi_4)$ is exactly of magnitude $v\gamma $ and the uncertainty principle is
$\gamma vT_0 \times \ell \ge \hbar_2 h^{-\delta'}r^{-1}$ which is due to (\ref{3-30}) as well.

So one can pick up
\begin{equation}
T_0 =\hbar_3 h^{-\delta'}v^{-1}\ell^{-1} \asymp h^{1-\delta'} \gamma ^{-1} r \ell^{-2}.
\label{3-31}
\end{equation}
Therefore the contribution of this partition subelement to the remainder does not exceed
\begin{equation}
C\hbar_3^{-1}\hbar_4^{\prime\,-1}\ell^4 \times T_0 T_1^{-1}= C\mu h^{-1-\delta'} \gamma^2r^{-2} \ell^2
\label{3-32}
\end{equation}
while in comparison its contribution to the main part of asymptotics does not exceed $C\hbar_3^{-1}\hbar_4^{\prime\,-1}\ell^4 = C\mu^2 h^{-2}\gamma^5 \ell^4$.

On the other hand, the contributions of the subelement with $\ell\le {\bar\ell}_0$ to both the main part and the remainder do not exceed
\begin{equation}
C\hbar_3^{-1}\hbar_4^{\prime\,-1}\ell^4 = C\mu^2 h^{-2} r^{-4} \gamma^5 {\bar\ell}_0^4
\label{3-33}
\end{equation}
which is less than (\ref{3-32}) (with a larger but still arbitrarily small exponent $\delta'$).

Let us notice that for given index $p$ and $\ell$-element operator $\cA_{pn}$ is elliptic as $|n-{\bar n}| \mu h\ge C_0\ell^2\gamma r^{-1}$ for some ${\bar n}$; therefore ellipticity fails for no more than
$\bigl(C_0\ell^2\gamma r^{-1}(\mu h)^{-1}+1\bigr)$ indices $n$.

Then the contribution to the remainder of all the indices $n$ (while index $p$ and $\ell$-subelement remain fixed) does not exceed
\begin{equation}
C\mu h^{-1-\delta'}\gamma^2\ell^2 \times
\bigl(C_0\ell^2\gamma (\mu h r)^{-1}+1\bigr).
\label{3-34}
\end{equation}
As
\begin{equation}
\ell \ge {\bar\ell}_1\Def \epsilon_1\mu ^{1/2}h^{1/2}r^{1/2} \gamma^{-1/2}
\label{3-35}
\end{equation}
this expression is $\asymp Ch^{-2-\delta'}\gamma^3 r^{-1}\ell^4$ and the sum over $\ell$-subelements in $B(0,1)$ with $\ell$ satisfying (\ref{3-35}) results in $Ch^{-2-\delta'} \gamma^3 r^{-1}$.

Then the sum over all such subelements residing in $\cZ_{(r,\gamma)}$ does not exceed
\begin{equation}
Ch^{-2-\delta'} \gamma^3 r^{-1}\times r^2\gamma^{-2} \times r^2\gamma^{-1}=Ch^{-2-\delta'} r^3
\label{3-36}
\end{equation}
where the second factor in the left-hand expression are the number of elements which appeared after rescaling (\ref{3-26}) and the corresponding partition I made in the beginning of subsubsection and the third factor is the number of the ``slabs'' corresponding to the final elements after the second transformation in the ``intermediate'' ball $B(0,1)$ obtained after the first transformation.

Then the sum over the whole $\cZ^*_{\out,I}$ does not exceed $Ch^{-2-\delta'}|\log h|$. So far index $p$ remains fixed; the sum
over $p\le C_0/(\mu h)$ does not exceed $C\mu^{-1}h^{-3-\delta'}|\log h|=O(\mu^{-1/2}h^{-3}$. Thus I arrive to

\begin{claim}\label{3-37} The total contribution to the remainder of all {\sf(pair,subelement)} combinations in zone $\cZ^*_{\out,I}$ with $\ell \ge {\bar\ell}_1$ does not exceed $C\mu^{-1/2}h^{-3}$.
\end{claim}
Therefore only subelements with $\ell\le {\bar\ell}_1$ remain to be considered.

For such subelements expression (\ref{3-34}) becomes
$C\mu h^{-1-\delta'}\gamma^2\ell^2$; however since without non-degeneracy condition the number of such subelements in $B(0,1)$ could be $\asymp \ell^{-4}$, the sum would become singular. To overcome this obstacle
let us redefine the scaling function (\ref{3-27}) in the manner more suitable for the analysis of elements with $\ell\le {\bar\ell}_1$:
\begin{multline}
\ell^*= \ell ^*(x'',\xi'')=\epsilon\min_{p,n}
\Bigl(\gamma^{-2} r^2|\nabla_{x'',\xi''} \cA_{pn}|^2+ r\gamma^{-1} | \cA_{pn}|\Bigr)^{1/2}+ {\bar\ell}_0, \\ \text{with \ \ } {\bar\ell}_0=\hbar_3^{1/2} h^{-\delta'}=\mu ^{-1/2}h^{1/2-\delta''} \gamma^{-3/2}r.
\label{3-38}
\end{multline}
Note that as the minimum in the right-hand expression is achieved for $(p,n)=({\bar p},{\bar n})(x'',\xi'')$, for every index $p$ inequality
$\ell \ge \ell_p=\ell^*+ \epsilon |p-{\bar p}| \mu h r$ holds as the $n=n(p,(x'',\xi''))$ is selected to break ellipticity condition; this index ${\bar n}$ is unique.

Then the contribution of all subelements with $\ell^*\asymp \lambda$ in $B(0,1)$ (obtained after rescaling (\ref{3-26})) to the remainder (after summation with respect to $n$ and $p={\bar p}$) does not exceed the sum with respect to $p$ of expressions (\ref{3-34}):
\begin{equation}
C\mu h^{-1-\delta'} \gamma^2 \sum_{p\ne {\bar p}} \bigl(\lambda+ \epsilon |p-{\bar p}| \mu h r\bigr)^{-2};
\label{3-39}
\end{equation}
I included in this sum only $p$ with $|p-{\bar p}|\ge 1$ leaving \emph{special\/} {\sf($({\bar n},{\bar p})$,subelement)} combinations for a special consideration. Expression (\ref{3-39}) obviously does not exceed
$C\mu ^{-1}h^{-3-\delta'}\gamma^2r^{-2}$ and therefore the contribution of all such {\sf(pair,subelement)} combinations residing in $\cZ_{(r,\gamma)}$ does not exceed
\begin{equation*}
C\mu ^{-1}h^{-3-\delta'}r^{-2}\times r^2\gamma^{-2} \times r^2\gamma^{-1}=
C\mu^{-1}h^{-3-\delta'}r^2\gamma^{-1}
\end{equation*}
with the same origin of the second and the third factors in the left-hand expression as in (\ref{3-36}).

Finally after summation over $\cZ^*_{\out,I}$ one gets $O(\mu^{-1/2}h^{-3})$.
Therefore I arrive to

\begin{claim}\label{3-40} The total contribution to the remainder of all {\sf($(p,n)$,subelement)} combinations in zone $\cZ^*_{\out,I}$ with $\ell \le {\bar\ell}_1$ and $(p,n)\ne ({\bar p},{\bar n})$ does not exceed $C\mu^{-1/2}h^{-3}$.
\end{claim}

\subsubsection
\label{sect-3-1-5}
So, in what follows I need to consider only \emph{special\/} {\sf($({\bar n},{\bar p})$,subelement)} combinations with $\ell^*\le {\bar\ell}_1$. Exactly for these combinations non-degeneracy condition becomes crucial. In the general case however $\cA_{pn}$ with $(p,n)=({\bar p},{\bar n})$ could be ``flat'' and then there would be no difference between estimates of the contribution of such subelement to the remainder and to the main part of the asymptotics which would be $C\hbar_3^{-1}\hbar_4^{\prime\,-1}\asymp \mu ^2h^{-2}\gamma^5 r^{-2}$.
Then the contribution of all such {\sf(pair,subelement)} combinations residing in $\cZ_{(r,\gamma)}$ does not exceed
\begin{equation}
C\mu ^2h^{-2}\gamma^5 r^{-2}\times r^2\gamma^{-2} \times r^2\gamma^{-1}=
C\mu^2h^{-2}\gamma^2 r^2
\label{3-41}
\end{equation}
with the same origin of the second and the third factors in the left-hand expression as in (\ref{3-36}).

Finally after summation over $\cZ^*_{\out,I}$ one gets $O(\mu^{-1/2}h^{-3})$.
Therefore I arrive to

\begin{claim}\label{3-42} The total contribution to the remainder of all {\sf($(p,n)$,subelement)} combinations in zone $\cZ^*_{\out,I}$ with $\ell \le {\bar\ell}_1$ and $(p,n)= ({\bar p},{\bar n})$ does not exceed $C\mu^2h^{-2}$.
\end{claim}

Combining with (\ref{3-37}),(\ref{3-40}) I arrive to

\begin{proposition}\label{prop-3-2} Under condition $(\ref{2-1})$ the total contribution to the remainder of the zone $\cZ^*_{\out,I}$ does not exceed
$C\mu^2 h^{-2}+ C\mu^{-1/2} h^{-3}$.
\end{proposition}

\begin{remark}\label{rem-3-3}
All the arguments leading to proposition \ref{prop-3-2} are applicable in
$\cZ^*_{\out,II}$ with the following minor modifications and simplifications:

\begin{enumerate}[label={(\roman*)}]
\item operators $\cA_{pn}$ are $\hbar_3$-PDOs from the very beginning (so there is no need in (\ref{3-26}) rescaling; there are no ``slabs'' etc, everything is homogeneous;

\item A factor $r^2$ in estimates translates into $\gamma$ which makes life easier; in particular there is no need in the special propagation analysis leading to (\ref{3-25}).
\end{enumerate}
So I arrive to
\begin{enumerate}[label={(\roman*)}]
\item $O(\mu^{-1/2}h^{-3})$ estimate for the contribution into the remainder of all {\sf($(p,n)$,subelement)} combinations in zone $\cZ^*_{\out,II}$ with either $\ell \ge {\bar\ell}_1$ or
$\ell \le {\bar\ell}_1$ and $(p,n)\ne({\bar p},{\bar n})$ and
\item
$O(\mu^2h^{-2}\gamma^3)$ estimate for the contribution into the remainder of all {\sf($(p,n)$,subelement)} combinations in zone $\cZ^*_{\out,II}\cap \{\dist (x,\Sigma)\le \gamma\}$ with $\ell \le {\bar\ell}_1$ and $(p,n)=({\bar p},{\bar n})$.
\end{enumerate}
\end{remark}

In particular:

\begin{proposition}\label{prop-3-4} Under condition $(\ref{2-1})$ the total contribution to the remainder of the zone $\cZ^*_{\out,II}$ does not exceed
$C\mu^2 h^{-2}+ C\mu^{-1/2} h^{-3}$.
\end{proposition}

\subsubsection
\label{sect-3-1-6}
Estimate (\ref{0-3}) $O\bigl(\mu^2 h^{-2}+ \mu^{-1/2} h^{-3}\bigr)$ is the best possible in the general case (see Appendix \ref{sect-A-3}) and it coincides with the best possible $O(\mu^{-1/2}h^{-3})$ as $\mu \le h^{-2/5}$; so $\mu \ge h^{-2/5}$ until the end of this subsection.

However, one can do better under some non-degeneracy conditions. Namely, as before one can estimate the contribution of the special pair {\sf($({\bar p},{\bar n})$,subelement)} to the remainder by
$C\mu h^{-1-\delta'}\gamma^2\ell^2$ as $\ell\ge {\bar\ell}_0$ (see (\ref{3-32})) and by $C\mu^2 r^{-2}\gamma^5$ (see (\ref{3-33})) where now $\ell=\ell^*$ defined by (\ref{3-38}). Therefore the contribution to the remainder of all the special subelements residing in the final element $B(0,1)$ (obtained after rescaling/partition (\ref{3-26}))
does not exceed
\begin{equation}
C\mu h^{-1-\delta'} \gamma^2r^{-2}\int_{\displaystyle\{{\bar\ell}_0\le \ell^*\le {\bar\ell}_1 \}}\ell^{*\,-2}\, dx''d\xi''
\ + \
C\mu^2h^{-2} r^{-4} \gamma^5\int_{\displaystyle\{ \ell^* \le {\bar\ell}_0\}} \, dx''d\xi''
\label{3-43}
\end{equation}
where I remind ${\bar\ell}_0= \mu^{-1/2}h^{1/2-\delta''}\gamma^{-3/2}r$.
${\bar\ell}_1= \mu ^{1/2}h^{1/2}r^{1/2} \gamma^{-1/2}$.

Then the total contribution to the remainder of all the special subelements with $\Psi$-image residing in $\cZ_{(r,\gamma)}\subset \cZ^*_{\out,I} $ to the remainder does not exceed
\begin{multline*}
C\mu h^{-1-\delta'}\gamma^2r^{-2}
\int_{\displaystyle\{{\bar\ell}_0\le \ell^* \circ \Psi^{-1}\le {\bar\ell}_1 \}\cap \cZ_{(r,\gamma)}}\ell^{*\,-2} \,|\det D\Psi|^{-1} dx
\\
C\mu^2h^{-2} r^{-4} \gamma^5\int_{\displaystyle\{ \ell^* \circ\Psi^{-1}\le {\bar\ell}_0\}\cap \cZ_{(r,\gamma)} } \, |\det D\Psi|^{-1} dx;
\end{multline*}
Since
\begin{equation}
|\det D\Psi |\asymp \gamma^4r^{-4}
\label{3-44}
\end{equation}
 one can rewrite this expression as
\begin{multline}
C\mu h^{-1-\delta'}
\int_{\displaystyle\{{\bar\ell}_0\le\ell^*\circ \Psi^{-1}\le {\bar\ell}_1 \} \cap \cZ}\gamma^{-2}r^2\ell^{*\,-2} \, dx
\ + \\
C\mu^2h^{-2}\int_{\displaystyle\{\ell^*\circ\Psi^{-1}\le {\bar\ell}_0\}\cap \cZ} \gamma dx
\label{3-45}
\end{multline}
with $\cZ=\cZ_{(r,\gamma)}$.

I leave to the reader the similar analysis in $\cZ^*_{\out,II}$ leading to the same estimate for the contribution of ${\bar\cZ}_{(r,\gamma)}\subset \cZ^*_{\out,II}$ with $r=\gamma^{1/2}$:

\begin{proposition}\label{prop-3-5} Under condition $(\ref{2-1})$ \nopagebreak

\smallskip
\noindent
(i) The total contribution of the special subelements with $\Psi$-image residing in $\cZ_{(r,\gamma)}\subset \cZ^*_{\out,I} $ does not exceed $(\ref{3-45})$ with $\cZ=\cZ_{(r,\gamma)}$;

\smallskip
\noindent
(ii) The total contribution of the special subelements with $\Psi$-image residing in ${\bar\cZ}_{(r,\gamma)}\subset \cZ^*_{\out,II} $ does not exceed $(\ref{3-45})$ with $\cZ={\bar\cZ}_{(r,\gamma)}$;
where $r=\gamma^{1/2}$.

\smallskip
\noindent
(iii) The total contribution of all the special subelements (with $\Psi$-image residing in $\cZ^*_\out$ does not exceed $(\ref{3-45})$ with $\cZ=\cZ^*_\out$,
\end{proposition}

Note that all elements of $D\Psi $ do not exceed $C\gamma r^{-1}$. Combining with (\ref{3-45}) and (\ref{3-16}) one can conclude that all elements of $D\Psi ^{-1}$ do not exceed $C\gamma ^{-1}r$.
Therefore if $\ell^*$ was defined by (\ref{3-38}) with $H_{pn}$ instead of $\cA_{pn}$ then the following inequality would hold
\begin{equation}
\ell^*\circ\Psi^{-1} \ge L\Def \epsilon |\nabla_\Sigma (V/f_2)|+\min _n |\nabla_\Sigma (V/f_2)-(2n+1)\mu h|.
\label{3-46}
\end{equation}
Then under condition (\ref{0-10}) $\ell^* \asymp 1$ even if $\ell^*$ are defined by $\cA_{pn}$ and the first term in (\ref{3-44}) does not exceed $C\mu h^{-1-\delta'}\gamma^{-1}r^4$ while the second term vanishes and the total contribution of $\cZ^*_{\out,I}$ to the remainder does not exceed $C\mu^{3/2}h^{-1}$.

Therefore I arrive to

\begin{proposition}\label{prop-3-6} Under conditions $(\ref{2-1})$ and $(\ref{0-10})$ the total contribution of $\cZ^*_\out$ to the remainder does not exceed $C\mu^{-1/2}h^{-3}$.
\end{proposition}

Furthermore, under assumption (\ref{3-46}) and condition \ref{0-8-q} one can see easily that (\ref{3-45}) with $\cZ=\cZ^*_\out$ does not exceed $C\mu^{3/2}h^{-1}$ plus
\begin{phantomequation}
\label{3-47}
\end{phantomequation}
\begin{equation}
C\mu^2h^{-2}\int \gamma {\bar\ell}_0^q \,dx
\tag*{$(3.47)_q$}\label{3-47-q}
\end{equation}
which is
\begin{equation*}
\asymp \iint_\Omega J_q \,\gamma^{-1}d\gamma\cdot r^{-1}dr, \qquad
J_q= h^{-2\delta''}\mu^2h^{-2} (\mu^{-1}hr^2\gamma^{-3})^{q/2}r^2\gamma^2
\end{equation*}
with integral taken over $\Omega=\{\mu^{\delta'-1/2}r^{1/2}\le \gamma \le r^2\le 1\}$.

Obviously $(\ref{3-47})_2$ does not exceed $J_2$ calculated as $\gamma=\mu^{\delta'-1/2}$, $r=1$ which is $O(\mu^{-1/2}h^{-3}$.
Further, $(\ref{3-47})_1$ does not exceed $J_1$ calculated as $\gamma=r=1$ which is $C\mu^{3/2}h^{-3/2-\delta ''}$ with arbitrarily small $\delta''>0$.

So, one arrives to the following statement (still under (\ref{3-46})

\begin{proposition}\label{prop-3-7} Let condition $(\ref{2-1})$ be fulfilled. Then

\smallskip
\noindent
(i) Under condition $(\ref{0-9})_2$ the total contribution to the remainder of the zone $\cZ^*_\out$ with does not exceed $C\mu^{-1/2} h^{-3}$;

\smallskip
\noindent
(ii) Under condition $(\ref{0-9})_1$ the total contribution to the remainder estimate the zone $\cZ^*_\out$ does not exceed $C\mu^{3/2}h^{-3/2-\delta} + C\mu^{-1/2} h^{-3}$ with arbitrarily small $\delta''>0$.
\end{proposition}

\begin{proof} To get rid of assumption (\ref{3-26}) one needs to take in account perturbation term $\cB$ in $\cA_{pn}$; then $\ell^*\circ\Psi^{-1}\ge L- |\nabla_\Sigma \cB_{pn}\circ \Psi^{-1}|$.
Thus it would be enough to replace ${\bar\ell}_0$ by ${\bar\ell}_0 + |\nabla_\Sigma \cB_{pn}\circ \Psi^{-1}|$ in (\ref{3-47}).

There are two leading terms in $\cB_{pn}$; the first one is
$\mu^{-2}r^2\gamma^{-4}((2p+1)\mu h\gamma)^2$; since $p\le C_0/(\mu h r)$ this term does not exceed $C\mu^{-2}\gamma^{-2}$. The second term is $O(\mu^{-1}hr^2\gamma^{-3})$ but looking its origin one can notice that it must contain $|\nabla_\Sigma (V/f_2)|^2$; then
\begin{equation*}
|\nabla_\Sigma \cB_{pn}\circ \Psi^{-1}|\le L_1\Def C\mu^{-2}\gamma^{-3}r.
\end{equation*}
Plugging ${\bar\ell}_0= \mu^{-2}\gamma^{-3}r$ into $(\ref{3-47})_1$ results in $O(\mu^{1/2}h^{-2})$.
\end{proof}

\begin{remark}\label{rem-3-8} Note that in the extra terms $O(\mu^2 h^{-2})$ and $O(\mu^{3/2}h^{-3/2-\delta''})$ appearing in the general case and under condition \ref{0-8-q} the main contributor is far outer zone where $\gamma\approx 1$ and $r\approx 1$ while in the sharp estimate $O(\mu^{-1/2}h^{-3})$ the main contributor is zone where $\gamma \approx \mu^{-1/2}$ and $r\approx 1$;

\smallskip
\noindent
(ii) Actually one can greatly improve remainder estimate in $(\ref{0-8})_1$ case using more refined partition-rescaling technique like in the proof of proposition \ref{sect-4-5}. I believe that it can be even brought to $O(\mu^{-1/2}h^{-3})$. However I think that really interesting are only general ($q=0$) and generic $(q=3$) cases.
\end{remark}

\subsection{Estimates in the near outer zone}
\label{sect-3-2}

After deriving remainder estimates in $\cZ^*_\out$ I need to consider the remaining part of $\cZ_\out$ (which will be the \emph{ near outer zone}
\begin{equation}
\cZ'_\out= \cZ_\out \cap \{\dist(x,\Sigma)\le \mu^{-1/2}h^{-\delta}\}
\label{3-48}
\end{equation}
and the \emph{inner zone}. However I also introduce and start from the \emph{inner core} where I am able to apply pretty non-sophisticated approach and still derive proper estimates, just to keep $r$ more disjoint from $0$ to be able to use canonical form (\ref{2-38}).

\begin{figure}[h!]
\centering
%
%
\includegraphics{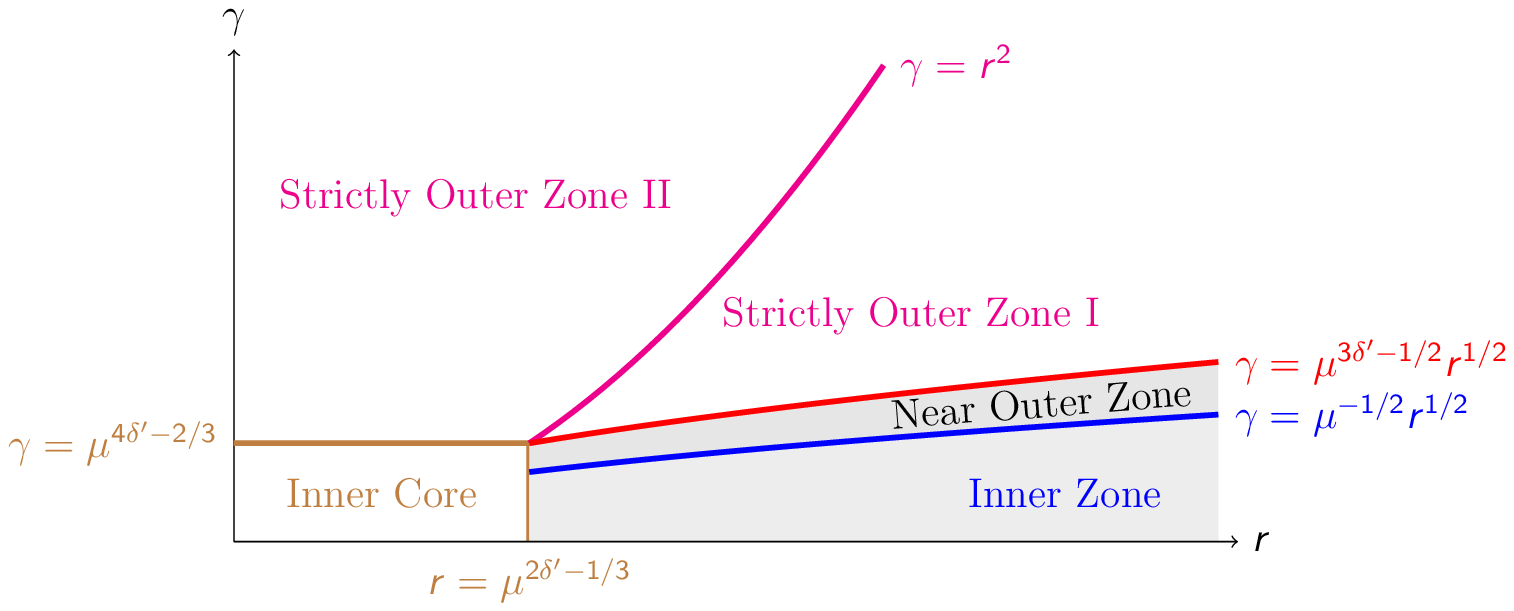}
\caption{\label{fig-3} Zones.}
\end{figure}

\subsubsection
\label{sect-3-2-1}

In this subsection my goal is to derive remainder estimate $O(\mu^{-1/2}h^{-3}+\mu^2h^{-2})$. Therefore my approach here would be pretty unsophisticated. Note first that

\begin{claim}\label{3-49}
The contribution of the domain
\begin{equation}
\bigl\{\dist (x,\Lambda)\le r, \ \dist(x,\Sigma)\le \gamma, \ |Z_1|\le \rho\bigr\}
\label{3-50}
\end{equation}
to the remainder does not exceed $C\mu h^{-3}\gamma r^2\rho^2$ (while the main part of asymptotics is given by (\ref{3-1}) with $T_0\ge Ch|\log h|$);
\end{claim}
therefore considering \emph{inner core} given by (\ref{3-50}) with $r=\mu^{-1/3+2\delta'}$, $\gamma=\mu^{-2/3+4\delta'}$ one gets $O(\mu^{-1/3+6\delta'}h^{-3})$.

This remainder estimate is very rough and should be improved despite rather unambitious goal of this subsubsection. However let me note first that considering the strip $\{|x_1|\le \mu^{-1/2+3\delta'}\}$ and thus taking $r=1$, $\rho=1$ and $\gamma=\mu^{-1/2+3\delta'}$ one would get the remainder estimate $O(\mu^{1/2+3\delta'}h^{-3})$ which is $O(\mu^2h^{-3})$ as as $\mu \ge h^{-2/3-\delta }$. Further, the contribution of this zone to the asymptotics is $O(h^{-4}\gamma r^2)$ and thus taking $r=1$, $\rho=1$ and $\gamma=\mu^{3\delta'-1/2}$ one would get $O(\mu^{-1/2+3\delta'}h^{-4})$ which is $O(\mu^2h^{-2})$ as $\mu \ge h^{\delta-4/5}$. Therefore in what follows until the end of subsubsection one can assume that
\begin{equation}
h^{-\delta}\le \mu \le h^{\delta-1}.
\label{3-51}
\end{equation}

Note first that
\begin{equation}
| F_{t\to h^{-1}\tau} \chi_T(t)\Gamma (u Q_y^\dag)|\le Ch^s
\label{3-52}
\end{equation}
as $Q$ is supported in ${\bar Z}_{(\gamma)}\Def \{|x_1|\le \gamma\}$ and $T_0\le T\le T_1$ with $T_0= C h|\log h|$ and $T_1=\epsilon \mu^{-1}$; this statement is empty in the case of the \emph{very strong magnetic field\/}
\begin{equation}
\epsilon (h|\log h|)^{-1} \le \mu \le C h^{-1}
\label{3-53}
\end{equation}
which is not the case now.

However looking at the propagation in direction $\bK_1$ as $|Z_1|\asymp \rho $ one can easily notice that the propagation speed is $\asymp \rho $ and the shift for time $T$ is $\asymp \rho T$; this shift is observable as $\rho^2 T\ge C h|\log h|$. Therefore (\ref{3-51}) holds as $T_0\le T\le T_1$ with $T_0= C\rho^{-2} h|\log h|$ and $T_1=\epsilon \mu^{-1}\gamma^{-1}$. Under assumption $\gamma \ge C_0\mu^{-1/2}$ this short evolution is contained in $\epsilon_1\gamma$-vicinity of the original point\footnote{\label{foot-25} One can improve this statement and the estimate following from it.}.

These two intervals $[Ch|\log h, \epsilon\mu^{-1}]$ and $[C\rho^{-2}h|\log h|, \epsilon \mu^{-1}\gamma^{-1}]$ overlap as
\begin{equation}
\rho\ge (C\mu h|\log h|)^{1/2}
\label{3-54}
\end{equation}
and therefore I arrive to

\begin{claim}\label{3-55} Under assumptions $(\ref{3-51})$ and $(\ref{3-54})$
estimate $(\ref{3-52})$ holds as $T\in [T_0, T_1]$ with $T_0=Ch|\log h|$, $T_1= \epsilon \min \bigl(\mu ^{-1}\gamma^{-1},\mu ^{-1/2}\bigr)$.
\end{claim}

Therefore the contribution of the domain
\begin{equation}
\bigl\{\dist (x,\Lambda)\le r, \ \dist(x,\Sigma)\le \gamma, \ |Z_1|\ge C(\mu h|\log h|)^{1/2}\bigr\}
\label{3-56}
\end{equation}
to the remainder does not exceed $C\mu h^{-3}r^2 \gamma (\gamma + \mu^{-1/2})$.

Combining with (\ref{3-49}) I arrive to

\begin{proposition}\label{prop-3-9}
Contribution of domain $(\ref{3-50})$ with $\rho\asymp 1$ to the remainder does not exceed
\begin{equation}
C\mu h^{-3} r^2\gamma (\gamma + \mu^{-1/2}) +C\mu ^2h^{-2} \gamma r^2 |\log h|
\label{3-57}
\end{equation}
while the main part of asymptotics is given by $(\ref{0-12})$ with (any) $T_0\ge C\mu h|\log h|$.
\end{proposition}

\begin{remark}\label{rem-3-10} Estimate (\ref{3-57}) is sufficient for my limited goal here. In particular it covers the inner core giving remainder estimate $O(\mu^{-1/2}h^{-3}+\mu^{2/3+\delta }h^{-2})$ which is not only better than $O(\mu^{-1/2}h^{-3}+\mu^2h^{-2})$ but than $O(\mu^{-1/2}h^{-3}+\mu^{3/2}h^{-3/2})$.
\end{remark}

So, I need to cover the \emph{near outer\/} and \emph{inner\/} zones; in both $\gamma \le r^2$ and therefore I can use canonical form (\ref{2-38}).
Even without it, according to section \ref{sect-1} the drift speed is
$\asymp \mu^{-1} \rho ^2 r \gamma^{-2}$ as $\rho^2 r \ge C_0 \gamma$.

Then the shift for time $T$ is observable as
$\mu^{-1} \rho ^2 r \gamma^{-2} T\times \rho \ge Ch|\log h|$; plugging $T=\mu^{-1}\gamma^{-1}$ one gets $(\ref{3-58})_2$ below and therefore \emph{in $(\ref{3-52})$ one can upgrade $T_1=\epsilon \mu^{-1}\gamma^{-1}$ to $T_1=\epsilon \mu \gamma^2$ provided}
\begin{phantomequation}\label{3-58}\end{phantomequation}
\begin{equation}
\rho^2 r \ge c \gamma,\qquad \rho^3 r \ge C\mu^2 \gamma^3 h|\log h|.
\tag*{$(3.58)_{1,2}$}
\label{3-58-12}
\end{equation}
However, the contribution to the remainder of the part of $\cZ'_\out$ where
$T_1=\epsilon \mu \gamma^2$ is $O(\mu^{-1/2}h^{-3})$. Therefore \emph{one needs to consider only part of $\cZ'_\out$ where \ref{3-58-12} is violated}. At this moment one can ignore subzone $\{\rho \le C (\mu h|\log h|)^{1/2}\}$ and therefore one can upgrade $T_1=\epsilon \mu^{-1}$ to $T_1=\epsilon \mu^{-1}\gamma^{-1}$ with no penalty.

Moreover, using $O(\mu h^{-3}r^2\gamma^2\rho^2)$ estimate (with $r=1, \gamma=h^{-1/2+3\delta'})$ one can see easily that the contribution to the remainder of the zone where condition $(\ref{3-58})_2$ is violated does not exceed $C\mu^{1/3}h^{-7/3-\delta}$ which $o(\mu^{-1/2}h^{-3}+\mu^2 h^{-2})$.

Therefore one needs to consider only zone where condition $(\ref{3-58})_1$ is violated; again due to $O(\mu h^{-3}\gamma^2 \rho^2r^2)$ estimate its contribution to the remainder does not exceed $C\mu^{-1/2}h^{-3-\delta}$ which is only marginally short of what I want (and only the case $\mu \le h^{-2/5-\delta}$ needs to be considered).

The proper estimate of the zone where condition $(\ref{3-58})_1$ is violated can be done easily on the base of proposition \ref{prop-1-6}. Really, the drift speed is $\asymp \mu^{-1}(\rho +\rho') \Delta r \gamma^{-2}$ with
$\rho' = O(\gamma r^{-1})$, $\Delta=|\rho -\rho'|$ and therefore one can upgrade $T_1=\epsilon \mu^{-1}\gamma^{-1}$ to $T_1=\epsilon \mu \gamma^2$ as
\begin{equation*}
\mu^{-1} (\rho +\rho') \Delta r \gamma^{-2} \times \mu^{-1}\gamma^{-1}\times \Delta \ge Ch|\log h|
\end{equation*}
with $\rho' =O(\gamma^{1/2}r^{-1/2})$. Therefore one can take
\begin{equation*}
{\bar\Delta } =
C\mu \gamma^{3/2} r^{-1/2} (h|\log h|)^{1/2} (\rho +\rho')^{-1/2};
\end{equation*}
then as $|\rho-\rho'|\ge {\bar\Delta}$ one can take $T_1=\epsilon \mu \gamma^2$.

On the other hand, the contribution of the zone
$\bigl\{(x,\xi): |\rho -\rho'|\le {\bar\Delta}\bigr\}\cap \cZ_{(r,\gamma)}$ to the remainder does not exceed
\begin{equation*}
C\mu h^{-3}r^2 \rho' \gamma r^2{\bar\Delta}= C\mu^2 h^{-3}r^2 \gamma^{5/2}
r^{3/2} (h|\log h|)^{1/2} \rho^{\prime\,1/2}.
\end{equation*}
Since $\rho' \le c(\gamma/r)^{1/2}$ the latter expression does not exceed
$C\mu^2 h^{-3} \gamma^{11/4} r^{5/4} (h|\log h|)^{1/2}$.
In particular, the contribution to the remainder of the part of $\cZ_{(\gamma)}$ where $(\ref{3-58})_1$ is violated does not exceed
$C\mu^2 h^{-5/2} \gamma^{11/4}|\log h|^{1/2} $ and therefore the
contribution of the corresponding part of $\cZ'_\out$ does not exceed
$C\mu^{5/8} h^{-5/2-\delta}$.

So I arrive to

\begin{proposition}\label{prop-3-11} Under condition $(\ref{2-1})$ the total contribution to the remainder of zone $\cZ'_\out$ does not exceed $C(\mu^{-1/2}h^{-3}+\mu^2h^{-2})$.
\end{proposition}

\subsubsection
\label{sect-3-2-2}
Now I need to improve the previous estimate under non-degeneracy condition.
However, let us summarize the remainder estimate I actually derived:
\begin{enumerate}[label=(\roman*)]
\item The inner core intersected with $\{\rho \ge C(\mu h|\log h|)^{1/2}\}$ contributed $O(\mu h^{-3}{\bar\gamma}_0^2{\bar r}_0^2)= o(\mu^{-1/2}h^{-3})$
where ${\bar\gamma}_0=\mu^{-2/3+4\delta'}$, ${\bar r}_0=\mu^{-1/3+2\delta'}$ ;
\item $\cZ'_\out\cap \{\rho \ge C(\mu h|\log h|)^{1/2}\}$ contributed
$O(\mu h^{-3}{\bar\gamma}_1^2 (\mu h{\bar\gamma}_1)^{2/3})=
O(\mu^{1/3}h^{-7/3-\delta''})$ where ${\bar\gamma}_1=\mu^{-1/2+4\delta'}$;
\item $\{|x_1|\le {\bar\gamma}_1\}\cap\{\rho \le C(\mu h|\log h|)^{1/2}\}$ contributed $O(\mu^2h^{-2}{\bar\gamma}_1|\log h|)=O(\mu^{1/2+\delta''}h^{-2})$;
\end{enumerate}
each of estimates in (ii),(iii) are less than $O(\mu^{-1/2}h^{-3})$ as $\mu \le h^{-1/2+\delta}$. Therefore in what follows one can assume that
\begin{equation}
h^{-1/2+\delta}\le \mu \le C h^{-1}.
\label{3-59}
\end{equation}
The better analysis is based on the canonical form (\ref{2-38}). However, this requires the inner core to be treated separately:

\begin{proposition}\label{prop-3-12} Under condition $(\ref{3-59})$ the contribution of
\begin{equation*}
\cZ^0_{(r)}\Def \bigl\{\dist (x,\Sigma)\le r^2,\dist (x,\Lambda)\le r\bigr\}
\end{equation*}with $r\ge \mu^{-1/3+2\delta'}$ does not exceed $C\mu h^{-3}r^6$ while the principal part of asymptotics is given by $(\ref{0-12})$ with $T={\bar T}\Def Ch|\log h|$.
\end{proposition}

\begin{proof} Looking at the canonical form (\ref{2-17}) and scaling $x''\mapsto x''r^{-1}$, $x_1\mapsto x_1 r^{-1}$, $\mu\mapsto \mu_*=\mu r^3$, $h\mapsto h_*=hr^{-1}$ one gets 2-dimensional magnetic Schr\"odinger operator with
parameters $\mu_*,h_*$, which also is 1-dimensional $\hbar_4$-PDO with respect to $x_4$ with $\hbar_4=\mu^{-1}hr^{-2}$; for such Schr\"odinger operators one can prove easily that the main part of asymptotics is given by Weyl formula and is of magnitude $\hbar_4^{-1}h_*^{-2}= \mu h^{-3}r^4$
while the remainder estimate gains the factor $\mu_*h_*=\mu h r^2$ and thus is of magnitude $\mu^2h^{-2}r^6$.
Really, if $\mu_*=1$ it would follow from the standard results (and no non-degeneracy condition is needed); as $\mu_*\ge 1$ one needs just apply extra rescaling $x'\mapsto x'\mu_*$ which is also rather standard.
\end{proof}

\begin{corollary}\label{cor-3-13} Therefore even without nondegeneracy conditions contribution of the (extended) inner core $\cZ^0_{(r)}$ with $r=\mu^{-1/4}$ to the remainder does not exceed $C\mu^{-1/2}h^{-3}$.
\end{corollary}

Thus in what follows I can apply canonical form (\ref{2-38}).
Let us fix some index $n$ and some final $B(0,1)$ element where in contrary to subsection \ref{sect-3-1} final means only after reduction to (\ref{2-38}). This element is either \emph{regular\/}, when on the original $B({\bar x}',\epsilon r)$ element
\begin{equation}
|Vf_2^{-1}-(2n+1)\mu h |+|\nabla_{x'} \bigl(Vf_2^{-1}-(2n+1)\mu h\bigr)|^2 \asymp L^2\ge r^2
\label{3-60}
\end{equation}
or \emph{singular\/} when on the original $B({\bar x}',\epsilon r)$ element
\begin{equation}
|Vf_2^{-1}-(2n+1)\mu h |+|\nabla_{x'} \bigl( V f_2^{-1}-(2n+1)\mu h\bigr)|^2 \le r^2.
\label{3-61}
\end{equation}
Since analysis of the regular elements is similar to the analysis of the singular elements but simpler, I consider only singular elements, leaving regular elements to the reader.

Let us introduce the scaling function with respect to $(x_3,x_4,\xi_4)$
\begin{equation}
\ell =\ell_n= \epsilon r^{-1}|\nabla_{x_3,x_4,\xi_4} Vf_2^{-1}|+{\bar\ell}
\label{3-62}
\end{equation}
where ${\bar\ell} \ge C\gamma r^{-1}|\log h|$ will be chosen later. With respect to $|Z_1^\#|$ I use scaling function $\rho =\epsilon \bigl(|Z_1|+ r^{1/2}\ell\bigr)$.  Then  ellipticity of $\cA_n$ is not broken unless $\bigl(Vf_2^{-1}-(2n+1)\mu h\bigr)\asymp \rho^2$ or $\rho \le cr^{1/2}\ell$.

Then it follows from the standard theory that for $u_n$ holds as $Q$ is supported in the corresponding $\ell$-partition element, $T\in [T_0,T_1]$, $\tau\le Cr \ell^2$ with
\begin{align}
&T_0 = Ch \bigl(\rho^2+r \ell^2)^{-1} |\log h|,\label{3-63}\\
&T_1 = \epsilon \mu^{-1}\gamma^{-1},\label{3-64}\\
&\rho +r^{1/2}\ell \ge C{\bar\rho}_0 \Def C (\mu h \gamma  |\log h|)^{1/2} \label{3-65}\\
\intertext{and moreover, under these conditions}
&|F_{t\to h^{-1}\tau}{\bar \chi}_T(t) \Gamma (u_n Qy^\dag)|\le C \mu h^{-2} r^3\ell^3\gamma\label{3-66}
\end{align}
Therefore

\begin{claim}\label{3-67}
While the contribution of this $(\ell;\rho)$-subelement into asymptotics does not exceed
$C\mu h^{-3} r^3\bigl(\rho^2+r\ell^2\bigr)\ell^3\gamma$ its contribution to the remainder so far does not exceed 
\begin{equation*}
C\mu h^{-3} r^3\bigl(\rho^2+r\ell^2\bigr)\ell^3\gamma \times T_0T_1^{-1}=C\mu^2 h^{-2}r^3 \ell^3\gamma^2.
\end{equation*}
\end{claim}
Note that for each $(\ell;\rho)$-partition subelement in $(x_3,x_4,\xi_4;\xi_1,\xi_3)$ the number of indices $n$ satisfying (\ref{3-62}) and also violating ellipticity does not exceed
\begin{equation}
M\Def \Bigl(C_0(\mu h)^{-1}(\rho^2+r \ell^2) +1\Bigr)
\label{3-68}
\end{equation}
Also note that under condition \ref{0-8-q} the total number of $\ell$-subelements is $O(\ell^{q-3})$. Then the total contribution of all $\ell$-elements to the asymptotics (as magnitudes of $r,\gamma$ are given) does not exceed
\begin{equation}
C\mu h^{-3}\gamma r^{(2-q)_+} (\rho^2+r \ell^2) \times \Bigl(C_0(\mu h)^{-1}(\rho^2+r\ell^2) +1\Bigr) \times \ell^q
\label{3-69}
\end{equation}
which one can  sum nicely to its value at maximal $\ell=1$, $\rho=1$ and then to maximal $r=1$, $\gamma=\mu^{\delta-1/2}$ resulting $Ch^{-4}\mu^{\delta-1/2}$. 

On the other hand, the total contribution of $\ell$-subelements to the remainder so far does not exceed
\begin{equation}
C\mu^2 h^{-2}\gamma^2 r^{(2-q)_+} \times \Bigl(C_0(\mu h)^{-1}(\rho+r\ell)^2+1\Bigr) \times \ell^q
\label{3-70}
\end{equation}
which one can also sum nicely according to the same principles but the sum is $C\mu h^{-3}\gamma^2$ not estimated that well and $\ell,r$ close to 1 are problematic. So I want to increase $T_1$.

I remind that  $W_n=\omega \bigl(Vf_2^{-1}-(2n+1)\mu h\bigr)$ with $\omega=(\sigma\sigma_1)^{-1}f_2$ and here (see proposition \ref{prop-2-8})
$\omega \sim |\{Z_1,f_1\}|^{-2/3}$ rescaled;
therefore for ``small'' $\rho$ calculation of $\nabla W_n$ and $\omega\nabla \bigl(Vf_2^{-1}-(2n+1)\mu h\bigr)$ are equivalent modulo $\omega \rho^2 r^{-1}$. Note that 

\begin{claim}\label{3-71} As $\bigl(Vf_2^{-1}-(2n+1)\mu h\bigr)\asymp \rho^2$ and $\ell \ge C\rho^2r^{-1}$ the propagation speed with respect to $x_4$ does not exceed $C\mu^{-1}r^{-1}\ell$ and is of this magnitude as $|\partial_{\xi_4}W_n|\asymp \ell$; then the shift with respect to $x_4$ is observable as
\begin{equation}
\mu^{-1}r^{-1}\ell T\times \ell \ge C\mu^{-1}r^{-2}h|\log h|;
\label{3-72}
\end{equation}
\end{claim}

\begin{claim}\label{3-73} As $\bigl(Vf_2^{-1}-(2n+1)\mu h\bigr)\asymp \rho^2$ and $\ell \ge C\rho^2r^{-1}$ the propagation speed with respect to $\mu^{-1}hD''$ does not exceed $C\mu^{-1}r^{-1}\ell$ and is of this magnitude as
$|\partial_{x''}W_n|\asymp \ell$; then the shift with respect to $\mu^{-1}hD''$ is observable under the same condition;
\end{claim}

\begin{claim}\label{3-74} Finally, in the short-term propagation (as $T\le \epsilon \mu^{-1}\gamma^{-1}$) the propagation speed with respect to $(x_1,x_3)$ is of magnitude $|Z_1| \asymp \rho$ and the shift is observable as $\rho T\times \rho \ge Ch|\log h|$. \end{claim}

On the other hand,

\begin{claim}\label{3-75} Let $\rho^2 \ge C\gamma$; then the drift speed (in the rescaled coordinates) with respect to $x''$ is of magnitude $C\mu^{-1}\gamma^{-2}\rho^2$ and the shift for time $T$ is observable as
\begin{equation}
\mu^{-1}\gamma^{-2}\rho^2 T \times \rho \ge Chr^{-1}|\log h|.
\label{3-76}
\end{equation}
\end{claim}

Therefore one can increase $T_1= \epsilon \mu^{-1}\gamma^{-1}$ provided one of the conditions (\ref{3-72}),(\ref{3-76}) holds with $T=T_1$.
Plugging $T_1= \epsilon\mu^{-1}\gamma^{-1}$ into (\ref{3-72}), (\ref{3-76}) I get
\begin{align}
&\rho \ge {\bar\rho}\Def C\bigl(\mu^2 h \gamma^3 r^{-1}|\log h|\bigr)^{1/3}+ C\gamma^{1/2},\label{3-77}\\
&\ell \ge {\bar\ell}\Def C\bigl(\mu h \gamma r^{-1}|\log h|\bigr)^{1/2}+ C\rho^2r^{-1}
\label{3-78}
\end{align}
respectively. 

Now consider increased $T_1$. I claim that

\begin{claim}\label{3-79}
Under condition (\ref{3-78}) one can take $T=\epsilon \mu \gamma^2 \ell$.
\end{claim}

Really, then magnitude of $|\nabla_{x_1,x_4,\xi_4}W_n|$ changes by no more than $\epsilon \ell$ as long as $\rho\le \epsilon_1$; 
however, as $\rho \ge \epsilon$ condition (\ref{3-77}) is fulfilled and one can take $T_2=\epsilon \mu \gamma^2$;

Therefore the contribution to the remainder of the zone with the given magnitudes of $\rho, \ell, \gamma,r$ satisfying (\ref{3-78}) does not exceed expression (\ref{3-69}) multiplied by
$T_0T_2^{-1}$; the resulting expression does not exceed 
\begin{multline}
C\mu h^{-3}\gamma r^{(2-q)_+} (\rho^2+r \ell^2) \times \Bigl(C_0(\mu h)^{-1}(\rho+r\ell)^2 +1\Bigr) \times \ell^q  \times\\ \hfil h (\rho^2+r\ell^2)^{-1} \times 
\mu ^{-1}\gamma^{-2} \ell^{-1}=\\
C\mu^{-1} h^{-3}\gamma^{-1} r^{(2-q)_+}  \Bigl((\rho+r\ell)^2 +\mu h\Bigr)\ell^{q-1};
\label{3-80}
\end{multline}
then summation with respect $\ell,\rho$ results in 
$C\mu^{-1} h^{-3}\gamma^{-1} r$ as $q>1$ (then one should take $q$ slightly larger than 1). As $q=1$ one gets 
$C\mu^{-1/2} h^{-3} (1+\mu h|\log h|)$ as provided summation is taken over $\rho \le c(\ell + |\log h|^{-1})$; but the latter case is covered by arguments below.

\begin{claim}\label{3-81} Under conditions (\ref{3-77}) and $\rho\ge r^{1/2}\ell$ one can take $T_2=\epsilon \mu \gamma^2 \rho r^{-1/2}$.
\end{claim}
Really, in this case due to bounds for propagation and drift speeds this guarantees that both $|\nabla_{x_1,x_4,\xi_4}W_n|$ and $W_n\asymp |Z_1|^2$ change by no more than $r^{-1/2}\rho$, $\epsilon_1\rho^2$ respectively during the evolution.

Now instead of (\ref{3-80}) one gets 
\begin{multline}
C\mu h^{-3}\gamma r^{(2-q)_+} (\rho^2+r \ell^2) \times \Bigl(C_0(\mu h)^{-1}(\rho+r\ell)^2 +1\Bigr) \times \ell^q  \times\\ \hfil  h (\rho^2+r\ell^2)^{-1} \times 
\mu ^{-1}\gamma^{-2} \rho^{-1}\le \\
C\mu^{-1} h^{-3}\gamma^{-1} r^{(2-q)_+}  \Bigl((\rho+r\ell)^2 +\mu h\Bigr)\ell^q \rho^{-1}
\label{3-82}
\end{multline}
and summation with respect to $\rho$, $\ell$, $\gamma$, $r$ results in $C\mu^{-1/2}h^{-3}$. Therefore

\begin{claim}\label{3-83}
The total contribution to the remainder of elements in $\cZ''_out$ satisfying (\ref{3-65}) and one of conditions (\ref{3-77}), (\ref{3-78}) does not exceed $C\mu^{-1/2}h^{-3}$.
\end{claim}

\subsubsection{}
\label{sect-3-2-3}
So, I am left with two types of elements, which are not covered by the above arguments:

\begin{claim}\label{3-84}
With  ${\bar\rho}_0 \le \rho \le {\bar\rho}= C(\mu^2h\gamma^3r^{-1}|\log h|)^{1/3}+C\gamma^{1/2}$ and $ \ell\le {\bar\ell}\asymp r^{-1}{\bar\rho}_0+\rho^2 r^{-1}$;
\end{claim}

\begin{claim}\label{3-85}
With   $\rho \le {\bar\rho}_0= C(\mu h\gamma)^{1/2}$ and $\ell\le {\bar\ell}\asymp r^{-1/2}{\bar\rho}_0$.
\end{claim}

For (\ref{3-84})-type elements $r\ell^2 \le \rho^2$ and expression (\ref{3-69}) becomes 
\begin{equation*}
C\mu h^{-3}\gamma r^{(2-q)_+} \rho^2 \times \Bigl((\mu h)^{-1}\rho^2 +1\Bigr) \times {\bar\ell}^q;
\end{equation*}
multiplying it by $T_0T_1^{-1}= Ch\rho^{-2}|\log h| \cdot \mu \gamma$ one gets
\begin{equation*}
C\mu ^2h^{-2}\gamma^2 r^{(2-q)_+}  \times \Bigl((\mu h)^{-1}\rho^2 +1\Bigr) \times {\bar\ell}^q|\log h|
\end{equation*}
which after summation over $\rho$ becomes
\begin{equation}
C\mu ^2h^{-2}\gamma^2 r^{(2-q)_+}  \times \Bigl((\mu h)^{-1}{\bar\rho}^2 +1\Bigr) \times {\bar\ell}^q|\log h|.
\label{3-86}
\end{equation}
Meanwhile for (\ref{3-85})-type elements expression (\ref{3-69}) becomes 
\begin{equation}
C\mu h^{-3}\gamma r^{(2-q)_+} {\bar\rho}_0^2  {\bar\ell}^q.
\label{3-87}
\end{equation}
Let us plug ${\bar\rho}$, ${\bar\rho}_0$, ${\bar\ell}$ and $\gamma$; one can take 
\begin{gather*}
\gamma=\mu^{-1/2}r^{1/2},\ {\bar\rho}_0=(\mu hr^{1/2})^{1/2},  {\bar\rho}=(\mu^{1/2} h r^{1/2})^{1/3}+\mu^{-1/4}r^{1/4}\ \text{ and }\\
{\bar\ell}=(\mu ^{1/2}hr^{-1/2})^{1/2}+(\mu^{1/2}h )^{2/3}r^{-2/3}+\mu^{-1/2}r^{-1/2},
\end{gather*}
 and compensate this by multiplication the result by $h^{-\delta''}$ with an arbitrarily small exponent $\delta''>0$.

Take first ${\bar\ell}={\bar\ell}_1\Def (\mu^{1/2}hr^{-1/2})^{1/2}$;  
 note that then in (\ref{3-86}), (\ref{3-87})  factor $r^{-q/4}$ coming from ${\bar\ell}^q$ is more than compensated by the factor $r$ coming from $\gamma^2$ or $\gamma {\bar\rho}_0^2$ respectively. Therefore summation with respect to $\gamma,r$ results in the same expressions calculated as $r=1$, $\gamma=\mu^{-1/2}$, ${\bar\ell}={\bar\rho}_0=(\mu^{1/2}h)^{1/2}$ and
${\bar\rho}=(\mu^{1/2} h )^{1/3}+\mu^{-1/4}$:
\begin{align}
&C h^{-3}\Bigl((\mu^{1/2}h)^{2/3}+\mu^{-1}+\mu h\Bigr) (\mu ^{1/2}h)^{q/2}\times h^{-\delta''},\label{3-88}\\
&C \mu^{1/2}h^{-3}(\mu^{1/2}h)^{(2+q)/2}\times h^{-\delta''}\label{3-89}.
\end{align}
One can see easily that both these expressions do not exceed

\begin{enumerate}[label=({\roman*})]
\item $C\bigl(\mu^{-1/2}h^{-3}+\mu^{3/2}h^{-3/2}\bigr)$ as $q=1$;
\item $C\bigl(\mu^{-1/2}h^{-3}+\mu^{1/2+\delta''}h^{-2}\bigr)$ as $q=2$;
\item $C\mu^{-1/2}h^{-3}$ as $q=3$.
\end{enumerate}

However one should consider (\ref{3-86}) with ${\bar\ell}={\bar\ell}_2\Def(\mu^{1/2}h )^{2/3}r^{-2/3}$ and
${\bar\ell}={\bar\ell}_3\Def\mu^{-1/2}r^{-1/2}$; as $q=1$ one can take again $r=1$ obviously; as $q=2$ it is not so but one can take then by $q=9/5$ arriving to $O(\mu^{-1/2}h^{-3})$.

Then combining with (\ref{3-83}) I arrive to

\begin{proposition}\label{prop-3-14} The total contribution to the remainder of the zone $\cZ'_\out$ does not exceed $C\mu^{-1/2}h^{-3}+C\mu^{3/2}h^{-3/2}$ under condition $(\ref{0-8})_1$, $C\mu^{-1/2}h^{-3}+C\mu^{3/2} h^{-1-\delta}$ under condition $(\ref{0-8})_2$ and $C\mu^{-1/2}h^{-3}$ under condition $(\ref{0-8})_3$.
\end{proposition}

\begin{remark}\label{rem-3-15} The main part of this asymptotics is given by expression (\ref{0-12}) with $T={\bar T}_0$.
\end{remark}

\subsection{Estimates in the inner zone}
\label{sect-3-3}

Now I want to consider the \emph{inner zone\/}
\begin{align}
&\cZ_\inn\Def \bigl\{\dist(x,\Sigma)\le {\bar\gamma}_0=c\mu^{-1/2}\cdot \dist(x,\Lambda)^{1/2}\bigr\},
\label{3-90}
\intertext{in term consisting off the \emph{inner bulk zone}}
&\cZ_{\inn, II}\Def\cZ_\inn\cap \bigl\{ |Z_1| \cdot \dist (x,\Lambda)\le \epsilon \mu \cdot \dist(x,\Sigma)^2\bigr\}\label{3-91}\\
\intertext{which mimics the outer (bulk) zone $\cZ_\out$ and the \emph{true inner zone}}
&\cZ_{\inn, I}\Def\cZ_\inn\cap \bigl\{ |Z_1| \cdot \dist (x,\Lambda)\ge \epsilon \mu \cdot \dist(x,\Sigma)^2\bigr\}\label{3-92}
\end{align}
which mimics the inner zone for 2-dimensional Schr\"odinger operator. One needs to remember that trajectories can leave one of them and enter another one after a while.

Note first that the inner core is already covered by proposition \ref{prop-3-12}, so in the definition of the inner zone I can include condition $\dist(x,\Lambda)\ge \mu^{-1/4}$ and therefore use canonical form (\ref{3-28}) there. Actually using the same method as in the proof of proposition \ref{prop-3-12}, one can provide estimate $O(\mu h^{-3}r^2\gamma^2)=O(h^{-3}r^3)$ thus covering even $\cZ_\inn\cap \{\dist (x,\Lambda)\le \mu^{-1/6}\}$ but I don't need this.

\subsubsection{}
\label{sect-3-3-1}
I start from the inner bulk zone $\cZ_{\inn,II}$.
Let $|Z_1|\asymp \rho$, $\rho r \le \epsilon \mu\gamma^2$ and condition (\ref{3-77}) be fulfilled. Then exactly as in the outer zone the drift speed is $\asymp \mu^{-1}\rho^2 r\gamma^{-2}\le c\rho $ and evolution speed with respect to $hD_4$ does not exceed $c$; then for time $T_1=\epsilon\rho r$ dynamics remains in $\cZ_{\inn, II}$. Therefore in the domain where  condition (\ref{3-77})  is fulfilled and $\rho \ge  C(\mu h|\log h|)^{1/2}$ one can take $T_1=\epsilon\rho r$ and $T_0=Ch|\log h|$; however in the estimates this logarithmic factor is not needed.

Therefore the contribution to the remainder of this domain intersected with $\cZ_{(r,\gamma,\rho)}$ does not exceed $Ch^{-3}\gamma T_1^{-1}r^2\rho^2= C\mu^{-1/2}h^{-3}r^2\rho$ and the summation with respect to $(r,\gamma,\rho)$ results in $O(\mu^{-1/2}h^{-3})$.

On the other hand, consider the domain where condition (\ref{3-77}) is violated, but still $\rho \ge {\bar\rho}_1$. The contribution to the remainder of of this domain intersected with $\cZ_{(r,\gamma,\rho)}$ does not exceed $C\mu h^{-3}\gamma^2r^2\rho^2 $
where two factors $\gamma$ come as the width of the strip and the part of $T_1^{-1}\asymp \mu\gamma$. The summation with respect to $(r,\gamma,\rho)$: $r\rho^2\le c\gamma$ results in $C\mu^{-1/2}h^{-3}$.

Further, consider domain where condition (\ref{3-78}) is violated but still $\rho \ge {\bar\rho}_1$. In the similar manner one can estimate its contribution to the remainder by $C\mu h^{-3}\gamma_0^2 (\mu^{1/2}h|\log h|)^{2/3}$ which is $o(\mu^2h^{-2}+\mu^{-1/2}h^{-3})$ for sure and is $O(\mu^{-1/2}h^{-3})$ as
$\mu\le C(h|\log h|)^{-4/5}$.

Finally, contribution to the remainder of the domain where $\rho \le {\bar\rho}_1$ does not exceed $C\mu \gamma {\bar\rho}_1^2= C\mu^{3/2}h^{-2}|\log h|$.

Thus I arrive to

\begin{proposition}\label{prop-3-16} Contribution of $\cZ_{\inn,II}$ to the remainder
is $O(\mu^2h^{-2}+\mu^{-1/2}h^{-3})$ for sure. Furthermore, this contribution is is $O(\mu^{-1/2}h^{-3})$ as $\mu\le C(h|\log h|)^{-1/2}$. Again the main part of this asymptotics is given by expression $(\ref{0-12})$.
\end{proposition}

\subsubsection{}
\label{sect-3-3-2}
Again, the estimate achieved in proposition \ref{prop-3-16} is good in the general case but should be improved under non-degeneracy condition as
\begin{equation}
c(h|\log h|)^{-1/2}\le \mu \le C h^{-1}.
\label{3-93}
\end{equation}
Again one should consider regular final elements and singular ones and I consider only latter. Further introducing the scaling function $\ell_n$ by (\ref{3-62}) I can define $T_0$ and $T_1$ by (\ref{3-63}),(\ref{3-64}) as (\ref{3-65}) holds; furthermore estimate (\ref{3-66}) holds.

Further, if  condition  (\ref{3-78}) is fulfilled, I can increase $T_1$ to $T_2=\epsilon \ell$ since on this time interval both $|\nabla_{x_1,x_4,\xi_4}|$ change by no more than $\epsilon_1\min (\ell,r)$.

Then repeating the arguments of subsubsections \ref{sect-3-2-2} I find that the contribution of this part of $\cZ'_{\inn,II}$ to the remainder does not exceed expression similar to (\ref{3-80}) but with $\mu ^{-1}\gamma^{-2}$ replaced by 1:
\begin{multline}
C\mu h^{-3}\gamma r^{(2-q)_+} (\rho^2+r \ell^2) \times \Bigl(C_0(\mu h)^{-1}(\rho+r\ell)^2 +1\Bigr) \times \ell^q  \times h (\rho^2+r\ell^2)^{-1} \times 
 \ell^{-1}=\\
C h^{-3}\gamma  r^{(2-q)_+}  \Bigl((\rho+r\ell)^2 +\mu h\Bigr)\ell^{q-1};
\label{3-94}
\end{multline}
Then the sum with respect to $(\ell,\rho,r)$-partition with $\rho\ge \ell$ does not exceed $Ch^{-3}\gamma$ as $q\ge 2$ and $Ch^{-3}\gamma |\log h|$ as $q=1$
and one can get rid off logarithmic factor exactly as in the analysis in $\cZ'{\out,II}$.

Furthermore, condition  (\ref{3-77}) is fulfilled instead of (\ref{3-78}) and $\rho\ge R^{1/2}\ell$, I can increase $T_1$ to $T_2=\epsilon \rho$  
and the contribution of this part of $\cZ'_{\inn,II}$ to the remainder does not exceed expression similar to (\ref{3-82}) but with $\mu ^{-1}\gamma^{-2}$ replaced by 1:
\begin{multline}
C\mu h^{-3}\gamma r^{(2-q)_+} (\rho^2+r \ell^2) \times \Bigl(C_0(\mu h)^{-1}(\rho+r\ell)^2 +1\Bigr) \times \ell^q  \times h (\rho^2+r\ell^2)^{-1} \times 
 \rho^{-1}=\\
C h^{-3}\gamma  r^{(2-q)_+}  \Bigl(\rho^2 +\mu h\Bigr)\ell^q \rho^{-1};
\label{3-95}
\end{multline}
then summation with respect to $(\ell,\rho,r)$-partition results in
$Ch^{-3}\gamma + C\mu h^{-2}\gamma |\log h|$. Therefore the sum with respect to $\gamma$-partition with $\rho\ge \ell$ does not exceed $C\mu^{-1/2}h^{-3}$ as $q\ge 2$ and $\mu^{-1/2}h^{-3} + C\mu^{1/2}h^{-2}|\log h|$ as $q=1$.

On the other hand, due to the same arguments as in subsubsection \ref{sect-3-2-3} the contribution to the remainder of the domain where both conditions (\ref{3-77}),(\ref{3-78}) are violated but l (\ref{3-65}) still holds does not exceed (\ref{3-86}). Then the summation with respect to $(\ell,\rho,r,\gamma)$ partition results again in
$C\mu h^{-2} (\mu ^{1/2}h |\log h|)^{q/2} + O(\mu ^{-1/2}h^{-3})$.

Finally, contribution of the domain where (\ref{3-65}) fails is tackled in the same way as in subsubsection in subsubsection \ref{sect-3-2-3} and estimated by (\ref{3-87}).

So, I arrive to the following copy-cat of proposition \ref{prop-3-14} and remark \ref{rem-3-15}:

\begin{proposition}\label{prop-3-17} The total contribution to the remainder of the zone $\cZ_{\inn,II}$ does not exceed $C\mu^{-1/2}h^{-3}+C\mu^{3/2}h^{-3/2}$ under condition $(\ref{0-8})_1$, $C\mu^{-1/2}h^{-3}+C\mu^{3/2} h^{-1}|\log h|$ under condition $(\ref{0-8})_2$ and $C\mu^{-1/2}h^{-3}$ under condition $(\ref{0-8})_3$.

Again the main part of this asymptotics is given by expression $(\ref{0-12})$.
\end{proposition}

\subsubsection
\label{sect-3-3-3} Now let us consider the \emph{true inner zone\/} $\cZ_{\inn,I}$ defined by (\ref{3-94}). The crucial difference between this zone and $\cZ_\out \cap \cZ_{\inn, II}$ is that in the classical evolution $\dist(x,\Sigma)$ is not necessarily preserved\footnote{\label{foot-26}
As $|Z_1|\cdot \dist (x,\Lambda)\ge C \mu \cdot\dist(x,\Sigma)^2$ it is not preserved for sure; in the rest of this zone variation of $\dist(x,\Sigma)$ is of the magnitude of $\dist(x,\Sigma)$.}

I remind that this subzone (as the whole inner zone) should be studied as $r=\dist (x,\Lambda)\ge c\mu^{-1/4}$.
Also I remind that the contribution of subzone
$\{ |Z_1|\le C(\mu h\log h)^{1/2}\}$ to the remainder does not exceed
$C\mu h^{-3}\gamma \times \mu h |\log h|=C\mu^{3/2} h^{-2} |\log h|)$ which is
$O(\mu^2 h^{-2})$ for sure and $O(\mu^{-1/2}h^{-3})$ as
$\mu \le c (h|\log h|)^{-1/2}$.

Let us introduce $\gamma = \mu ^{-1/2} (r\rho)^{1/2}$; then in the classical evolution $\gamma$ is the magnitude of $\max \dist(x,\Sigma)$ along the cyclotron movement. Note that the contribution of subzone
$\{ (\mu h|\log h|)^{1/2} \le |Z_1| \le \rho\}$ to the remainder does not exceed
$C\mu h^{-3}\gamma^2 \rho^2= C h^{-3} \rho^3$ which sums to
$O(\mu^{-1/2}h^{-3})$ in subzone $\{ \rho \le c\mu^{-1/6}|\log h|^{-1/3}\}$.

Consider first the general case of $q=0$ when our goal is $O(\mu^{-1/2}h^{-3}+\mu^2h^{-2})$ estimate. Then the total contribution of $\cZ_{\inn,I}\cap\{|Z_1|\ge C(\mu h|\log h|)^{1/2}\}$ to the remainder is $O(h^{-3})$. Therefore in this case one needs to analyze only $\mu \le h^{-1/2}$ and only subzone $\{ \rho \ge c(\mu^{-1/6}+\mu^{2/3}h^{1/3}) \}$.

I am going to prove that the contribution of
$\cZ_{\inn, I}\cap \{\rho \ge C(\mu h|\log h|)^{1/2}\}$ to the remainder does not exceed $C\mu^{-1/2}h^{-3}+C(\mu^{1/2}h|\log h|)^{1/2} h^{-3}$ (see proposition \ref{prop-3-19}) but I prefer a bit less direct way to do this; I am concentrating on this zone due to the previous remarks and because I do not have non-degeneracy condition at the moment. Thus one can take $T_0=Ch|\log h|$ and $T_1=\epsilon \mu^{-1/2}\gamma^{-1}=\epsilon \mu^{-1}(\rho r)^{-1}$. Let us try to increase it.

According to the series of propositions of subsubsection \ref{sect-1-4-3} the drift speed in the inner zone is $\kappa r\rho (\xi_3 - k^*\rho) + O(\mu^{-1/2})$.
Then the shift during one cyclotron movement is
$\asymp r\rho (\xi_3 - k^*\rho) T_1$ provided
\begin{equation}
r\rho |\xi_3 - k^*\rho|\ge C \mu^{-1/2}
\label{3-96}
\end{equation}
and it is observable as
\begin{equation*}
r\rho |\xi_3 - k^*\rho| \times \mu ^{-1/2}\times |\xi_3 - k^*\rho| \ge Ch|\log h|
\end{equation*}
or equivalently
\begin{equation}
|\xi_3-k^*\rho| \ge
C(r\rho)^{-1/2} (\mu h\gamma |\log h|)^{1/2}= C(r\rho)^{-1/4} (\mu^{1/2} h |\log h|)^{1/2}.
\label{3-97}
\end{equation}
Both conditions (\ref{3-96}), (\ref{3-97}) are fulfilled as
\begin{equation}
|\xi_3-k^*\rho| \ge \Delta \Def
c(r\rho)^{-1/4} (\mu^{1/2} h|\log h|)^{1/2} + c\mu^{-1/2} (r\rho)^{-1/2}.
\label{3-98}
\end{equation}
Then as (\ref{3-98}) is fulfilled\footnote{\label{foot-27} This makes sense as $\Delta \le \rho$: i.e.
\begin{equation}
\rho \ge {\bar\rho}_2\Def r^{-1/5}(\mu^{1/2} h |\log h|)^{2/5}+\mu^{-1/3}r^{-1/3};
\label{3-99}
\end{equation}
otherwise one needs to take $\Delta=c\rho$, thus redefining
$\Delta\redef \min(\Delta,{\bar\rho}_2)$. } one can upgrade $T_1$ to $T_2=\epsilon r |\xi_3-k^* \rho| $ because then not only magnitude of $\rho$ but also of $|\xi_3-k^*\rho|$ will be preserved during the evolution. However taking the evolution in the right time direction it could be increased to
$T_2=\epsilon r |\xi_3-k^* \rho| ^{1-\delta'}$ exactly as in \cite{IRO6}. The standard calculations imply 

\begin{claim}\label{3-100} The contribution of subzone of $\cZ_{\inn,I}$ where condition (\ref{3-98}) is fulfilled intersected with $\{\rho \ge C (\mu h|\log h|)^{1/2}\}$ to the remainder does not exceed $C\mu^{1/2}h^{-3}$.
\end{claim}

So, one needs to study
\emph{near periodic zone\/}
\begin{equation}
\cZ_{\per}=\bigl\{ |\xi_3-k^*\rho| \le \Delta \bigr\},
\label{3-101}
\end{equation}
intersected with $\{ |Z_1|\asymp \rho \}\cap \{\dist (x,\Lambda)\asymp r\}$
as
\begin{equation*}
\rho \ge C(\mu h|\log h|)^{1/2}+{\bar\rho}_2.
\end{equation*}
Note however that as $r\le \epsilon$ then $|\nabla \omega|\asymp r^{-1}\omega$ and $|\Hess \omega|\asymp r^{-2}\omega$ and therefore in the rescaled coordinates $|\nabla _{x_3,x_4,\xi_4}W_n|\asymp \rho^2 $ as $r\le \epsilon$ and $\rho^2\ge Cr$ and $|\Hess _{x_3,x_4,\xi_4}W_n|\asymp \rho^2 $ as $r\le \epsilon$ and $\rho \ge Cr$. In the first case periodicity is broken for sure and in the second case periodicity is not broken only as $|\nabla W_n|\le C(\mu \gamma h )^{1/2}$ which in intersection of $\cZ_\per$ has measure not exceeding 
$C\mu^{1/2}h\gamma$ instead of $C(\mu^{1/2}h)^{1/2}\gamma$ and since one can take $T_1=\epsilon \mu^{-1} \gamma^{-1}$ here for sure the contribution of the reduced periodic zone to the remainder estimate does not exceed
$C  h^{-4} \times (\mu^{1/2}h|\log h|) \gamma\times \mu h \gamma \le C\mu ^{1/2}h^{-2}|\log h|$ thus implying

\begin{proposition}\label{prop-3-18} The contribution of
$\cZ_{\inn, I}\cap \{ |Z_1|\ge C(\mu h|\log h|)^{1/2}+ Cr \}$ to the remainder does not exceed $C\mu^{-1/2}h^{-3}$.

Again the main part of this asymptotics is given by expression $(\ref{0-12})$.
\end{proposition}

So, only $\rho \le Cr$ and thus $\gamma \le C\mu^{-1/2} r$ should be analyzed.

\begin{proposition}\label{prop-3-19} The contribution of
$\cZ_{\inn, I}\cap \{ |Z_1|\ge C(\mu h|\log h|)^{1/2}\}$ to the remainder does not exceed $C\mu^{-1/2}h^{-3}+C(\mu^{1/2}h|\log h|)^{1/2} h^{-3}$.

Again the main part of this asymptotics is given by expression $(\ref{0-12})$.
\end{proposition}

\begin{proof} Now I am concentrating only on even more reduced zone. In particular (\ref{3-99}) implies $r\ge (\mu ^{1/2}h |\log h|)^{1/3}+ \mu^{-1/4}$.

One can take $T_1= \epsilon \mu^{-1/2}(r\rho)^{-1/2}$ as condition (\ref{3-98}) fails. Then the contribution of $\cZ_\per$ to the remainder does not exceed
\begin{multline}
C\mu h^{-3}\gamma^2 r^2\rho \Delta = C h^{-3}r^3\rho^2\Delta=\\
C h^{-3}r^3\rho^2 \Bigl(c(r\rho)^{-1/4} (\mu^{1/2} h|\log h|)^{1/2} + c\mu^{-1/2} (r\rho)^{-1/2}\Bigr).
\label{3-102}
\end{multline}
Obviously when taking sum with respect to $(r,\rho)$ partition the latter expression sums to itself with
${\bar\Delta}=\Delta|_{r=1, \rho =1,\gamma=\mu^{-1/2}}$:
\begin{equation}
{\bar\Delta} =C (\mu^{1/2}h|\log h|)^{1/2} + C\mu^{-1/2}.
\label{3-103}
\end{equation}
Then (\ref{3-102}) becomes $O\bigl(\mu^{-1/2}h^{-3}+ h^{-3}(\mu^{1/2}h|\log h|)^{1/2}\bigr)$.

One can check easily that the contribution of subzone
$\bigl\{{\bar\rho}_2 \ge \rho \ge C(\mu h|\log h|)^{1/2}\bigr\}$ is just $O(\mu h^{-3}\gamma^2 {\bar\rho}_2^2)$ which is less than this.
\end{proof}

\subsubsection
\label{sect-3-3-4}

Again I would like to improve the remainder estimate achieved in proposition \ref{prop-3-19} under non-degeneracy condition \ref{0-8-q} with $q\ge 1$ also  getting rid off assumption $|Z_1|\ge C(\mu h|\log h|)^{1/2}$. In the analysis below instead of $W_n$ I consider its pullback $\omega (V-(2n+1)f_2\mu h)$ in the rescaled coordinates $x'$ on $\Sigma$; it is useful for Hessians but it transforms into properties of $W_n$ since I am interested only in the measures of sets $|\nabla W_n|\asymp \lambda$.

\begin{proposition}\label{prop-3-20} As $\mu \le ch^{-1}|\log h|^{-K}$ the contribution of $\cZ_{\inn, I}\cap \{ |Z_1|\ge  Cr \}$ to the remainder does not exceed $C\mu^{-1/2}h^{-3}$.
\end{proposition} 

\begin{proof} Let us introduce 
\begin{gather}
\ell =\ell_n=\epsilon r^{-1}|\nabla_{x_3,x_4,\xi_4} W_n|+{\bar\ell}, \label{3-104}\\
\intertext{and}
\varrho =\varrho_n =\epsilon |\xi_3-k^*W_n^{1/2}|+\Delta.
\label{3-105}
\end{gather}
Then the contribution of  $(r,\rho,\ell)$-elements to the asymptotics does not exceed  (\ref{3-69})-like expression
\begin{equation}
C\mu h^{-3}\gamma r^2 (\rho^2+r\ell^2) \times \Bigl(C_0(\mu h)^{-1}(\rho+r\ell^2) +1\Bigr); 
\label{3-106}
\end{equation}
if $\rho^2+r\ell^2\le C\mu h\gamma |\log h|$ then also $r^2\le C\mu h\gamma |\log h|$ and (\ref{3-106} does not exceed \newline $C\mu ^3 h^{-1}\gamma^3 |\log h|^2\le C\mu ^{3/2}h^{-1}|\log h|^2$. So let us consider $\rho^2+r\ell^2\ge C\mu h\gamma |\log h|$.

As $\rho^2 \ge Cr$ then periodicity is broken for sure (since $\nabla W_n|\ge \epsilon \rho^2-Cr$ and one can instantly upgrade $T_0= Ch|\log h| (\rho^2+r\ell^2)^{-1}$ to $T_2=\epsilon$. Otherwise it could be done as $\ell\ge C(\mu h |\log h|)^{1/2}$. Then contribution of such elements to the remainder does not exceed 
\begin{multline}
C\mu h^{-3}\gamma r^2 (\rho^2+r\ell^2) \times \Bigl(C_0(\mu h)^{-1}(\rho+r\ell^2) +1\Bigr)\times h|\log h|(\rho^2+r\ell^2) ^{-1}=\\
C\mu h^{-2}\gamma r^2  \times \Bigl(C_0(\mu h)^{-1}(\rho+r\ell^2) +1\Bigr) |\log h|
\label{3-107}
\end{multline}
because $|\Hess W_n|\ge \epsilon \rho^2- Cr^2$ and one can delete safely logarithmic factor because it is not in the estimate of the Fourier transform.

Otherwise one can take $T_1= \epsilon \mu^{-1}\gamma^{-1}$ and then one gets
\begin{multline}
C\mu h^{-3}\gamma r^2 (\rho^2+r\ell^2) \times \Bigl(C_0(\mu h)^{-1}(\rho+r\ell^2) +1\Bigr)\times h|\log h|(\rho^2+r\ell^2) ^{-1} \mu  \gamma \times \ell \Delta=\\
C\mu^2 h^{-2}\gamma^2 r^2  \times \Bigl(C_0(\mu h)^{-1}(\rho+r\ell^2) +1\Bigr) {\bar\ell}\Delta  |\log h|\le Ch^{-3}\times (\mu^{1/2}h)= \\ C\mu^{1/2}h^{-2}|\log h|^2
\label{3-108}
\end{multline}
where one copy of $(\mu^{1/2}h|\log h|)^{1/2}$ comes as the width of zone with respect to $\xi_3$ and another as the width with respect to $|\nabla W_n|$. 
\end{proof}

\begin{proposition}\label{prop-3-21}Under  condition $(\ref{0-8})_1$ be fulfilled. Then

\smallskip
\noindent
(i) The total contribution of $\cZ_{\inn,I} \cap \{|Z_1|\le \epsilon r\}$ to the remainder does not exceed $C\mu^{-1/2}h^{-3}$ as $\mu \le h^{-1}|\log h|^{-K}$.

\smallskip
\noindent
(ii) The total contribution of $\cZ_{\inn,I} \cap \{|Z_1|\ge \epsilon r\}$ to the remainder does not exceed
\begin{equation}
C\mu^{-1/2}h^{-3}+C \mu^{5/4}h^{-3/2} |\log h|^{1/2};
\label{3-109}
\end{equation}
in particular it is $O(\mu^{-1/2}h^{-3})$ as
$\mu \le h^{-6/7}|\log h|^{-2/7}$\,\,\footnote{\label{foot-28} Actually  I am going to prove later that one can drop logarithmic factors. In this case estimate cannot be improved without either additional assumptions or correction terms.}.

\smallskip
Again the main part of this asymptotics is given by expression $(\ref{0-12})$.
\end{proposition}

\begin{proof} 
 Note that the most delicate  is the case $\epsilon r\le \rho \le cr$ when even condition $(\ref{0-8})_2$ does not yield $(\ref{0-8})_1$ for $W_n$\,\footnote{\label{foot-29} However $(\ref{0-8})_1$ with the first and the second order derivatives taken only with respect to $\Lambda$ implies such condition for $W_n$ thus eliminating vicinity of $\Lambda$ from the contributing more than $C\mu^{-/2}h^{-3}$ to the remainder.}

\smallskip
\noindent
(a) Consider first $r\le \epsilon_1$, $\rho \asymp r$, $\gamma \asymp\mu^{-1/2}r^{1/2}\rho^{1/2}\asymp  \mu^{-1/2}r$. 

Let us introduce a scaling function
\begin{equation}
\ell_n =\epsilon r^{-1}|\nabla W_n| +{\frac 1 2}{\bar\ell},\qquad {\bar\ell}=C(\mu^{1/2} h|\log h|)^{1/2}
\label{3-110}
\end{equation}
and the corresponding subpartition as $\rho \le cr$. Note that actually $\ell \ge {\bar\ell}=C(\mu \gamma r^{-1}h|\log h|)^{1/2}$ is the condition that for time $T_1=\epsilon \mu^{-1}\gamma^{-1}$ the shift is observable. Then exactly as before

\begin{claim}\label{3-111} The total contribution to the remainder of all $\ell$-elements with $\ell \ge {\bar\ell}$, $\rho \asymp r$ does not exceed $C\mu^{-1/2}h^{-3}$.
\end{claim}

So only 
$\ell\le {\bar\ell}$ subelements should be considered.
Note that after rescaling $|\nabla \omega f_2|\asymp 1$ and therefore 

\begin{claim}\label{3-112}
For each ${\bar\ell}$ subelement there exists ${\bar n}$ such that 
$|\nabla \omega (V-(2n+1)f_2\mu h|\ge {\bar\ell}$ as $|n-{\bar n}|\ge
C{\bar\ell}(\mu h)^{-1}$ and therefore over each ${\bar\ell}$-element there are no more than $M= C{\bar\ell}(\mu h)^{-1}$ other subelements.
\end{claim}
Therefore the total contribution of them to the remainder does not exceed
\begin{multline}
C\mu h^{-3}\gamma  r^2 (\rho^2+r\ell^2) \times \Bigl(C_0(\mu h)^{-1}{\bar\ell} +1\Bigr) \times h(\rho^2+r\ell^2)^{-1} \times \mu \gamma \Delta\asymp
\\
C \mu h^{-2}r^3 \times \Bigl({\bar\ell}^2(\mu h)^{-1}  +{\bar\ell} \Bigr)
\label{3-113}
\end{multline}
since $\Delta \asymp {\bar\ell}$. After summation over $r$-partition one gets 
$C\mu ^{1/2} h^{-3}|\log h| + C\mu h ^{-2}{\bar\ell}$ which is exactly (\ref{3-109}).

\smallskip
\noindent
(b) Consider now $r\asymp 1$, $\rho \asymp 1$, $\gamma \asymp\mu^{-1/2}$.
One needs to recover an extra factor $\mu h$ in comparison with estimate of proposition \ref{prop-3-19} and only case $\mu \ge h^{\delta-2/3}$ needs to be considered.

Let us assume first that
\begin{equation}
|\Hess\bigl({\frac V {f_2}}\bigr)|\ge \epsilon_0
\label{3-114}
\end{equation}
and let us fix
\begin{equation}
\zeta \ge {\bar\zeta}=C
\bigl(\mu ^{1/2}h |\log h|\bigr)^{\frac 1 3} +\epsilon \mu h
\label{3-115}
\end{equation}
and examine subzone
\begin{equation}
\Omega_\zeta\Def \bigl\{ \exists n\in \bZ^+:\ |\Hess W_n| \asymp \zeta \bigr\}.
\label{3-116}
\end{equation}
Let us introduce $\zeta$-admissible partition in this zone; as $r\le \zeta$ one should replace it by $r$-admissible partition. Note that one can assume that

\begin{claim}\label{3-117}
For each element of such partition no more than $M\Def C_0\zeta (\mu h)^{-1}$ magnetic numbers $n$ satisfy condition $|\Hess W_n| \asymp \zeta$.
\end{claim}

Really, if $|\nabla \omega f_2|\asymp 1$ then arguments of (a) work. If 
$|\nabla \omega f_2|\asymp \le \epsilon_1$ but $|\nabla f_2^{-1}V|\asymp 1$
then $|\nabla \omega V|\asymp 1$ and then $|\nabla W_n|\asymp 1$ for all $n \le C(\mu h)^{-1}$ and everything is easy. So, let both $|\nabla \omega f_2|\asymp \le \epsilon_1$ and
$|\nabla f_2^{-1}V|\le \epsilon_1$. Then as $|\Hess \omega f_2|\asymp 1$ then 
(\ref{3-117}) is obviously true. Finally, if $|\Hess \omega f_2|\le \epsilon_1$
then (\ref{3-114}) implies that $|\Hess \omega V|\asymp 1$ and then
$|\Hess W_n|\asymp 1$ for all $n \le C(\mu h)^{-1}$.

Then for each  index $n$ described (\ref{3-117})  all the above arguments could be repeated but with 
\begin{equation}
\ell = \zeta^{-1} |\nabla W_n|+{\frac 1 2}{\bar\ell},\qquad {\bar\ell}\Def C(\zeta^{-1} \mu h |\log h|)^{1/2}
\label{3-118}
\end{equation}
leading to the contribution to the remainder
\begin{equation}
C\mu^2 h^{-3}\gamma^2 \times \zeta (\mu h)^{-1} \times h|\log h| \times {\bar\ell}\Delta \asymp C \mu^{1/2} h^{-2} |\log h|^2
\label{3-119}
\end{equation}
since $\zeta {\bar\ell}=\Delta$.

So I am left with $\zeta={\bar\zeta}$ and only with indices $n$ 
$|\Hess W_n| \le {\bar\zeta}$. Here one needs to consider only subelements with $\ell \le {\bar\ell}$ i.e. with 
\begin{equation}
|\nabla W_n \le C{\bar \zeta}^2.
\label{3-120}
\end{equation}
So, I am looking on ${\bar\zeta}$ partition elements satisfying this condition.

Assume first that $|\Hess \omega f_2|\asymp 1$. Let us introduce scaling function
\begin{equation}
\ell =  |\nabla \omega f_2|+{\frac 1 2}{\bar\ell};
\label{3-121}
\end{equation}
then as $\ell \ge {\bar\ell}$ over each such subelement leave no more than 
$C{\bar\zeta}^2/(\mu h \ell)$ indices $n$ satisfying (\ref{3-121}) and then (\ref{3-119}) is replaced by
\begin{align}
&C\mu^2 h^{-3}\gamma^2 \times {\bar\zeta}^2\ell ^{-1}(\mu h)^{-1} \times h|\log h| \times \Delta \asymp  C  h^{-3}\Delta  {\bar\zeta}^2\ell^{-1}|\log h|;
\label{3-122}\\
\intertext{alternatively one can replace (\ref{3-119}) by }
&C\mu^2 h^{-3}\gamma^2 \times {\bar\zeta} (\mu h)^{-1} \times h|\log h| \times \ell\Delta \asymp Ch^{-3}\Delta {\bar\zeta}\ell |\log h|
\label{3-123}
\end{align}
where $\ell$ is a width of of the strip where $|\nabla \omega f_2|\le \ell$.

Comparing (\ref{3-122}) and (\ref{3-123}) one can see that the best choice of ${\bar\ell}$ is ${\bar\ell}={\bar\zeta}^{1/2}$; then both of them become
$Ch^{-3}\Delta {\bar\zeta}^{3/2} |\log h|\asymp C\mu^{1/2}h^{-3}|\log h|^2$.

On the other hand, as $|\Hess \omega f_2|\le \epsilon_1$ then due to the arguments proving (\ref{3-117}) this analysis is not needed.

On the other hand, if (\ref{3-114}) is replaced by condition
\begin{equation}
|\nabla\bigl({\frac V{f_2}}\bigr)|\ge \epsilon_0
\label{3-124}
\end{equation}
the analysis is essentially the same, but simpler, with (\ref{3-116})--(\ref{3-117}) replaced by
\begin{align}
&\zeta \ge {\bar\zeta}=c (\mu^{1/2} h |\log h|)^{\frac 1 2} +\mu h
\label{3-125},\\
&\Omega_\zeta \Def \bigl\{ \exists n\in \bZ^+:\ |\nabla W_n| \asymp \zeta \bigr\},\label{3-126}
\end{align}
respectively and in (\ref{3-118}) $\Hess W_n|\asymp \zeta$ is replaced by $|\nabla W_n|\asymp \zeta$.

\smallskip
\noindent 
(c) Finally as $\rho \le \epsilon r$ condition $(\ref{0-8})_1$ translates for a similar condition to $W_n$ and then one can apply the arguments similar to those of proposition \ref{prop-3-20} which imply (i). I leave easy details to the reader.

\smallskip
\noindent 
(d) Analysis of $\rho^2 \le c(\mu h|\log h|)^{1/2}$ is also easy. Let us introduce $\ell$ by (\ref{3-110}). Then  as $\rho^2+r\ell^2 \le C \mu^{1/2}h|\log h|$ factor $(\rho^2+r\ell^2)$ takes care of everything; otherwise factor $\Delta$ is retained but factor $(\rho^2+r\ell^2)$ translates in the end into extra factor $\mu h|\log h$. I leave easy details to the reader.
\end{proof}

\subsubsection
\label{sect-3-3-5}
Now the better remainder estimate should be pursued only as
$\mu \ge h^{-6/7}|\log h|^{-2/7}$. Further, as $q=0,1$ the contribution of $\cZ_\per$ to the remainder is better than the remainder estimates $O(\mu^{-1/2}h^{-3}+\mu^2h^{-2})$,
$O(\mu^{-1/2}h^{-3}+\mu^{3/2}h^{-2-\delta})$ and therefore only $q=2,3$ should be considered.

Moreover, it follows from the proof of proposition \ref{prop-3-21} that if summation is taken over $\cZ_\per\cap \{\dist(x,\Lambda)\le r\}$ then the second term in (\ref{3-109}) would get an extra factor $r^4$ (from $r^2\gamma^2$) and  therefore only zone 
$\cZ_\per\cap \{\dist(x,\Lambda)\ge {\bar r}=\mu^{-1/16}|\log h|^{-1/8}\}$
 should be considered.
 
 Moreover, then ${\bar\zeta}\asymp \mu h$ (which is much larger than the first term in (\ref{3-114}) and therefore for each ${\bar\zeta}$-partition subelement only one magnetic number $n$ should be considered. 
 
 There is a minor problem as $\mu \ge h^{-1}|\log h|^{-K}$ because then  contributions of $\{\rho\le \epsilon r\}$ and $\{\rho\ge C r\}$
were estimated by $\mu^{1/2}h^{-2}|\log h|^K$; however in this case there also factor $r^2$ at least coming if looking at $\cZ_\per\cap \{\dist(x,\Lambda)\le r\}$ and therefore one needs to consider only 
$\cZ_\per\cap \{\dist(x,\Lambda)\ge |\log h|^{-K}\}$. Then in the first case condition $(\ref{0-8})_2$ supplies factor $(\mu^{1/2}h|\log h)^{1/2}|\log h|^K$ which is more than enough to take care of $|\log h|^K$; however I will just include both cases in the final analysis,

Then repeating all arguments of subsections 2.8, 2.9 of \cite{IRO6} with the obvious modifications, due to the presence of $r$ (which is greater than ${\bar r}\ge {\bar r}$) and of $(x_2,\mu^{-1}hD_2)$ in all the operators, one can prove easily the following

\begin{proposition}\label{prop-3-22} Under condition $(\ref{0-8})_2$ the total contribution of $\cZ_{\inn, I} \cap \{|Z_1|\ge \epsilon\}$ to the remainder does not exceed $C\mu^{-1/2}h^{-3}$
while the main part of the asymptotics is given by $(\ref{0-12})$ plus a correction term
\begin{align}
&\int \cE_\corr^\MW (x',\tau)\psi(x')\,dx',\label{3-127}\\
\intertext{where temporarily}
&\cE^\MW_\corr \Def h^{-1}\int \int_{-\infty}^0 \Bigl( F_{t\to h^{-1}\tau}
\bigl({\bar\chi}_T(t)-{\bar\chi}_{T_0}(t)\bigr)\bigl(\Gamma_x (uQ^t_y)\bigr)\Bigr)\,d\tau dx_1
\label{3-128}
\end{align}
and in this formula one can choose arbitrarily
\begin{align}
&T_0\in [Ch|\log h|, \epsilon '\mu^{-1/2}\gamma^{-1}],\label{3-129}\\
&T\ge C r^{-1/2} (h |\log h|)^{1/2} \bigl(\max (\bigl|\xi_3- k^*|Z_1|\bigr|,h^{1/2})\bigr)^{-1}\label{3-130}
\end{align}
and pseudo-differential $Q\equiv I$ in $\cZ_\per \cap \{\dist(x,\Lambda)\ge {\bar r}\}$ and supported in the same zone, slightly inflated.
\end{proposition}

In the next section this correction term will be rewritten in more explicit form.

\subsubsection{\sc Conclusion.\ }
\label{sect-3-3-6} So, with the remainder estimates described in theorems \ref{thm-0-1}, \ref{thm-0-3} asymptotics are derived where the main parts are given by rather implicit formula (\ref{4-1}).

The analysis in frames of theorem \ref{thm-0-5} is similar but much simpler. I leave it to the reader.

\section{Calculations}
\label{sect-4}
So, according to the previous section, with the remainder estimates described there the main part of the asymptotics of $\int e(x,x,0)\psi (x)\, dx $ is given by
\begin{equation}
\sum _\iota h^{-1}\int_{-\infty}^0 \Bigl( F_{t\to h^{-1}\tau} {\bar\chi}_{T_\iota}(t) \Gamma \bigl( u\psi Q_{\nu y}^t\bigr)\Bigr)\,d\tau
\label{4-1}
\end{equation}
where $Q_\iota$ form an appropriate partition of unity and $T_\iota \ge T_{(\iota)0}$
also described there.

In this section I will make all the calculations transforming implicit expression (\ref{4-1}) to a more explicit form, namely to
\begin{equation}
\int \cE^\MW (x,0)\psi (x)\, dx.
\label{4-2}
\end{equation}
There are many methods to make a reduction and I will employ them depending on the power of the magnetic field and non-degeneracy condition.

Also I will calculate (\ref{3-122}) in a more explicit form.

\subsection{Moderate magnetic field. I}
\label{sect-4-1}

I remind that I proved that for $\mu \le h^{-\delta}$ the remainder estimate $O(\mu^{-1/2}h^{-3})$ holds while the main part is given by (\ref{0-12}) with $T={\bar T}\Def Ch|\log h|$. Now I want to prove the same under a larger upper bound for $\mu$; further, increasing $\mu$ I want to recover a remainder estimate consistent with those in theorems \ref{thm-0-1}, \ref{thm-0-3}, \ref{thm-0-5} while keeping the main part given by (\ref{0-12}) with $T={\bar T}$. To do this I need to estimate properly (\ref{1-22})-like correction term
\begin{equation}
h^{-1}\int_{-\infty}^0 \Bigl( F_{t\to h^{-1}\tau} \bigr({\bar\chi}_{T'}(t) - {\bar\chi}_{T''}(t)\bigr)\bigl(\Gamma uQ_y^t\bigr)\Bigr)\,d\tau
\label{4-3}
\end{equation}
with $T'=T_0$ derived in section \ref{sect-3} and
$T''={\bar T}$\,\footnote{\label{foot-30} However in what follows I will consider some other parameters $T'>T''$.}. This will let me to rewrite (\ref{0-12}) as a Weyl expression, may be with $O(\mu^2h^{-2})$ correction. Here $Q=\psi$ but in the analysis below I will need more general PDOs.

This expression (\ref{4-3}) is the sum of expressions of the type
\begin{equation}
h^{-1}\int_{-\infty}^0 \Bigl( F_{t\to h^{-1}\tau} \bigr({\bar\chi}_{2T}(t) - {\bar\chi}_T(t)\bigr)\Gamma uQ_y^t\Bigr)\,d\tau=
i T^{-1}\Bigl( F_{t\to h^{-1}\tau}\chi_T(t) \Gamma uQ_y^t\Bigr)\Bigr|_{\tau=0}
\label{4-4}
\end{equation}
with $\chi (t)= t^{-1}\bigl({\bar\chi}({\frac 1 2}t) - {\bar\chi}(t)\bigr)$
and $T=2^kT''$, $k$ running from 2 to $\lceil\log_2 (T'/T'')\rceil$.

\begin{proposition}\label{prop-4-1} Let assumptions of theorem \ref{thm-0-1} with
\begin{equation}
h^{-\delta}\le \mu \le \epsilon (h|\log h|)^{-1}
\label{4-5}
\end{equation}
be fulfilled and let $Q$ be supported in $\{\dist(x,\Sigma)\le \gamma\}$,
$\gamma \ge \mu^{-1/2}$.

Then expression $(\ref{4-3})$ with $T'=\epsilon \mu^{-1}\gamma^{-1}$ and $T''={\bar T}$ does not exceed $C\mu^{-1/2}h^{-3}+C\mu^2 h^{-2}$\, \footnote{\label{foot-31} Term $C\mu^{-1/2}h^{-3}$ is probably of no need in some estimates here and below but it is in the remainder estimate anyway.}.
\end{proposition}

\begin{proof} Note first that expression (\ref{4-4}) with
$T\in [{\bar T}, \epsilon \mu^{-1}]$ is negligible due to condition (\ref{0-2}).
So, until the end of the proof I can consider $T''=\epsilon \mu^{-1}$.

Further note that, as $Q$ is supported in $\{|Z_1|^2\ge C\mu h |\log h|\}$ expression (\ref{4-4}) with $T\in [T'',T']$ is negligible as well. Then for a general PDO $Q$ expression (\ref{4-4}) does not exceed $C\mu^2 h^{-2}|\log h|$.

Then proposition is proven as $\mu \le c(h|\log h|)^{-2/5}$ and until the end of the proof I can assume that
\begin{equation}
(h|\log h|)^{2/5}\le \mu \le \epsilon (h|\log h|)^{-1}.
\label{4-6}
\end{equation}
Further, as $Q$ is supported in $\{\dist(x,\Sigma)\le {\bar\gamma}_3= c|\log h|^{-1}\}$, expression (\ref{4-3}) would not exceed
$C\mu^2h^{-2}{\bar\gamma}_3|\log h|\asymp \mu^2 h^{-2}$
and as $Q$ is supported in $\{|Z_1|^2\le c\mu h\}$ expression (\ref{4-3}) would not exceed $C\mu^2 h^{-2}$ as well.

So, the only zone in question is
\begin{equation}
\bigl\{ C_0 \mu h \le\rho^2 \le C\mu h |\log h| \bigr\}\cap
\bigl\{\dist(x,\Sigma)\asymp\gamma, |Z_1|\asymp\rho \bigr\}
\label{4-7}
\end{equation}
with $|\log h|^{-1}\le \gamma \le \epsilon$. In this zone one needs to consider expression (\ref{4-4}).

Using precanonical form and decomposition (\ref{2-16})
one can rewrite (\ref{4-4}) as the sum of the similar expressions but with $u, Q, \Gamma$ replaced by $u_n$, $Q_n$ and $\Gamma'$ respectively:
\begin{equation}
i T^{-1}\Bigl( F_{t\to h^{-1}\tau}\chi_T(t) \Gamma' u_nQ_{ny}^t\Bigr)\Bigr|_{\tau=0}.
\label{4-8}
\end{equation}

From the decomposed problem
\begin{align}
&(hD_t-\cA_n)u_{nn'}\equiv 0,\label{4-9}\\
&u_{nn'}|_{t=0}\equiv \delta_{nn'}\delta (x'-y')\label{4-10}
\end{align}
one can prove not only that expression (\ref{4-8}) is negligible as
$T\ge T_\rho \Def Ch\rho^{-2}|\log h|$ but also that its absolute value does not exceed $C\mu h^{-2}\gamma \rho^2 \bigl(h\rho^{-2} T\bigr)^{-s}$ as
$t\le T_\rho$. The summation of this expression with respect to $T$ from $\epsilon \mu^{-1}$ to $\infty$ results in
$C\mu h^{-2}\gamma \rho^2 \bigl(\mu h\rho^{-2} \bigr)^{-s}$.

Further, summation of the latter expression with respect to $n$ results in
\begin{equation*}
C\mu h^{-2}\gamma \rho^2 \bigl(\mu h \rho^{-2} \bigr)^{-s}\times
\rho^2(\mu h)^{-1} = Ch^{-3}\rho ^4 \bigl(\mu h \rho^{-2} \bigr)^{-s}
\end{equation*}
since no more than $\rho^2 (\mu h)^{-1}$ numbers $n$ violate the ellipticity of $\cA_n$.

Finally, the summation with respect to $\rho\ge c(\mu h)^{1/2}$ results in $C\mu^2 h^{-2}\gamma$. So, one can conclude that expression (\ref{4-3}) with indicated $T'$ and $T''$ by absolute value does not exceed $C\mu^2h^{-2}$.
\end{proof}

One can improve the above error estimate under nondegeneracy condition.

\begin{proposition}\label{prop-4-2} Let assumptions of theorem \ref{thm-0-1} and condition $(\ref{4-6})$ be fulfilled. Further, let $Q$ be supported in
$\{\dist(x,\Sigma)\le \gamma\}$ with $\gamma \ge \mu^{\delta-1/2}$. Finally, let condition \ref{0-9-q} with $q=1,2$ be fulfilled.

Then expression $(\ref{4-3})$ with $T'=\epsilon \mu^{-1}\gamma^{-1}$, $T''={\bar T}$ and $Q$ supported in $\{\dist(x,\Sigma)\le \gamma\}$ does not exceed
\begin{phantomequation}\label{4-11}\end{phantomequation}
\begin{equation}
C\mu^{-1/2}h^{-3} + C (\mu h)^{2+q/2}|\log h|^K h^{-4}\gamma.
\tag*{$(4.11)_q$}
\label{4-11-q}
\end{equation}
\end{proposition}

\begin{proof}
There is no need to consider the inner core. Further, due to proposition \ref{prop-4-1} the case $\mu \le Ch^{-2/5}$ is already covered. So, let us assume that $\mu\ge h^{-2/5}$.

Consider precanonical form (\ref{2-17}) and define $\ell = |\nabla_V^\#$, $\rho =|Z_1^\#|$. Then even as $ \rho\le {\bar\rho}\Def C(\mu h |\log h|)^{1/2}$ but
$\ell \ge {\bar\ell}\Def {\bar\rho}$ one can still
trade $T'=\epsilon \mu^{-1}\gamma^{-1}$ to $T''=\epsilon \mu^{-1}$ because then the shift for time $T\in [\epsilon \mu^{-1},T'']$ is still observable.
Here perturbation terms and their derivatives are at most $C\mu^{-2}\ll {\bar \ell}$.

On the other hand, the total contribution to the remainder of all partition elements with $\rho \le {\bar\rho}$, $\ell \le {\bar\ell}$, $\dist(x,\Sigma)\le \gamma$ does not exceed
$C{\bar\rho}^2 {\bar\ell}^{2+q} \gamma h^{-4}$ which is exactly \ref{4-11-q}.
\end{proof}

\begin{remark}\label{rem-4-3}
While the error estimate of proposition \ref{prop-4-1} is as good as I need in the general case, the error estimate of proposition \ref{prop-4-2} is as good as I need under condition \ref{0-9-q} only under some restrictions to $\mu$ and/or $\gamma$.
\end{remark}

\subsection{Moderate magnetic field. II}
\label{sect-4-2}

Now I want to trade $T''= \epsilon \mu ^{-1}\gamma^{-1}$ to a larger value $T_0$
which provides the remainder estimate derived in section \ref{sect-3}. Now I will need to use some canonical forms and thus different zones will be treated separately.

\subsubsection{}
I will start from the strictly outer zone.
\begin{proposition}\label{prop-4-4} Let assumptions of theorem \ref{thm-0-1} and
condition $(\ref{4-5})$ be fulfilled and let $Q$ be supported in
$\cZ^*_\out\cap \{\dist(x,\Sigma)\le \gamma\}$ with $\gamma \ge C_0\mu^{\delta-1/2}$.
Then

\smallskip
\noindent
(i) In the general case expression $(\ref{4-3})$ with $T'=T_0$ and
$T''={\bar T}$ does not exceed
\begin{phantomequation}\label{4-12}\end{phantomequation}
\begin{equation}
C\mu^{-1/2}h^{-3}+C(\mu h)^{(4+q)/2}h^{-4}\gamma |\log h|^K+
C(\mu h)^{(3+q)/2}h^{-4}\gamma^2|\log h|^K
\tag*{$(4.12)_q$}
\label{4-12-q}
\end{equation}
with $q=0$;

\smallskip
\noindent
(ii) Under condition \ref{0-9-q} this expression does not exceed \ref{4-12-q}.
\end{proposition}

\begin{proof} One needs to consider only $Q$ supported in a zone in which in the
proofs of propositions \ref{prop-4-1}, \ref{prop-4-2} term (\ref{4-3}) with $T'=\epsilon \mu^{-1}$ and $T''={\bar T}$ was negligible. Then I am going to prove \ref{4-12-q} without the second term.

(i) Assume first that $Q$ is supported in $\cZ_{(r,\gamma)}\subset \cZ^*_{\out,I}$.

Let us consider canonical form (\ref{2-54}) and introduce function $\ell^*$ by (\ref{3-38}). If $\ell \ge 1$ (i.e. $p\ge Cr^{-1}(\mu h)^{-1}$) I do not partition the ``final'' ball $B(0,1)$ corresponding to the original $(\gamma; \gamma r^{-1})$ element; otherwise I make $\ell $-subpartition. Then on any given $\ell $ subpartition element the propagation speed with respect to $x_3$ does not exceed $ \mu^{-1}\gamma^{-2}r\ell_p$; further the propagation speed is of this magnitude as $|\partial_{\xi_3} H_{pn}| \asymp\ell$; in this case the shift for time $T$ is observable as it satisfies the logarithmic uncertainty principle
\begin{equation*}
\mu^{-1}\gamma^{-2} \ell \times T \times \ell \ge
Cr^2\mu^{-1}h\gamma^{-3}|\log h|.
\end{equation*}
Plugging $T\asymp \mu^{-1}\gamma^{-1}$ one arrives to
\begin{equation}
\ell \ge {\bar\ell}_1 \Def C\bigl( \mu h|\log h|\bigr)^{1/2}.
\label{4-13}
\end{equation}
Similarly, considering propagations with respect to other variables one arrives to the same conclusion as $|\nabla W_{np}|\asymp \ell$.

Note that for given index $p$ and for index $n$ violating ellipticity
$\ell \asymp \ell_p= \bigl(\ell^*+ r |p-{\bar p}|\mu h \bigr)$ with some ${\bar p}\le C_0(\mu h)^{-1}$.

Further, on each $\ell^*$ subpartition element condition (\ref{4-13}) is fulfilled for all indices $p$ as $\ell^*\ge {\bar\ell}_1$ and condition (\ref{4-13}) is violated for no more than $C_0{\bar\ell}_1 (r\mu h)^{-1}$ ($\gg 1$) indices $p$ as $\ell^*\le {\bar\ell}_1$. Furthermore for each such index $p$ ellipticity is violated by no more than $C_0\bigl({\bar\ell}_1^2\gamma(r\mu h)^{-1}+1\bigr)$ indices $n$.

Therefore the total contribution to (\ref{4-3}) of all singular subelements (i.e. subelements with $\ell^*\le C_0{\bar\ell}_1$) belonging to $\cZ_{(r,\gamma)}$ does not exceed
\begin{equation}
C\mu^2 h^{-2} \gamma^2 r^2\times \Bigl({\bar\ell}_1 (r\mu h)^{-1}+1\Bigr) \times \Bigl({\bar\ell}_1^2\gamma (r\mu h)^{-1}+1\Bigr) \times \ell_1^q ;
\label{4-14}
\end{equation}
since the second and the third factors do not exceed $C{\bar\ell}_1(r\mu h)^{-1}$ and $C|\log h|^K$ respectively, I get 
$C h^{-4} \bigl( \mu h\bigr)^{(3+q)/2} \gamma^2|\log h|^K$. This expression sums with respect to $r, \gamma$ to the same expression calculated at the maximal values of $\gamma$ and $r=1$, and multiplied by $|\log h|$, which is exactly expression \ref{4-12-q}. I increase $K$ as needed.

\smallskip
\noindent
(ii) Zone $\cZ^*_{\out,II}$ is treated in the same way.
\end{proof}

While estimates \ref{4-11-q} with $q\ge 0$ and even \ref{4-12-q} with $q\ge 1$ are sufficient for my needs, $(\ref{4-12})_0$ is not and needs to be improved to \ref{4-12-q} with $q$ arbitrarily close to 1:

\begin{proposition}\label{prop-4-5} Let assumptions of theorem \ref{thm-0-1} and
condition $(\ref{4-5})$ be fulfilled and let $Q$ be supported in
$\cZ^*_\out\cap \{\dist(x,\Sigma)\le \gamma\}$ with $\gamma \ge C_0\mu^{\delta-1/2}$.
Then expression $(\ref{4-3})$ with $T'=T_0$ and $T''={\bar T}$ does not exceed
\begin{equation}
C\mu^{-1/2}h^{-3}+C\mu^2 h^{-2}\gamma |\log h| + C(\mu h)^{(3+q)/2}h^{-4}\gamma^2|\log h|^K
\label{4-15}
\end{equation}
with $q<1$ arbitrarily close to $1$.
\end{proposition}

\begin{proof} Again, one needs to consider only $Q$ supported in a zone $\{|Z_1|\ge C (\mu h |\log h|)^{1/2}\}$ in which in the proof of proposition \ref{prop-4-1} term (\ref{4-3}) with $T'=\epsilon \mu^{-1}$ and $T''={\bar T}$ was negligible. Then I am going to prove (\ref{4-15}) without the second term.

(i) Assume first that $Q$ is supported in $\cZ_{(r,\gamma)}\subset \cZ^*_{\out,I}$. One needs to consider only $p\le C_0(r\mu h)^{-1}$.

Let us introduce the scaling function and $\ell$-admissible partition the in the ``final'' ball $B(0,1)$ corresponding to the original $(\gamma; \gamma r^{-1})$ element:
\begin{equation}
\ell^*_k =\ell_k \Def \min_{n,p}\epsilon \Bigl(\zeta^{-1}\sum _{j\le k} |\nabla ^j \cA_{np}|^{(m+1)/(k+1-j)}\Bigr)^{1/(m+1)}+{\bar\ell}_k
\label{4-16}
\end{equation}
with $k=m$, $\zeta=\gamma r^{-1}$ and
\begin{equation}
{\bar\ell}_k = \bigl( \mu h |\log h|^K\bigr)^{1/(k+1)}.
\label{4-17}
\end{equation}
As $m=1$ this function coincides with $\ell^*$ given by (\ref{3-38}).

Consider first partition \emph{$m$-singular\/} groups {\sf($(p,n)$,element)} i.e. groups such that
\begin{equation}
|\partial^\alpha \cA_{np}|\le c\zeta \ell^{k-|\alpha|}\qquad \forall \alpha:|\alpha|\le k
\label{4-18}
\end{equation}
with $k=m$, $\ell ={\bar\ell}_m$ and $\zeta=\zeta_m$.

Then contribution of all $m$-singular groups in $\cZ_{(r,\gamma)}$ could be estimated in the same manner as (\ref{4-14}):
\begin{equation}
C\mu^2 h^{-2} r^2 \gamma^2\times \Bigl(\zeta \ell ^k (\gamma\mu h)^{-1} +1\Bigr) \times \Bigl(\zeta \ell ^{k+1} (\mu h)^{-1}+1\Bigr)
\label{4-19}
\end{equation}
with the same $k$, $\ell$ and $\zeta$ as in (\ref{4-18}). Here the second factor estimate the number of indices $p$ involved and the third factor indicates the number of indices $n$ violating ellipticity of $\cA_{pn}$ for given $p$.

The groups which are not $k$-singular are \emph{$k$-regular\/}, i.e. they satisfy (\ref{4-18}) with some $\ell\ge 2{\bar\ell}_k$ and also
\begin{equation}
|\partial^\alpha \cA_{np}|\asymp \zeta_k\ell ^{k-|\alpha|}\qquad \text{for some }\alpha:|\alpha|\le k.
\label{4-20}
\end{equation}
Then one can rescale such element to $B(0,1)$. After this rescaling conditions (\ref{4-18}),(\ref{4-20}) are fulfilled with $\ell\redef 1$ and $\zeta\redef \zeta\ell^k$. So, one can apply the same rescaling and partition to it with $(k-1)$ instead of $k$ and with crucial $\ell$ (separating $(k-1)$-singular from $(k-1)$-regular) equal to ${\bar\ell}_{k-1}/\ell_k$. This means that if one returns to the original element, the radius of the $k$-singular element would be ${\bar \ell}_k$ every time.

This continuation however has a new property: every time when regular or singular $\ell'$-subelements appear inside of some regular elements, their relative densities do not exceed $\ell'/\ell$ (all the time I refer to the same scale) and therefore the absolute density of all $\ell'$-subelements does not exceed $C\ell' {\bar\ell}_m^{-1}$. 

Therefore the total contribution of all $k$-singular elements does not exceed (\ref{4-19}) multiplied by ${\bar\ell}_k{\bar\ell}_m^{-1}$
\begin{equation}
C{\bar\ell}_m^{-1}\mu^2 h^{-2} r^2 \gamma^2\times \Bigl(\zeta \ell ^{k+1} (\gamma\mu h)^{-1} +\ell \Bigr) \times \Bigl(\zeta \ell ^{k+1} (\mu h)^{-1}+1\Bigr)
\label{4-21}
\end{equation}
with $\ell= {\bar\ell}_k$ and $\zeta=\gamma r^{-1}$
since I know that $\zeta$ does not increase (even if it is different in the different elements). Then two last factors do not exceed $|\log h|^K$ for sure and the total contribution of all singular {\sf($(p,n)$,subelement)} groups does not exceed ${\bar\ell}_m^{-1} \mu^2 h^{-2}\gamma^2r^2 |\log h|^K$ and summation with respect to $(r,\gamma)$ results in the same expression with the maximal values of $\gamma$ and $r=1$.

Microlocally this is sound as logarithmic uncertainty principle
\begin{equation}
{\bar\ell}_k^2 \ge \mu^{-1}h r^2\gamma^{-3}|\log h|^K
\label{4-22}
\end{equation}
and to satisfy it for any $\gamma\gg \mu^{-1/2}r^{1/2}$ one needs to take $k=1$.

This leaves me with 1-regular groups and there are three types of them:
\begin{enumerate}[label=(\alph*)]
\item With $\zeta \ell^2 \le \gamma r^{-1} \mu h |\log h|^K$, which are covered by the same estimate (\ref{4-21});

\item With $\zeta \ell^2 \ge \gamma r^{-1} \mu h |\log h|^K$;
for them the shift $\mu^{-1}\gamma^{-2}r\zeta \ell T$ with $T=\mu^{-1}\gamma^{-1}$ is observable since
$\mu^{-1}r ^2\gamma^{-3} \zeta \ell \times \mu^{-1}\gamma^{-1} \times \ell \ge Cr^2\mu ^{-1}h\gamma^{-3}$;

\item $0$-regular groups but then $\cA_{pn}$ is just elliptic on the corresponding elements.
\end{enumerate}

\smallskip
\noindent
(ii) Zone $\cZ^*_{\out,II}$ is treated in the same way.
\end{proof}

\subsubsection{}\label{sect-4-2-2}
In the general case my purpose is the remainder estimate $O(\mu^{-1/2}h^{-3}+\mu^2h^{-2})$.

\begin{proposition}\label{prop-4-6} Let assumptions of theorem \ref{thm-0-1} be fulfilled. Then

\smallskip
\noindent
(i) As
\begin{equation}
h^{-\delta}\le \mu \le {\bar\mu}_0\Def=h^{\delta-2/5}
\label{4-23}
\end{equation}
$\int e(x,x,\tau)\psi (x)\,dx$ is given modulo $O(\mu^{-1/2}h^{-3})$ by expression $(\ref{0-12})$ with $T={\bar T}$;

\smallskip
\noindent
(ii) Moreover, as
\begin{equation}
h^{\delta-2/5}\le \mu \le h^{\delta-1}
\label{4-24}
\end{equation}
and $\psi(x)$ is a nice\footnote{\label{foot-32} F.e. $\psi(x)=\psi_1 (x_1/{\bar\gamma}_2)\psi_2(x_2,x_3,x_4)$.} function supported in $\{|x_1|\le {\bar\gamma}_2=\mu^{-\delta'}\}$
$\int e(x,x,\tau)\psi (x)\,dx$ is given modulo $O(\mu^{-1/2}h^{-3}+\mu^2h^{-3})$ by expression $(\ref{0-12})$ with $T={\bar T}$
\end{proposition}

\begin{proof}
Combining propositions \ref{prop-3-2}, \ref{prop-3-4} and \ref{prop-4-5} I conclude that this estimate holds in frames of (i) for $\psi$ supported in
$\cZ^*_\out$ and in frames of (ii) for $\psi$ supported in
$\cZ^*_\out\cap\{|x_1|\le \mu^{-\delta}\}$.

In propositions \ref{prop-3-11}, \ref{prop-3-14}, \ref{prop-3-16} and \ref{prop-3-19} it was essentially proven that this estimate holds for
$\psi$ supported in the near outer zone $\cZ'_\out$, inner core, inner bulk zone $\cZ_{\inn,II}$ and true inner zone $\cZ_{\inn,I}$ respectively.
\end{proof}

\subsection{Intermediate magnetic field. I}
\label{sect-4-3}

\subsubsection{}\label{sect-4-3-1}
Now I want just to calculate (\ref{0-12}) with $T=Ch|\log h|$. It is well-known (see for example \cite{IRO5}) that

\begin{proposition}\label{prop-4-7} Let $\mu \le h^{-1+\delta}$ with arbitrarily small $\delta>0$. Then under condition $(\ref{0-2})$

\smallskip
\noindent
(i) Asymptotics holds
\begin{equation}
h^{-1}\int_{-\infty}^\lambda \Bigl( F_{t\to h^{-1}\tau} {\bar\chi}_T(t) \Gamma _x uQ_y^t\Bigr)\,d\tau \sim
\sum_{(n, m)\in \bZ^{+\,2}} \kappa_{n,m, Q} (x,\lambda) h^{-4+2m +2n} \mu^{2n}
\label{4-25}
\end{equation}
as $|\lambda|\le \epsilon$, ${\bar T}\le T \le \epsilon \mu^{-1}$.

\smallskip
\noindent
(ii) Moreover, as $Q=Q(x)$ \ \ $\kappa_{n,0, Q} (x,.)=\kappa_{n,0, I}(x,.)Q(x)$ and
\begin{equation}
Th^{-1} \int_{-\infty}^\infty {\widehat{\bar\chi}}\bigl((\lambda -\tau) Th^{-1}\bigr) \cE^\MW (x,\tau)\,d\tau \sim
\sum_{n\in \bZ^{+}} \kappa_{n,0, I} (x,\lambda) h^{-4+2m +2n} \mu^{2n}
\label{4-26}
\end{equation}
with the Standard Weyl Expression
\begin{equation}
\kappa_{0,0,I} (x,\lambda) h^{-4}=\cE^\W (x,\tau)\Def {\frac 1 {32\pi^2}} (2\tau + V)^2 \sqrt g .
\label{4-27}
\end{equation}
\end{proposition}
Combining propositions \ref{prop-4-6}(i) and \ref{prop-4-7}(i) I arrive to

\begin{corollary}\label{cor-4-8} Under conditions $(\ref{0-2})$ and $(\ref{4-23})$  modified asymptotics $(\ref{0-3})$ (with $\cE^\MW$ replaced by $\cE^\W$) holds.
\end{corollary}

To prove theorem \ref{thm-0-1} under condition (\ref{4-23}) it is sufficient to estimate properly
\begin{equation}
|\int \Bigl( \cE^\MW (x,\lambda)-\cE^\W (x,\lambda)\Bigr)\psi(x)\,d\tau|
\label{4-28}
\end{equation}
which due to proposition \ref{4-7}(ii) is equivalent to
\begin{equation}
|\int \Bigl( \cE^\MW (x,\lambda)- Th^{-1} \int_{-\infty}^\infty {\widehat{\bar\chi}}\bigl((\lambda -\tau) Th^{-1}\bigr) \cE^\MW (x,\tau)\,d\tau \Bigr)\psi(x)\,d\tau|.
\label{4-29}
\end{equation}

\begin{proposition}\label{prop-4-9} (i) Under conditions $(\ref{0-2})$ and $(\ref{4-23})$ both expressions $(\ref{4-28})$, $(\ref{4-29})$ with $\lambda=0$ do not exceed $C\mu^{-1/2}h^{-3}$.

\smallskip
\noindent
(ii) Under conditions $(\ref{0-2})$ and $(\ref{4-24})$ both expressions $(\ref{4-28})$, $(\ref{4-29})$ with $\lambda=0$ and a nice$^{\ref{foot-32}}$ function $\psi$ supported in $\{|x_1|\le {\bar\gamma}_2=\mu^{-\delta'}\}$ do not exceed $C\mu^{-1/2}h^{-3}+C\mu^2h^{-2}$.
\end{proposition}

\begin{proof} Proof is standard, based on the same scaling functions and partition as in the proof of proposition \ref{prop-4-5} as $\psi$ is supported in $\cZ^*_\out$. These arguments work in the other zones as well since I do not need uncertainty principle anymore. I leave the easy details to the reader.
\end{proof}

Now propositions \ref{prop-4-6}(i), \ref{prop-4-7} and \ref{prop-4-9}(i) imply 

\begin{corollary}\label{cor-4-10}
Under conditions $(\ref{0-2})$ and $(\ref{4-23})$   asymptotics $(\ref{0-3})$  holds.
\end{corollary}

\subsubsection{}
\label{sect-4-3-2}
In view of theorem \ref{thm-0-1} nondegeneracy condition should be used only as 
\begin{equation}
c  h^{\delta-2/5}\le \mu \le Ch^{-1}
\label{4-30}
\end{equation}
with $\delta=0$. Then, under condition $(\ref{0-9})_1$ the remainder estimate derived in section \ref{sect-3} is $O(\mu^{-1/2}h^{-3}+\mu^{3/2}h^{-3/2-\delta})$ i.e. it is
$O(\mu^{-1/2}h^{-3})$ as $\mu \le h^{\delta-3/4}$ and then there is no need to invoke \ref{0-9-q} with $q\ge 2$. However it makes life much easier and I will do it right now leaving more difficult analysis for the later; so $\delta>0$ is arbitrarily small in this subsubsection.

First of all one needs to prove

\begin{proposition}\label{prop-4-11} Let conditions $(\ref{0-2})$ and \ref{0-9-q} be fulfilled and $\mu\le h^{\delta-1}$. Let $\psi$ be a nice function supported in $\{|x_1|\le \gamma\}$ with $\mu^{1/2-\delta}\le \gamma \le\epsilon$. Then
$\int e(x,x,\tau)\psi (x)\,dx$ is given by expression $(\ref{0-12})$ with $T={\bar T}$ with an error not exceeding \ref{4-12-q}.
\end{proposition}

\begin{proof}
Combining propositions \ref{prop-3-7} and \ref{prop-4-4} I conclude that the statement holds for $\psi$ supported in
$\cZ^*_\out\cap \{|x_1|\le\gamma\}$.

In propositions \ref{prop-3-14}, \ref{prop-3-12}, \ref{prop-3-17} and \ref{prop-3-22} it was essentially proven that this estimate holds for
$\psi$ supported in the near outer zone $\cZ'_\out$, inner core, inner bulk zone $\cZ_{\inn,II}$ and true inner zone $\cZ_{\inn,I}$ respectively.
\end{proof}

Proposition \ref{prop-4-9} and \ref{prop-4-7} imply

\begin{corollary}\label{cor-4-12} Let conditions $(\ref{0-2})$ and \ref{0-9-q} be fulfilled and $\mu\le h^{\delta-1}$. Then

\smallskip
\noindent
(i) One can define the main part of asymptotics by $(\ref{0-12})$ with $T=Ch|\log h|$ (or equivalently by $(\ref{4-26}$) while the the remainder does not exceed \ref{4-12-q} with $\gamma=1$;

\smallskip
\noindent
(ii) In particular, the remainder is $O(\mu^{-1/2}h^{-3})$ as
\begin{phantomequation}\label{4-31}\end{phantomequation}
\begin{equation}
\mu \le {\bar\mu}_q \Def h^{-{\frac {q+1}{q+4}}}|\log h|^{-K}
\tag*{$(4.31)_q$}
\label{4-31-q}
\end{equation}
under condition \ref{0-9-q} with $q\ge 1$; in particular ${\bar\mu}_2=\mu^{-1/2}|\log h|^{-K}$ and ${\bar\mu}_3= \mu^{-4/7}|\log h|^{-K}$.
\end{corollary}

I refer to the case of \ref{4-31-q} as a \emph{moderate magnetic field\/}. While I can get rid off logarithmic factors, there is no point to do it right now.

To prove theorem \ref{thm-0-3} one needs to estimate properly expressions (\ref{4-28}) and (\ref{4-29}).

\begin{proposition}\label{prop-4-13} Let conditions $(\ref{0-2})$ and \ref{0-9-q} be fulfilled and $\mu\le h^{\delta-1}$. Let $\psi $ be a nice function supported in $\{|x_1|\le\gamma\}$ with $\mu^{\delta-1/2}\le\gamma\le \epsilon$. Then

\smallskip
\noindent
(i) Both expressions $(\ref{4-28})$ and $(\ref{4-29})$ do not exceed \ref{4-12-q}.

\smallskip
\noindent
(ii) In particular as $\gamma=1$ both expressions $(\ref{4-28})$ and $(\ref{4-29})$ do not exceed $C\mu^{-1/2}h^{-3}$ as $\mu \le \mu_q$.
\end{proposition}

\begin{proof}
Proof is standard, based on the same scaling functions and partition as in the proof of proposition \ref{prop-4-4} as $\psi$ is supported in $\cZ^*_\out$. These arguments work in the other zones as well since I do not need uncertainty principle anymore. I leave the easy details to the reader.
\end{proof}

\begin{corollary}\label{cor-4-14}
 Theorem \ref{thm-0-3} is proven as $\mu\le\mu_q$.\
 \end{corollary}

\subsection{Intermediate magnetic field. II}
\label{sect-4-4}

\subsubsection{}\label{sect-4-4-1}
Now I want to prove theorems \ref{thm-0-1}, \ref{thm-0-3} as $\mu\ge {\bar\mu}_q$ defined by (\ref{4-23}), \ref{4-31-q} for $q=0$ and $q\ge 1$ respectively. To do this I need to consider expression (\ref{0-12}) localized in this singular zone because in the regular zone one can always apply (\ref{0-12}) with $Q=I$ even if the expression is not as explicit as in the case $Q=I$.

\begin{proposition}\label{prop-4-15} One can rewrite with the same error as a remainder estimate derived in section \ref{sect-3} expression $(\ref{0-12})$ as
\begin{equation}
\int {\widehat\cE}^\MW_{\bar Q} (x,0)\psi(x)\,dx
\label{4-32}
\end{equation}
with
\begin{equation}
{\widehat\cE}^\MW_{\bar Q} (x,0) \Def  \const \sum_{n,p} \Bigl(\theta (\cA_{pn}){\bar Q}_{pn}\Bigr)\Bigr|_{(x'',\xi'')=\Psi^{-1}(x)} \times 
f_1(x)f_2(x) \mu^2h^{-2}
\label{4-33}
\end{equation}
with exactly the same constant as in $\cE^\MW$.
\end{proposition}

\begin{proof} Proof is rather standard: I consider canonical form of the operator and in this form I apply the method of successive approximations with operator in question $\cA_{p,n}(x,\mu^{-1}hD_x)$ and unperturbed operator
$\cA_{p,n}(y,\mu^{-1}hD_x)$. I leave details repeating those in my multiple papers to the reader.
\end{proof}

\begin{remark}\label{rem-4-16}
Obviously, defining ${\widehat\cE}^\MW_{\bar Q}$ by (\ref{4-33}) but with $\cA_{pn}$ replaced by $\cA^0_{pn}$ one would get exactly 
\begin{multline}
\cE^\MW_{\bar Q}(x,0)\Def \\\const \sum_{n,p} \theta\Bigl(V(x)- (2n+1)\mu hf_2(x)- (2p+1)\mu h f_1(x)\Bigr) f_1(x)f_2(x){\bar Q}_{pn}\mu^2 h^{-2}
\label{4-34}
\end{multline}
where $\cA^0_{pn}$ is obtained from $\cA_{pn}$ by replacing perturbation terms $\cB_{pn}$ by $0$.
\end{remark}

While one can take $Q=I$ and then ${\bar Q}_{np}$ will be the diagonal elements of $\psi$ transformed according to (\ref{2-55}) the  expression (\ref{4-33}) is not very explicit either because of presence of $\cB_{pn}$ and the similar terms in (\ref{2-55}). Because of this I want to take $Q$ supported in as small zone as possible; the only restriction so far is that $Q\equiv I$ in the singular zone
\begin{multline}
\Omega_\sing \Def \bigl\{ \ell^* \le {\bar \ell}_1=C(r \mu h |\log h|)^{1/2}, \quad
|p-{\bar p}|\le r^{-1}{\bar\ell}_1 (\mu h)^{-1},\\ |n-{\bar n}|\le C\gamma r^{-1}{\bar\ell}_1^2 (\mu h)^{-1}+1
\Bigr\}.
\label{4-35}
\end{multline}
Therefore I will take $Q$ supported in this zone (with increased $C$).

Furthermore I actually need to consider not (\ref{4-32}) but only its correction with respect to what is given by magnetic Weyl formula; namely
\begin{equation}
\int \Bigl({\widehat\cE}^\MW_{\bar Q} (x,0)- \cE^\MW_{\bar Q} (x,0)\Bigr)\psi(x)\,dx  .
\label{4-36}
\end{equation}

\subsubsection{}
\label{sect-4-4-2} 
Now my goal is to estimate expression (\ref{4-36}). I  do it first under condition \ref{0-8-q} (including formally $q=0$) and in the next subsubsection I improve case $q=0$ in the same manner as proposition \ref{prop-4-5} improves proposition \ref{prop-4-4}.

\begin{proposition}\label{prop-4-17}  Let norms of perturbation operators do not exceed $c\varepsilon$, 
\begin{equation}
\varepsilon\le \mu h
\label{4-37}
\end{equation}
and $Q$ be supported in $\cZ^*_\out$.   Then expression $(\ref{4-36})$ does not exceed
\begin{phantomequation}\label{4-38}\end{phantomequation}
\begin{equation}
\left\{\begin{aligned}
&C\varepsilon \gamma^2(\mu h)^{(1+q)/2}h^{-4} |\log h|^K &&q\ge 2,\\
&C\gamma^2\Bigl(\varepsilon(\mu h) h^{-4} + C\varepsilon ^{1/2} (\mu h)^2 h^{-4} \Bigr)|\log h|^K &&q=1,\\
&C\varepsilon^{1/2}(\mu h )h^{-4} \gamma^2|\log h|^K
 +C\gamma^2(\mu h )^2h^{-4}\qquad \qquad &&q=0 ;
\end{aligned}\right.
\tag*{$(4.38)_q$}
\label{4-38-q}
\end{equation}
\end{proposition}

\begin{proof} Let us introduce a scaling function by (\ref{3-38}). Further I replace $\ell$ by $\min (\ell,r)+\varepsilon^{1/2}$. 
Then in  as $\ell=\ell_{pn}\ge C{\bar\ell}_1$, the contribution of all $\ell_{pn}$ groups to (\ref{4-36}) does not exceed the left-hand expression of (\ref{4-14}) with ($\ell$ instead of ${\bar\ell}_1=C(r \mu h |\log h|)^{1/2}$), multiplied by $C\varepsilon \ell^{-2}$:
\begin{equation}
C\mu^2 h^{-2} \gamma^2 r^{(2-q)_+}\times \Bigl(\ell (r\mu h)^{-1}+1\Bigr) \times \Bigl(\ell ^2\gamma (r\mu h)^{-1}+1\Bigr) \times \ell^{q-2}\varepsilon ;
\label{4-39}
\end{equation}
one can prove it easily by considering such elements and integrating by parts if $\ell\ge c\varepsilon^{1/2}$.
Here the third factor does not exceed $|\log h|^K$ for sure as $\ell^2\le \mu h|\log h|^K$ and $\gamma \le r^2$. 

Then as $q\ge 2$ (\ref{4-39}) sums with respect to $\ell$ to its value as $\ell$ reaches its maximum $\min({\bar\ell}_1,r)$; then it sums with respect to $r$  s to the the first line in \ref{4-38-q}. 

As $q=1$ one gets
\begin{equation*}
C (\mu h)h^{-4}r\gamma^2\Bigl( \varepsilon  +(\mu h) \varepsilon \ell^{-1}\Bigr) |\log h|^K;
\end{equation*}
one can estimate $\ell ^{-1}\varepsilon$ by $\varepsilon^{1/2}$ as $\ell\ge \epsilon^{1/2}$ and the summation results in the second line in \ref{4-38-q}.

The same approach works as $q=0$ as well but results in an unwanted now factor $|\log h|^K$ at $C(\mu h)^2h^{-4}\gamma^2$. However, as $\ell \ge \varepsilon^{1/2} |\log h|^{K_1}$ with the large enough exponent $K_1$ one can estimate the fourth factor in (\ref{4-39}) by $|\log h|^{-K_1}$ and all summations will result in 
$C(\mu h)^2h^{-4}\gamma^2 |\log h|^{K_0-K_1}$ with $K_0$ independent on $K_1$; choosing $K_1$ large enough one gets $C\mu^2h^{-2} \gamma^2$. 

On the other hand, one can estimate contribution of all groups with $\ell \le \varepsilon^{1/2} |\log h|^{K_1}$ does not exceed the product of the first three factors in (\ref{4-39}) with $\ell=\varepsilon^{1/2} |\log h|^K_1$ (since one does not need to use an integration by parts then)  which does not exceed after summation with respect to $r$, $\gamma$ the third line in \ref{4-38-q}.
\end{proof}

I remind that  the ``perturbation'' terms in $(\ref{2-50})$ are
\begin{align*}
&\smashoperator{\sum_{2q+2l+2j\ge 3}^{} }\qquad B^\w_{qjl}\times 
\bigl((2p+1)\mu h\gamma\bigr)^q \, \nu^{2-2q-2j}\hbar ^l+\notag\\
&\smashoperator{\sum_{2k+2q+2m+2s+2l\ge 3}^{}}\qquad B^\w_{kqjmsl}\times \bigl((2n+1)\mu h\bigr)^k
\bigl( (2p+1) \mu h\gamma \bigr)^q \,\mu ^{2-2k-2m-2q-s} h ^s r^{-4q-4m-4s} \nu^{-2j} \hbar^l\Bigr).\notag
\label{2-53}
\end{align*}
and due to (\ref{2-51}) their norms do not exceed $C\mu^{-2}  r \gamma^{-3}$
and thus one can take
\begin{equation}
\varepsilon  =   C\min( \mu^{-2}\gamma^{-3}, \mu^{-1-\delta}\gamma^{-1})
\label{4-40}
\end{equation}
as $|x_1|\asymp \gamma$ since $r\le \min (1, \mu^{1-\delta} \gamma^2)$ in $\cZ^*_\out$.
Then condition (\ref{4-37})  translates into 
\begin{equation}
\gamma \ge {\bar\gamma}_2 \Def C \min\bigl(\mu^{-1}  h^{-1/3},\mu^{\delta-2}h^{-1}\bigr)=\left\{\begin{aligned} 
&\mu^{-1}h^{-1/3}\qquad &&\text{as\ } \mu \le h^{\delta'-2/3},\\
&\mu^{\delta -2}h^{-1}&&\text{as\ } \mu \le h^{\delta'-2/3}.
\end{aligned}\right.
\label{4-41}
\end{equation}

\subsubsection{}
\label{sect-4-4-3}
Consider $q=0$ first. Then plugging (\ref{4-40}) into $(\ref{4-38})_0$ and taking the sum with respect to $\gamma$ one gets the same expression calculated as $\gamma=1$ i.e. $Ch^{-3}|\log h|^K + C\mu^2h^{-2}$ which is $O(\mu^2h^{-2})$ as $\mu \ge h^{-1/2}|\log h|^K$.  Therefore due to due to propositions \ref{prop-4-6},  \ref{prop-4-11},  \ref{prop-4-15} and  \ref{prop-4-17} I arrive to

\begin{corollary}\label{cor-4-18}
Estimate $(\ref{0-3})$ holds as $\mu \ge h^{-1/2}|\log h|^K$.\
\end{corollary}

\begin{remark}\label{rem-4-19} As $h^{\delta-1}\le \mu \le ch^{-1}$ the contribution of zone $\{|x_1|\le h^{\delta'}\}$ to the asymptotics is $O(h^{-4+\delta'})=O(\mu^2h^{-2})$ and one does not need to use \ref{prop-4-6} and \ref{prop-4-11}.
\end{remark}

So estimate (\ref{0-3}) remains to be proven as  $h^{\delta-2/5}\le \mu \le h^{-1/2}|\log h|^K$; I will do it in the next subsubsection.

Consider case $q\ge 1$ now.  Plugging (\ref{4-41}) into (the first term of) \ref{4-38-q} one gets $C(\mu h)^{(1+q)/2}h^{-4}\min (\mu^{-2}\gamma^{-1},\mu ^{2-\delta}\gamma)$; after summation with respect to $\gamma\ge {\bar\gamma}_2$ one the same expression \ref{4-38-q}  calculated at 
\begin{equation}
\gamma={\bar\gamma}_3\Def \max({\bar\gamma}_2,\mu^{\delta-1/2})=
\left\{\begin{aligned} 
&\mu^{-1}h^{-1/3}\qquad &&\text{as\ } \mu \le h^{\delta'-2/3},\\
&\mu^{\delta -1/2}&&\text{as\ } \mu \le h^{\delta'-2/3}.
\end{aligned}\right.
\label{4-42}
\end{equation}
which is 
$C(\mu h)^{(q-1)/2}h^{-8/3)}|\log h|^K=O(\mu ^{-1/2})$ as $\mu \le h^{\delta'-2/3}$ and
$C(\mu h)^{(q-2/2}h^{-5/2}\mu^{-\delta}$ as $\mu \ge h^{\delta'-2/3}$; the latter expression is $O(\mu^{-1/2}h^{-3})$ for $q\ge2$.

So, as $q= 2$ one gets $O(\mu^{-1/2}h^{-3})$ as an estimate of the contribution $\cZ^*_\out\cap \{|x_1|\ge {\bar\gamma}_2\}$ into 
\begin{equation}
|\int \bigl( e (x,x,0)-\cE^\MW (x,0)\bigr)\psi (x)\,dx| 
\label{4-43}
\end{equation}
On the other hand, due  to propositions \ref{prop-4-9},  \ref{prop-4-11} and  \ref{prop-4-13} and their corollaries and also propositions \ref{prop-3-20}, \ref{prop-3-21}, \ref{prop-3-22} contribution of zone $\cZ^*_\out \cap \{|x_1|\le {\bar\gamma}_2\}\cup \cZ'_\out \cup \cZ_\inn$ to (\ref{4-42}) is estimated by \ref{4-12-q} calculated at $\gamma={\bar\gamma}_3$.  
One can see easily that the result is $O(\mu^{-1/2}h^{-3})$ as $\mu \le h^{\delta-2/3}$, $q\ge 2$; as $\mu \ge h^{\delta-2/3}$ one gets 
\begin{equation*}
C\mu^{-1/2}h^{-3}+C(\mu h)^{(2+q)/2}h^{-3}\mu^{\delta+1/2}+
C(\mu h)^{(1+q)/2}h^{-3}\mu^\delta
\end{equation*}
where the third term is far less than the second one.

Then for $q=2$ this result does not exceed $C\mu^{-1/2}h^{-3}$ as $\mu \le h^{\delta' -2/3}$ and statement (i) of the proposition below is proven. On the other hand, as $\mu \ge h^{\delta'-4/5}$ ${\bar\gamma}_2$ is below of the bottom of $\cZ^*_\out$ (see \ref{fig-2}) and therefore statement (ii) is proven:

\begin{proposition}\label{prop-4-20} Let conditions $(\ref{0-2})$ and $(\ref{0-8})_2$ be fulfilled. Then

\smallskip
\noindent
(i) As $\mu \le h^{\delta'-2/3}$ expression $(\ref{4-42})$ does not exceed $C\mu^{-1/2}h^{-3}$ and therefore estimate $(\ref{0-6})$ holds with $\cE_\corr^\MW=0$;

\smallskip
\noindent
(ii) As $h^{\delta'-2/3} \le \mu \le ch^{-1}$ contribution of $\cZ^*_\out$ into the $(\ref{4-43})$ does not exceed $C\mu^{-1/2}h^{-3}$ while contribution of 
$\cZ'_\out\cup\cZ_\inn $ (including the inner core) does not exceed
\begin{equation}
C\mu^{-1/2}h^{-3}+  C\mu^{5/2+\delta}h^{-1}.
\label{4-44}
\end{equation}
\end{proposition}

\begin{remark}\label{rem-4-21}
Plugging (\ref{4-40}) into the second term of $(\ref{4-38})_1$ one gets after summation its value as $\gamma=1$ i.e. $C\mu h^{-2}|\log h|^K$. Because this and other ugly terms and because estimate under condition $(\ref{0-8})_1$ is not a part of my core theorems,   I am no more considering $q=1$, leaving to the reader either to derive some estimate of (\ref{4-43}) or to derive $O(\mu^{-1/2}h^{-3}+\mu^{3/2}h^{-\delta-3/2})$ estimate of some more complicated expression.
\end{remark}

\subsubsection{}
\label{sect-4-4-4} To finish the proof of theorem \ref{thm-0-1} I need to prove

\begin{proposition}\label{prop-4-22} Let condition $(\ref{0-2})$ be fulfilled,
$h^{\delta-2/5}\le \mu \le Ch^{-1}$ and $\psi$ be a nice function supported in
$\{ {\bar\gamma}_2=\mu^{-\delta'}\le |x_1|\le \epsilon\}$ with small enough exponents $\delta>0$, $\delta'>0$. Then
\begin{equation}
|\int \Bigl( {\widehat\cE}(x,0)-\cE^\MW (x,0)\Bigr)\psi (x)\,dx|\le C\mu^2h^{-2}.
\label{4-45}
\end{equation}
\end{proposition}

\begin{proof} One needs to consider only zone $r\ge \mu^{-\delta'}$. Here and in the assumption $\delta'>0$ could be taken arbitrarily small. 

Let me introduce scaling function $\ell_m$ as in the proof of proposition \ref{prop-4-5} and the corresponding partition. Then due to the same analysis as in proposition \ref{prop-4-5} the contribution of all the $m$-regular groups with   to the remainder is less than  $C\mu^{-1/2}h^{-3}+C\mu^{-2-\delta''}h^{-2}$.

So one needs to consider only $m$-singular \sf{($(p,n)$,subelement)} groups. Corresponding subelements are of the size ${\bar\ell}=(\mu h|\log h|^K)^{1/(m+1)}$ and on each such group $|\partial^2 \cA_{pn}|\le C (\mu h|\log h|)^{(m-1)/(m+1)}$.  Consider $\epsilon \mu h$-subpartition. Then for each element of it and for each index $p$
\begin{equation*}
|V (x) - f_1(x)(2p+1)\mu h - f_2(x)(2n+1)\mu h |\ge \epsilon_0 \mu h |n-n(p)|
\end{equation*}
and the same is true for this expression perturbed by $\cB_{pn}=O(\varepsilon)$, $\varepsilon=\mu^{-2+4\delta'}$, and also for each element there exists index ${\bar p}$
\begin{equation*}
|\nabla\bigl(V (x) - f_1(x)(2p+1)\mu h - f_2(x)(2n+1)\mu h\bigr) |\ge \epsilon_0 \mu h |p-{\bar p}|
\end{equation*}
as $n=n(p)$.

Then repeating arguments of the proof of proposition \ref{prop-4-20} one can see easily that the left hand expression of (\ref{4-45}) does not exceed $C\mu^{-2}h^{-2}\bigl(1+\varepsilon \sum_{k\ge 1} (\mu h k)^{-1}\bigr)\le C\mu^{-2}h^{-2} $ since $\varepsilon \ll \mu h $  so perturbation does not violate ellipticity of elliptic $\cA_{pn}$. 
\end{proof}

\begin{corollary}\label{cor-4-23}
Theorem \ref{thm-0-1} is proven completely.
\end{corollary}

\subsection{Strong magnetic field}
\label{sect-4-5}

To prove theorem \ref{thm-0-3} or slightly worse estimate under condition 
$(\ref{0-8})_2$ one needs to improve estimate (\ref{4-44}) of the contribution of $\cZ'_\out\cup\cZ_\inn $ (including the inner core) into remainder as $\mu \ge h^{\delta-2/3}$ and to calculate correction term (\ref{3-127})-(\ref{3-128}) as $\mu \ge h^{-6/7}|\log h|^{-2/7}$.

\subsubsection{}
\label{sect-4-5-1}

To improve estimate (\ref{4-44}) of the contribution of $\cZ'_\out\cup\cZ_\inn $ (including the inner core) into remainder without calculating the correction term one can notice that this expression is given by (\ref{0-12}) with $T=T_0= Ch|\log h|(\rho^2+\ell^2)^{-1}$ matching one for 2-dimensional magnetic Schr\"odinger operator (\ref{2-26}) and therefore the asymptotics with the Weyl expression for (\ref{2-26}) 
which one can rewrite easily as
\begin{equation}
(2\pi)^{-2}\mu h^{-3}
\theta \Bigl(2\tau+V - (2n+1)\mu hf_2 -(2p+1)\mu hf_1\Bigr) f_1 f_2 \sqrt g
\label{4-46}
\end{equation}
where I skipped perturbation $\cB'_n=O(\mu^{-2})$ which under condition $(\ref{0-8}_2)$ leads to the relative error $O(\mu^{-2})$. 

While  each expression (\ref{4-46}) is of magnitude $\mu h^{-3}$ and after integration over $\bigl(\cZ'_\out\cup\cZ_\inn \bigr)$ it acquires factor $\gamma=\mu^{-1/2}r^{1/2}$ and after summation over $n$ it acquires factor $(\mu h)^{-1}$; so I get a required expression
\begin{equation*}
\int  \cE^\MW (x,\tau)\psi (x)\,dx
\end{equation*}
 of magnitude $\mu^{\delta-1/2}h^{-4}$ as $\psi $ is supported in
$\bigl(\cZ'_\out\cup\cZ_\inn \bigr)$, an error has magnitude
\begin{equation*}
\mu^{\delta-1/2}h^{-4}r^{5/2}\times \mu ^{-2}= \mu^{\delta-5/2}h^{-4}=
O(\mu^{-1/2}h^{-3}).
\end{equation*}

Therefore I arrive to (almost) final result:

\begin{theorem}\label{thm-4-24} Let $F$ be of Martinet-Roussarie type and condition $(\ref{0-2})$be fulfilled.  Let $\psi$ be supported in the small vicinity of $\Sigma$. Then 

\smallskip
\noindent
(i) Under condition $(\ref{0-8})_3$ estimate $(\ref{0-6})$ holds;

\smallskip
\noindent
(ii) Under condition $(\ref{0-8})_2$ the left hand expression of estimate $(\ref{0-6})$ does not exceed $C\mu^{-1/2}h^{-3}+Ch^{\delta-5/2}$ with an arbitrarily small exponent $\delta>0$;

\smallskip
\noindent
So far $\cE^\MW_\corr$ is defined by $(\ref{3-127})-(\ref{3-128})$ and is
$O(\mu^{5/4}h^{-3/2}|\log h|^{1/2})$; a more explicit representation is given by $(\ref{4-51})$.
\end{theorem}

\subsubsection{}
\label{sect-4-5-2}

First of all note that formula (\ref{IRO6-3-57}) \cite{IRO6} one can rewrite in the form which does not require $g^{jk}|_\Sigma =\delta_{jk}$:

\begin{lemma}\label{lem-4-25} One can rewrite formula $(\ref{IRO6-3-57})$ {\rm\cite{IRO6}} as
\begin{equation}
\cE^\MW_{\corr, d=2}\equiv(2\pi )^{-3/2}h^{-1}\hbar^{1/2}\kappa^{-1/2} V^{(\nu-1)/4\nu} \phi ^{1/2\nu} 
G\Bigl({\frac {S_0V^{(\nu+1)/(2\nu)}\phi^{-1/\nu}}{2\pi \hbar}}\Bigr)\sqrt {g'}\Bigr|_\Sigma
\label{4-47}
\end{equation}
 where $\hbar= \mu^{1/\nu}h$, $g'= g_{22}=g^{11}g$ so $\sqrt {g'}\, dx_2$ is a Riemannian density on $\Sigma$,
 \begin{equation}
 \phi =\bigl(f \cdot \dist (x,\Sigma)^{1-\nu}\bigr)\bigr|_\Sigma,
 \label{4-48}
 \end{equation}
$\dist (.,.)$ is calculated according to the Riemannian metrics $(g_{jk})$ and function $G$ is given by $(\ref{IRO6-3-53})$ {\rm\cite{IRO6}}. 
\end{lemma}

Really, in this formula $V$, $\phi$ and $\sqrt{g'}\,dx_2$ are invariant with respect to change of the coordinates while multiplication of operator by $\omega^2$ is equivalent to substitution $g^{jk}\mapsto \omega^2g^{k}$ and $V\mapsto\omega^2V$ which leads to  $f\mapsto\omega^2f$, $g'\mapsto \omega^{-2}g'$, $\dist (x,\Sigma)\mapsto \omega^{-2}\dist (x,\Sigma)$  (up to a factor equal 1 at $\Sigma$) and $\phi \mapsto \omega^{\nu+1}\phi$; so $V^{(\nu+1)/(2\nu)}\phi^{-1/\nu}$ and $V^{(\nu-1)/4\nu} \phi ^{1/2\nu}\sqrt{g'}$ are invariants.

As $\nu=2$ and $\mu \le Ch^{-1}$ (\ref{3-47})-(\ref{3-48})  become
\begin{equation}
\cE^\MW_{\corr, d=2}\equiv(2\pi )^{-3/2}h^{-1}\hbar^{1/2}\kappa^{-1/2} V^{1/8} \phi ^{1/4} 
G\Bigl({\frac {S_0V^{3/4}\phi^{-1/2}}{2\pi \hbar}}\Bigr)\sqrt {g'}\Bigr|_\Sigma
\label{4-49}
\end{equation}
with an error $O(h^{-1}\hbar)=O(\mu^{-1/2}h^{-1})$ and
 \begin{equation}
 \phi =\|\nabla f\|=\bigl(\sum_{j,k}g^{jk}\partial_j f\partial_kf\bigr)^{1/2}\bigr|_\Sigma,
 \label{4-50}
 \end{equation}
 where norm $\|.\|$ corresponds to the Riemannian metrics $(g_{jk})$.

Similarly, repeating arguments of \cite{IRO6} leading to calculation of the correction term (minimal modifications I leave to the reader) one can calculate 
$(\ref{3-127})-(\ref{3-128})$ modulo $O(\mu^{-1/2}h^{-3})$ deriving the following modification of (\ref{4-49})-(\ref{4-50}):
\begin{multline}
\cE^\MW_\corr(x')=\\
(2\pi)^{-{\frac 3 2}} \mu h^{-2}\hbar^{\frac 1 2}\kappa ^{-{\frac 1 2}}
\sum_{n\in \bZ^+}
\bigl(V-(2n+1)f_2\mu h\bigr)^{1/8} \phi^{1/4}G\bigl({\frac {S_{0}
\bigl(V-(2n+1)f_2\mu h\bigr)^{3/4}\phi^{-1/2}}{2\pi\hbar}}\bigr)f_2\sqrt{g'} 
\Bigr|_\Sigma\label{4-51}
\end{multline}
where $g'=g^{11}g$ and $\sqrt{g'}\,dx'$ is a Riemannian density on $\Sigma$,
$\phi$ is now defined by (\ref{4-50}) with derivatives taken only along $\bK_1$,
$f=f_1$ 
\begin{equation}
 \phi =\| \nabla_{\bK_1} f_1|_{\bK_1} \|
 \label{4-52}
 \end{equation}
where $\|.\|$ corresponds to the Riemannian metrics $(g_{jk})$ restricted to $\bK_1$. One can see easily that only zone $\{\dist (x,\Lambda)\asymp 1\}$ contributes after integration since otherwise\newline 
$|\nabla_\Sigma \bigl(\bigl(V-(2n+1)f_2\mu h\bigr) \phi^{-2/3}\bigr)|$ is disjoint from $0$.

\begin{theorem}\label{thm-4-26} Statement of theorem \ref{thm-4-24}
holds with $\cE^\MW_\corr=O(\mu^{5/4}h^{-3/2})$ defined by $(\ref{4-51})-(\ref{4-52})$.
\end{theorem}

\appendix

\section{Additional Results}
\label{sect-A}

\subsection{Proof of Theorem \ref{thm-0-5}}

Proof of Theorem \ref{thm-0-5} cannot be obtained by the simple rescaling since $f_2\mu$ and $f_1\mu/|x_1|$ scale differently.
However, arguments of sections \ref{sect-1}, \ref{sect-3},\ref{sect-4} work with little or no modifications:

\subsubsection\label{sect-A-1-1} In section \ref{sect-1} propositions \ref{prop-1-5}--\ref{prop-1-9} do not require condition (\ref{0-2}) while proposition \ref{prop-1-10} holds under condition (\ref{0-10}) instead of (\ref{0-2}) where
${\bar T}= Ch|\log h|/\ell$ with
\begin{equation}
\ell= \epsilon\max\bigl(|x_1|, |V|\bigr)+{\hat\ell},\qquad {\hat\ell}=C\mu h|\log h|,
\label{A-1}
\end{equation}
instead of ${\bar T}=Ch|\log h|$.

Further, propositions \ref{prop-1-11}, \ref{prop-1-12}(i) do not require assumption (\ref{0-2}) and proposition \ref{prop-1-12}(ii) holds with modified ${\bar T}$ under condition (\ref{0-10}). Therefore proposition \ref{prop-1-13} also remains true: one needs to estimate
$|F_{t\to h^{-1}\tau}{\bar\chi}_T(t)\Gamma (u\psi_yQ_y^t)|$ as $T={\bar T}$ and this is done easily as $\ell \ge 3{\hat\ell}$ by rescaling and standard Schr\"odinger operator analysis; as $\ell \ge 3{\hat\ell}$ it can be done easily by more crude approach as well.

Furthermore, propositions \ref{prop-1-14}, \ref{prop-1-15} do not require assumption (\ref{0-2}) and due to the previous modifications arguments of proposition \ref{prop-1-16} lead to the same estimate $O(\mu^{-1/2}h^{-3})$ for the contribution of $\cZ_\out$ to the remainder estimate under condition (\ref{0-10}).

Finally, propositions \ref{prop-1-17}--\ref{prop-1-18} also do not require assumption (\ref{0-2}) and proposition \ref{prop-1-19} holds under condition (\ref{0-10}) instead of (\ref{0-2});
proposition \ref{prop-1-20} also does not require (\ref{0-2}) and propositions \ref{prop-1-21}--\ref{prop-1-22} hold under condition (\ref{0-10}).

\subsubsection\label{sect-A-1-2} Section \ref{sect-2} does not require condition (\ref{0-2}) at all.

\subsubsection\label{sect-A-1-3} In section \ref{sect-3} analysis in the strictly outer zone $\cZ^*_\out$ (subsection \ref{sect-3-1}) leading to propositions \ref{prop-3-2}, \ref{prop-3-4} does not require condition (\ref{0-2}); the analysis  and condition (\ref{0-10}) leads to proposition \ref{prop-3-6} which is the worthy substitution for proposition \ref{prop-3-7}.

In the arguments subsections \ref{sect-3-2} and \ref{sect-3-3} in the absence of condition (\ref{0-2}) one again should pick up ${\bar T}= Ch|\log h|/\ell$ with $\ell$ defined by  (\ref{A-1}); in comparison with the arguments of these subsections the factor $\ell^{-1}$ appears but it is compensated by a factor $(|W|+\mu h)$ which appears due the change of range of factors $\rho^2$ or  $(2n+1)\mu h$. This would add an extra term $C\mu ^{1/2}h^{-2}|\log h|$ to the remainder estimate which is subordinate as
$\mu \le C(h|\log h|)^{-1}$.

 On the other hand, as $\mu \ge (h|\log h|)^{-1}$ I can take $\epsilon\ell$-admissible partition with respect to $x'$ and note that as $W\le \epsilon_1\mu h$ at some partition element, it will be classically forbidden and therefore its contribution to the asymptotics would be $0$;  moreover, with the rescaling method I can in the right-hand expression of estimates and using condition (\ref{0-10}) I can take ${\bar T}=Ch$ there thus getting rid off the logarithmic factor:

\begin{proposition}\label{prop-A-1} The total contribution of zones $\cZ'_\out$
and $\cZ_\inn \setminus \cZ_\per$ to the remainder does not exceed $C\mu^{-1/2}h^{-3}$ under condition $(\ref{0-10})$ while the main part of asymptotics is given by $(\ref{0-12})$ with $T=Ch|\log h|$.
\end{proposition}

Moreover
in the vicinity of $\{W=0\}$ condition (\ref{0-10}) also reads as
\begin{equation}
|\nabla_\Sigma \bigl(  V \phi ^{-3/2} \bigr)|\ge \epsilon_0
\label{A-2}
\end{equation}
which kills periodicity for all $n$ not just for all $n$ but one as it would be the case far from $\{W=0\}$ thus leading to

\begin{proposition}\label{prop-A-2} The total contribution of zone $\cZ_\per$
 to the remainder does not exceed $C\mu^{-1/2}h^{-3}$ under condition $(\ref{0-10})$ while the main part of asymptotics is given by $(\ref{0-12})$ with $T={\bar T}$ modified as above and no correction term is needed.
\end{proposition}

\subsubsection\label{sect-A-1-4}
In section \ref{sect-4} condition (\ref{0-10}) allows me to trade easily $T= Ch|\log h|$ to $T=\epsilon \mu^{-1}$ to $T=\epsilon \mu^{-1}\gamma^{-1}$ to $T=\epsilon \mu \gamma^2$ as each next expression is larger than the previous one thus leading to theorem \ref{thm-0-5} as $\mu \le h^{\delta-1}$ covering everything by a moderate magnetic field approach.

Alternatively as $\mu \ge h^{\delta-1}$ the strong magnetic field approach works well: in $\cZ^*_\out$ one can just skip $O(\mu^{-2}\gamma^{-3})$ perturbation thus replacing ${\widehat\cE}^\MW$ by $\cE^\MW$ as in proposition \ref{prop-4-17} with an error well below $\mu^{-1/2}h^{-3}$; similarly in $\cZ'_\out\cup \cZ_\inn$ one can just skip $O(\mu^{-2})$ perturbation thus getting magnetic Weyl expression again.

\smallskip

I leave details to the reader.

\subsection{Special case}\label{sect-A-2}

Consider  the special case of operator $A$ defined as $A_{I}+A_{II}$ where
$A_I$, $A_{II}$ are operators of type (\ref{0-1}) in variables $x_I=(x_1,x_2)$
and $x_{II}=(x_3,x_4)$ respectively with magnetic intensities $f_1(x_I)= |x_1|$ and $f_2(x_{II})= 1$:
\begin{align}
&A_I= h^2D_1^2 + \bigl(hD_2-\mu x_1^2/2\bigr)^2-1,\label{A-3}\\
&A_{II}=h^2D_3^2 + \bigl(hD_4A-\mu x_3\bigr)^2.\label{A-4}
\end{align}
One can separate variables and prove spectral asymptotics with the remainder estimate $h^{-2}R_1$ where $R_1$ is the remainder estimate for operator $A_{I}+\tau$ for the worst possible $\tau>0$.  According to \cite{IRO6}  $R_1=O(\mu^{-1/2}h^{-1})$.

Then the main part of asymptotics is exactly as in estimate (\ref{0-6}) with the
 correction term is given by $(\ref{4-51})$ 
and this term is $O(h^{-3}\hbar^{1/2})$\,\footnote{\label{foot-33} So one loses factor $\mu h$ in comparison with theorem \ref {thm-4-26} in the estimate of the correction term.}.

The leading part of the correction term is rather irregular: it is
$(2\pi)^{-{\frac 3 2}} h^{-3}\hbar^{1/2}\kappa ^{-{\frac 1 2}} I_n$ with
\begin{equation}
I (\varepsilon,\hbar)\Def
\sum_{n\in \bZ^+}
\bigl(V-(2n+1)f_2\varepsilon \bigr)^{1/8} \phi^{1/4}G\bigl({\frac {S_{0}
\bigl(V-(2n+1)f_2\varepsilon\bigr)^{3/4}\phi^{-1/2}}{2\pi\hbar}}\bigr)f_2\sqrt{g'} 
\Bigr|_\Sigma\varepsilon
\label{A-5}
\end{equation}
with $\hbar=\mu^{1/2}h\ll \varepsilon =\mu h\ll 1$.

Obviously that $I_n=O(1)$ and but it is not clear that no better estimate like
$I_n=O(\varepsilon^\sigma)$ is possible; one can see easily that
$\sup_{\varepsilon' \asymp \varepsilon , \hbar'\asymp \hbar} |I(\varepsilon',\hbar')|\ge c^{-1} \varepsilon$.

I hope that some readers will be able make a numerical experiments.

\subsection{About term \texorpdfstring{$C \mu^2h^{-2} $}{%
C\textmu\texttwosuperior /h\texttwosuperior\ } in the remainder estimate}\label{sect-A-3}

Is term $O(\mu^{-2}h^{-2})$ really needed in the general case?
The answer most likely is positive. Consider operator $A$ defined as $A_I+A_{II}$ where 
\begin{equation}
A_I= h^2D_1^2 + \bigl(hD_2-\mu x_1^2/2\bigr)^2-1-kx_1
\label{A-6}
\end{equation}
and $A_{II}$ is given by (\ref{A-4}). Then 
\begin{equation}
W_{pn}\Def V-(2n+1)f_2\mu h-(2p+1)f_1\mu h=1+kx_1 -(2n+1)\mu h-(2p+1)|x_1|\mu h
\label{A-7}
\end{equation}
and as $k=(2p+1)\mu h$, $1=(2n+1)\mu h$ then $\{W_{pn}=0\}=\{x_1>0\}$.

While it does not mean that $\{\cA_{pn}=0\}=\{x_1>0\}$ because of perturbation, I believe that perturbing slightly $A_I$ one get achieve the latter result. I leave it to the reader.

\bibliographystyle{alpha}

\providecommand{\bysame}{\leavevmode\hbox to3em{\hrulefill}\thinspace}

\vglue .06truein

\hfill\hfill {\sl   \today \/}

\vglue .06truein

\begin{tabular}{rrl}
&{\hskip 220 pt} &Department of Mathematics,\cr
&&University of Toronto,\cr
&&40, St.George Str.,\cr
&&Toronto, Ontario M5S 2E4\cr
&&Canada\cr
&&ivrii@math.toronto.edu\cr
&&Fax: (416)978-4107\cr
\end{tabular}

\end{document}